\documentclass[11pt,
]{article}

\usepackage{amssymb,amsfonts,amstext,amsmath,amsthm,
latexsym,mathrsfs,amsbsy}
\usepackage{bbold}

\usepackage{marginnote}

\usepackage{mparhack}
\usepackage{makeidx}
\usepackage[a]{esvect}

\input{38e.sty}

\begin{document}

\title
{Non-uniformizable sets with  
countable cross-sections 
on a given level of the projective hierarchy\thanks
{Vladimir Kanovei's work was supported in part 
by RFBR grant 17-01-00705. 
Vladimir Kanovei is also grateful to the 
Ervin Schroedinger Institute  
%(Vienna, Austria) 
for their support 
%of his research  
during his visit 
%to Vienna 
in December 2016.
Vassily Lyubetsky's work was supported in part 
by RSF grant 14-50-00150.}}

\author 
%, Lyubetsky]
{Vladimir Kanovei\thanks{IITP RAS and RTU (MIIT),
 \ {\tt kanovei@googlemail.com} --- 
contact author. 
}  
\and
Vassily Lyubetsky\thanks{IITP RAS and Moscow University,
\ {\tt lyubetsk@iitp.ru} 
}}

\date 
{\today}

\maketitle

%$\,$\hfill UDC 510.225 and 510.223

\begin{abstract}
We present a model of set theory, in which, for a given $n\ge2$, 
there exists a non-ROD-uniformizable planar $\ip1n$ set, whose 
all vertical cross-sections are countable 
sets (and in fact Vitali classes), 
while all planar $\fs1n$ sets with countable  
vertical cross-sections are $\fd1{n+1}$-uniformizable. 
Thus it is true in this model, that the ROD-uniformization  
principle for sets with countable cross-sections  
fails for a given projective level but holds on  
all lower levels.

%{\bf Keywords:} uniformization, forcing, Vitali class.
\end{abstract}

{\scriptsize
\def\contentsname{\normalsize Contents}
\tableofcontents
}

\parf{Introduction}
\las{ko}

{\em The uniformization problem\/} 
was introduced into descriptive set theory 
by Luzin in a short note~\cite{lhad1} 
and in a more detailed paper~\cite{lhad2}.\snos 
{These notes were not included in Volume~II 
of Luzin's Collected works~\cite{L2L}.   
However its main elements were considered, 
partially translated, and analyzed in detail 
by V.\,A.~Uspensky in~\cite{uspen}. 
Luzin provides in~\cite{lhad1} a quote 
from a Hadamard's letter     
in the well-known ``Five letters''~\cite{cinq}, 
which can be understood so that Hadamard makes 
distinction between the pure Zermelo-style choice 
and choice of elements in non-empty sets by means of 
a concrete effectively defined function. 
This gave Luzin an occasion to connect the uniformization 
problem with the name of Hadamard in the titles 
of papers \cite{lhad1,lhad2}. 
Uspensky argues in~\cite[Section 4]{uspen} that 
the role of Hadamard is definitely exaggerated here, 
while the priority  with regard to the uniformization 
problem and related notions belongs to Luzin himself.} 
A planar set~$Q$ of the real number plane 
$\ves\ti\ves$ is called 
\textit{uniform} (or \textit{single-valued}),
if it intersects every vertical straight line in 
at most one point. 
If $Q\sq P\sq \ves\ti\ves$,
the set~$Q$ is uniform, and its projection onto the 
first axis is equal to the same projection of the 
set $P$, then by Luzin the set $Q$
\textit{uniformizes} the set $P$. 
Saying it differently,
to uniformize a given planar set $P$ means to choose 
a point~$q_x$
in every non-empty vertical cross-section $P_x$ 
of~$P$,
and then gather all chosen points~$q_x$, 
or more precizely, all pairs of the form 
$\langle x,q_x\rangle$,
in a common uniformizing set $Q\subseteq P$.
By Luzin, \textit{the uniformization problem} 
consists in the question 
\textit{is it possible or not to define a point set~$E$ 
for which we cannot name any uniformizing set $E'$}. 
(The translation is quoted from \cite[p.~120]{uspen},   
the italic text by Luzin and Uspensky.)

In modern set theoretic terminology, there exist exact 
definitions for such notions of the ``na\"\i ve'' set 
theory as `to define', `to name', `to give an effective 
construction', and the like. 
The largest class of effectively defined sets   
is the class~$\operatorname{ROD}$ 
of real-ordinal definable sets. 
It consists of all sets definable by a formula with 
real numbers and ordinals as parameters of the definition.
The class~$\operatorname{ROD}$ contains the subclass 
$\operatorname{OD}$ 
of all \textit{ordinal-definable} sets, that is, 
sets definable by a formula with 
ordinals (but not reals) as parameters.

There are more special subclasses of~$\operatorname{ROD}$
and~$\operatorname{OD}$, that is, 
\textit{projective classes} 
$\fs1n$, $\fp1n$, and $\fd1n=\fs1n\cap\fp1n$ and resp.\ 
\textit{effective projective classes}
$\is1n$, $\ip1n$, and $\id1n=\is1n\cap\ip1n$; here $n\ge1$. 
See \cite{mDST} in detail, as well as 
\cite{kl1}, \cite{kl10}, \cite{kl16}, \cite{skmlD}, 
\cite{kanB}, 
on projective hierarchy. 
Recall that 
$\boldsymbol\Delta_1^1=$ Borel sets,
$\boldsymbol\Sigma_1^1=$ Suslin, or A-sets,
$\boldsymbol\Pi_1^1=$ co-Suslin, or CA-sets, 
at the level $n=1$.

The following is considered as the most important 
uniformization theorem in classical 
descriptive set theory.

\bte
[Novikov -- Kondo -- Addison]
\lam{nka}
If\/ $P$ is a planar set in one of the classes\/ 
$\fp11$, $\ip11$, $\fs12$, $\is12$, 
then it can be uniformized by a set in the same class.
\ete

The key ingredient here was the method 
of effective choice of a point in a non-empty 
$\boldsymbol\Pi_1^1$ set by P.\,S.~Novikov,
introduced in~\cite{luno}.
On the base of this method, Kondo \cite{kond} obtained 
the result for $\fp11$. 
Addison \cite{add1,add2} transfered it to the effective  
class $\ip11$. 
The results for classes $\fs12$, $\is12$ are obtained 
by an elementary argument. 
On these and other theorems on uniformization and 
related questions see references above, as well as 
\cite{wood,p12ress,hashi,kash,fis1,fis3} 
with respect to modern studies, and also in the  
introductory section of our paper \cite{kl30}.

As for $\fp12$ and higher projective classes,  
similar uniformization theorems are not available  
since there exist models of set theory in which   
this or another $\fp12$ set  $P$   
is not uniformizable not only by a projective 
(of any class) set, 
but in general by a \rod\ set. 
The first such a model was defined by Levy in 
\cite[Theorem 3]{levy}, where the counter-example 
required is a pnanar $\ip12$ set 
$ 
P=\ens{\ang{x,y}\in\ves^2}{y\nin\rL[x]}
$, 
which is not uniformizable by a ROD set in the model. 
Recall that $\rL[x]$ contains all sets G\"odel  
constructible \rit{relative to} $x$. 

Note that every vertical cross-section 
$
P_x%=\ens{y}{\ang{x,y}\in P}
=\dR\bez\rL[x] 
$ 
of the set $P$ is either empty
(provided $\ves\sq\rL[x]$),  
or else uncountable set, so that it can never  
non-empty finite or countable. 
(Moreover, all cross-sections $P_x$ can be even 
co-countable, as \eg\ in the Solovay model \cite{sol}.)
The problem of the existence of non-uniformizable $\fp12$ 
sets with \rit{countable} vertical cross-sections 
was solved in \cite{kl28} by a model containing such a set. 
Then a more precise result was obtained:

\bte
[proved in \cite{kl30}, equal to the case $\nn=2$ 
in the next Theorem~\ref{Tun}]
\lam{kl30t}%
There exists a model of\/ \zfc, in which it is true 
that there is a planar\/ $\ip12$ set\/ $W\sq\dR^2$  
whose all non-empty vertical cross-sections\/ $W_x$ 
are Vitali classes\snos
{A \rit{Vitali class} in $\dR$ 
is any set of the form $x+\mathbb Q$, 
that is, a shift of the set $\mathbb Q$ of rationals.}, 
which is not uniformizable by a\/ $\rod$ set. 
\ete

The proof involves a generic extension of 
the constructible universe $\rL$ by a forcing 
defined as an uncountable product of an invariant 
version of the Jensen 
minimal forcing \cite{jenmin}.  
\index{forcing!Jensen forcing}%
%\index{zzJ@$\jf$}%
(See also 28A in \cite{jechmill} on the Jensen forcing.) 
Some other results obtained by this method include 
a countable $\ip12$ set containing no definable  
elements \cite{kl27}, 
a Vitali class with the same properties \cite{kl22},
and an $\ip12$ Groszek -- Laver pair of  
Vitali classes. 
See \cite[2.6]{kl30} on the interest in Vitali classes 
in the context of these results.

\parf{The main results}
\las{gr}

In continuation of this research line, 
we prove here the following theorem.

\bte
\lam{Tun}
Let\/ $\nn\ge3$. 
There is a model of the\/ $\zfc$ set theory, 
in which the following is true$:$
\ben
\renu
\itlb{Tun1}% 
there exists a\/ $\ip1\nn$ set\/ $P\sq\ves\ti\ves$, 
such that all cross-sections\/ 
$P_x=\ens{y}{\ang{x,y}\in P}$ 
are Vitali classes, and\/ $P$ is not uniformizable 
by a ROD set$;$

\itlb{Tun2}% 
if\/ $p\in\ves$ then every\/ $\is1\nn(p)$ set\/ 
$P'\sq\ves\ti\ves$, whose all vertical cross-sections 
are countable sets, 
is uniformizable by a\/ $\id1{\nn+1}(p)$ set, 
hence, by a ROD set. 
\een
\ete

Following 
the modern style in descriptive set theory based 
on certain technical advantages,
we shall consider the Cantor discontinuum $\dn$
\kmar{fnote eo} 
with a special equivalence relation\snos
{\label{eo}%
The relation $\Eo$ is defined on $\dn$ 
so that $x\Eo y$ iff the equality $x(n)=y(n)$ holds  
for all but finite indices $n$. 
If $X,Y\sq\dn$ then $X\eeo Y$ means that every  
element $a\in X$ is $\Eo$ equivalent to some  
$b\in Y$, and vice versa. 
See more on this in \cite{kl8,kl10,kanB} or elsewhere.} 
$\Eo$,
instead of the real line $\ves$ with the 
Vitali equivalence relation, 
in the substantial part of the proof.
\index{zzequivE0@$\eeo$}% 
\index{zzE0@$\Eo$}% 
\index{equivalence relation!E0@$\Eo$}% 
\index{equivalence relation!equivE0@$\eeo$}% 
\index{equivalence relation!class $\eko x$}% 
Thus the following theorem will be proved:

\bte
\lam{tun}
Let\/ $\nn\ge3$. 
There is a model of\/ $\zfc$   
in which the following holds$:$
\ben
\renu
\itlb{tun1}% 
there exists a\/ $\ip1\nn$ set\/ $W\sq\dn\ti\dn$, 
such that all cross-sections\/ 
$W_x=\ens{y}{\ang{x,y}\in W}$ 
are\/ $\Eo$-classes, and\/ $W$ is not uniformizable 
by a ROD set$;$

\itlb{tun2}% 
if\/ $p\in\ves$ then every\/ $\is1\nn(p)$ set\/ 
$W'\sq\dn\ti\dn$, whose all sections\/  
$W'_x=\ens{y}{\ang{x,y}\in W'}$ are countable sets, 
is uniformizable by a\/ $\id1{\nn+1}(p)$ set. 
\een
\ete

\noi
{\ubf Theorem~\ref{tun} implies Theorem~\ref{Tun}.}
The transformation of a set $W$ as in \ref{tun}\ref{tun1} 
into a set $P$ as in \ref{Tun}\ref{Tun1}
is carried out by means of elementary topological  
arguments, close to a similar transformation 
in \cite[\S\,17]{kl30}, so we skip this argument. 
The derivation of \ref{Tun}\ref{Tun2} 
from \ref{tun}\ref{tun2} is carried out by means 
of an effective homeomorphism between $\ves$ 
and the co-countable set 
$X=\ens{x\in\dn}{\kaz m\,\sus j\ge m\,(x(j)=0)}$. 
\qed

\parf{Structure of the paper}
\las{str}

The proof of Theorem~\ref{tun}  
is organized as follows. 

The notions related to perfect trees in the set of  
all dyadic strings $\bse$, 
are introduced in sections \ref{tre1},\ref{tre2}. 
We consider a collection $\lt$ of all 
\rit{large} trees ---   
essentially those ones, the relation $\Eo$ on which 
does not admit a Borel transversal. 
Every set $\rP\sq\lt$ closed under truncating trees 
at strings, and \rit{\dd\Eo invariant}, \ie, 
\index{invariant!$\Eo$-invariant}%
\index{e0invariant@$\Eo$-invariant}%
\index{zze0invariant@$\Eo$-invariant}%
invariant relative to that action of finite 
strings which induces the relation $\Eo$ 
(Remark~\ref{dej}), 
is considered 
(Section~\ref{ltf}) 
as a forcing notion adding 
a $P$-generic real $x\in 2^\omega$.
In fact, as $P$ is \dd\Eo invariant, 
an entire \dd\Eo equivalence class 
\index{zzE0x@$\eko x$}% 
\index{equivalence class!E0@$\eko x$}% 
$\eko x=\ens{y\in\dn}{x\Eo y}$ 
of generic reals is adjoined. 

Then in Section \ref{mud} we define the set  
$\mtt$ of all \rit{\mut s}, equal to the 
countable-support product $\lt^{\omi}$. 
We study  \mut s  
(including the behaviour of continuous functions on \mut s) 
in Sections \ref{nf}--\ref{opo}. 

Arguing in the constructible universe $\rL$, we  
define a forcing notion for Theorem \ref{tun} 
in Section \ref{jfor} as the countable-support product  
$\fP=\prod_{\xi<\omi}\yp\xi\sq\mtt$, 
where each factor $\yp{\xi}\sq\lt$ 
has the form of a union  
$\yp\xi=\bigcup_{\xi\le\al<\omi}\yyp\al\xi$,  
where all summands are countable \dd\Eo invariant  
sets $\yyp\al\xi\sq\lt$   
in $\rL$, pre-dense in $\yp{\xi}$. 
\dd\fP generic extensions of $\rL$ will be models 
for Theorem \ref{tun}. 
It turns out that each factor $\yp{\xi}$ adjoins a 
\dd{\yp{\xi}}generic real $x_\xi$, so that the  
whole extension is equal to $\rL[\sis{x_\xi}{\xi<\omi}]$. 
The following is the first key property 
of the forcing notion $\fP$:
\ben
\nenu
\itlb{gr0}%
if $\xi<\omi$ then the set $\yp\xi$ is \dd\Eo invariant.
\een

The next principal issue in the construction of 
forcing notions $\yp\xi$ is similar to the construction  
of Jensen's forcing in \cite{jenmin} 
and in some other cases. 
It consists in the definition of every  
successive  \lap{level} $\yyp\al\xi$ as  
\rit{generic} in some sense over the  
\lap{levels} $\yyp\ga\xi$, $\ga<\al$, already defined. 
This involves a fairly complex construction in   
Sections \ref{gex} -- \ref{pres}, based on the 
splitting technique for perfest trees.
This implies the cardinal preservation (Lemma \ref{wcc}), 
continuous reading of names (Lemma \ref{repd}), 
as well as the following:
\ben
\nenu
\atc
\itlb{gr1}%
for every index $\xi<\omi$,
the set of all \dd{\yp\xi}generic reals in the extension   
is equal to the \dd\Eo class 
$\eko{x_\xi}$ of the generic real $x_\xi$, 
and also is equal to the intersection 
$Y_\xi=\bigcap_{\xi\le\al<\omi}\bigcup_{T\in\yyp\al\xi}[T]$.
\een
Basically we need here only the equality  
$\eko{x_\xi}=Y_\xi$ (Theorem \ref{X=Y}). 
The trans\-form\-ation from a single generic real, 
as in 
Jensen, to a \dd\Eo class of generic reals is 
implied here by the \dd\Eo invariance property  
as in \ref{gr0}.
As a corollary, the definability of the set 
$W=\ens{\ang{\xi,y}}{\xi<\omi\land y\in\eko{x_\xi}}$ 
(the base for a counter-example for \ref{tun}\ref{tun1}) 
in a \dd\fP generic extension 
follows from the definability of the indexed set 
$\sis{\yyp\al\xi}{\xi\le\al<\omi}$ 
in $\rL$ (Section \ref{nounc}).

Following this idea, we proved Theorem \ref{kl30t} in 
\cite{kl30} (= case $\nn=2$ in Theorem \ref{tun}). 
By the way, the    
\rod-non-uniformizability of $W$ follows 
from the \dd\Eo invariance of each component  
of the forcing notion $\fP$ by \ref{gr0},  
both in \cite{kl30} and here.

The main case $\nn\ge3$ in Theorem \ref{tun} differs 
in that it is necessary to prove claim \ref{tun2} in 
the extension, that immediately holds for $\nn=2$   
by Theorem \ref{nka}. 
We get \ref{tun}\ref{tun2} via the following property 
true in \dd\fP generic extensions:
\ben
\atc
\atc
\nenu
\itlb{gr2}%
if $x\in\dn$ and $X\sq\dn$ is a  
countable $\is1\nn(x)$ set then $X\sq\rL[x]$.
\een
This property holds in Cohen and some other generic 
extensions even for $\od(x)$ sets $X$, see \cite{kl31}.
It also holds in  
\dd{\mtt}generic extensions of $\rL$, where it is 
implied by the permutation invariance of the forcing  
notion $\mtt=\lt^{\omi}$ and by a 
very special feature of those extensions, 
namely, 
\ben
\atc
\atc
\atc
\nenu
\itlb{gr3}%
if $x,y\in\dn$ in a \dd{\mtt}generic extension  
$\rL[\sis{x_\xi}{\xi<\omi}]$, and $y\nin\rL[x]$, then 
there exists an ordinal $\xi$ such that $x_\xi\in\rL[y]$ 
but $x\in\rL[\sis{x_\et}{\et\ne\xi}]$ 
\een
(compare to \cite[Theorem 20]{nwf} for the \dd\omi product  
of the Sacks forcing). 
\dd\fP generic extensions satisfy \ref{gr3} as well.  
(Theorem~\ref{fg}, based on the study of continuous  
functions defined on \mut s in Section~\ref{nf}.) 
Yet this does not directly imply \ref{gr2}  since the 
forcing notion $\fP=\prod_\xi\yp\xi$ is not  
permutation-invariant as the components $\yp{\xi}$ 
are pairwise different. 

This leads to the following modification of the forcing 
construction. 
Generally, the construction of $\fP$ can be viewed as 
the choice of a maximal chain in a certain partially 
ordered set $\cP$ of cardinality $\ali$ in $\rL$.
\ben
\atc
\atc
\atc
\atc
\nenu
\itlb{gr5}%
We require that this maximal chain intersects  
all sets dense in $\cP$ which belong to the  
definability class $\fs1{\nn-1}$.
(Theorem \ref{ep}, item \ref{ep2} of which contains  
a property more flexible than this 
straightforward genericity, 
but also more difficult for direct formulation.)
\een
Theorem \ref{ep} evaluates the definability of this 
construction. 
This leads to the definability class $\ip1\nn$ of 
the set $W$ 
(see above) in suitable generic extensions. 
In addition, the forcing notion $\fP$ turns out to  
be enough \lap{generic} in $\mtt$, so that it 
intersects all sets of definability class $\fs1{\nn-1}$, 
dense in $\mtt$ (Lemma \ref{blo}).
This implies a degree of \lap{similarity}
of \dd\fP generic and permutation-invariant  
\dd\mtt generic extensions, up to the \dd\nn th 
level of the projective hierarchy.
And further, by fairly complicated arguments in
Sections \ref{fs}--\ref{zav}  
(which also make use of \ref{gr3}), we obtain \ref{gr2}
in \dd\fP generic extensions, circumwenting the  
above-mentioned problem of the permutation noninvariance 
of $\fP$ and leading to item \ref{tun2} 
of Theorem \ref{tun}.

\parf{Trees and large trees} 
\las{tre1}

Here and in the next section, we reproduce, in brief form,   
some definitions and results from \cite{kl25} related 
to perfect and large trees and their transformations.

{\ubf Strings.} 
$\bse$ is the set of all strings  
\index{zz2om@$\bse$}%
(finite sequences) 
\index{string}%
of numbers $0,1$, 
including  \rit{the empty string} $\La$.
\index{zzLa@$\La$, the empty string}% 
\index{string!$\La$, the empty string}%
If $t\in\bse$ and $i=0,1$ then  
$t\we i$ is the extension of $t$ by $i$ 
as the rightmost term. 
\index{zztiw@$t\we i$}% 
If $s,t\in\bse$ then $s\sq t$ means that the 
string $t$ 
\index{zzstsq@$s\sq t$}% 
\index{zzstsu@$s\su t$}% 
extends $s$ (including the case $s=t$), while 
$s\su t$ means proper extension. 
\index{string!extension}%
The length of $s$ is $\lh s$,  
\index{zzlhs@$\lh s$, length}% 
\index{string!length}%
\index{string!$2^n,$ strings of length $n$}%
and  
$2^n=\ens{s\in\bse}{\lh s=n}$ (strings of length $n$).%

{\ubf Action.} 
Every string $s\in\bse$ {\it acts\/} on $\dn$ so that  
if $x\in\dn$ then
$(s\app x)(k)=x(k)+s(k)\pmod 2$ for $k<\lh s$, and   
$(s\app x)(k)=x(k)$ otherwise.
\index{string!action!$s\app x$}%
\index{string!action!$s\app X$}%
\index{string!action!$s\app t$}%
If $X\sq\dn$ and $s\in\bse$ then let 
$s\app X=\ens{s\app x}{x\in X}$.
\index{zzE0@$\Eo$}% 
\index{equivalence relation!E0@$\Eo$}% 

\bre
\lam{dej}
This action of strings on $\dn$ induces the relation $\Eo$ 
(footnote \ref{eo}), 
so that if $x,y\in\dn$ then $x\Eo y$ 
is equivalent to $y= s\ap x$ for a string  
$s\in\bse$.  
\ere

Similarly if $s\in 2^m,$ $t\in2^n,$ $m\le n$ then define   
a string $s\app t\in 2^n$ so that  
$(s\app t)(k)=t(k)+s(k)\pmod 2$ for $k<m$, and  
$(s\app t)(k)=t(k)$ for $m\le k<n$.
But if $m>n$ then let $s\app t=(s\res n)\app t$.
In both cases, $\lh{s\app t}=\lh t$. 

If $T\sq\bse$ then we let $s\app T=\ens{s\app t}{t\in T}$.\vom

{\ubf Trees.} 
A set $T\sq\bse$ is a \rit{tree}, if  
\index{tree}%
for any strings $s\su t$ in $\bse$, $t\in T$ implies $s\in T$.
If $T\sq\bse$ is a tree and $u\in T$, 
then define \rit{a trunkated subtree}  
\index{tree!truncated, $T\ret u$}%
\index{zzTIu@$T\ret u$}%
$T\ret u=\ens{t\in T}{u\sq t\lor t\sq u}$ of $T$. 
Clearly if $\sg\in\bse$ then  
$\sg\app(T\ret u)=(\sg\app T)\ret {\sg\app u}$.
\index{tree!perfect, $\pet$}%
\index{zzPT@$\pet$}%
%$\pu\ne T\sq \bse$. 
\kmar{pet}%

A non-empty tree $T\sq\bse$ is 
\rit{perfect}, in symbol $T\in\pet$,  
if it has no endnodes and no isolated branches. 
In this case, there is a longest string $s=\roo T\in T$ 
satisfying $T=T\ret s$   
({\it the stem\/} of $T$); 
\index{tree!stem, $\roo T$}%
\index{stem, $\roo T$}%
\index{zzstemT@$\roo T$}%
then $s\we 0\in T$ and $s\we 1\in T$.
If  $T\in\pet$ then the set   
$ 
[T]=\ens{a\in\dn}{\kaz n\,(a\res n\in T)} 
$  
\index{set!$[T]$}% 
\index{zzTII@$[T]$}% 
of all \rit{branches of $T$} 
\index{tree!branch}%
is a perfect set in $\dn$. 
%Note that $[S]\cap[T]=\pu$ if and only if $S\cap T$ 
%is finite.
\vom

{\ubf Large trees.} 
A tree $T\in\pet$ is \rit{large}, 
\index{tree!large, $\lt$}%
\index{zzLT@$\lt$}%
\index{string!uikT@$\qq ik[T]$}%
\kmar{qq ik[T]}%
$T\in\lt$, if there exists a system of strings  
$\qq ik=\qq ik[T]\in\bse,$ $k<\om$ and $i=0,1$, such that  
\ben
\nenu
\itlb{qq1}%
$\lh{\qq0k}=\lh{\qq1k}\ge1$ and $\qq 0k(0)=0$, $\qq 1k(0)=1$ 
for all $k$; 

\itlb{qq2}%
$T$ consists of all strings of the form  
$s=r\we \qq{i_0}0\we \qq{i_1}1\we \qq{i_2}2\we 
\dots\we \qq{i_n}n$ 
and their substrings, where $n<\om$, $r=\roo T$, $i_k=0,1$ 
for all $k$. 
\een
% $k\le n$. 
It this case, the set $[T]$ consists of all infinite 
strings of the form \break 
$a=r\we \qq{i_0}0\we \qq{i_1}1\we \qq{i_2}2\we \dots
\we \qq{i_n}n \we\dots\in\dn,$ 
where $i_k=0,1$, $\kaz k$.
We let %$\oin T0=\lh r$ and then by induction
$$
\oin T{n}=\lh{r}+\lh{\qq{i_0}0}+\lh{\qq{i_1}1}+\dots+
\lh{\qq{i_{n-1}}{n-1}} 
$$ 
(independent of the values of $i_k=0,1$). 
In particular, $\oin T{0}=\lh{r}$. 
Thus 
\index{zzsplnT@$\oin Tn$}%
$\oi T=\ens{\oin Tn}{n<\om}\sq\om$ is the set of all  
\index{zzsplT@$\oi T$}%
\rit{splitting levels} of $T$.

\bre
\lam{BD}
If $T\in\lt$ then the set 
$[T]$ is \rit{\dd\Eo nonsmooth}, 
that is, there is no Borel map $f:[T]\to\dn$  
satisfying $x\Eo y\eqv f(x)=f(y)$ for all  
$x,y\in [T]$.\pagebreak[2]
Conversely, every \dd\Eo nonsmooth Borel set $X\sq\dn$ 
contains a subset of the form $[T]$, where $T\in\lt$.
See \cite{conl}, \cite[7.1]{ksz}, \cite[10.9]{kanB} 
on this category of sets.  

We also note that \rit{Silver trees} are exactly those 
trees $T\in\pel$ satisfying $\qq0k[T](j)=\qq1k[T](j)$ 
for all $k$ and $1\le j<\lh{\qq0k[T]}=\lh{\qq1k[T]}$.
\ere

\parf{Splitting} 
\las{tre2}

%\bdf
%\lam{dspl}
\rit{The simple splitting\/} of a tree $T\in\pel$ consists 
\index{splitting}%
of subtrees
$\raw Ti=\req T{r\we i}$, $i=0,1$, where $r=\roo T$,
\index{splitting!zzT->i@$\raw T i$}%
\index{zzT->i@$\raw T i$}%
so that $[\raw Ti]=\ens{x\in[T]}{x(\lh{r})=i}$.
Then $\raw Ti\in\pel$,  
$\roo{\raw Ti}=r\we \qq i0(T)$,  
$\qq jk(\raw Ti)=\qq j{k+1}(T)$ for all $k$ and $j=0,1$, 
and 
$\oi{\raw Ti}= \oi T\bez\ans{\oin T0}$.

Splittings can be iterated. 
We let $\raw T\La=T$ for the empty string $\La$, 
and if
$s\in 2^n$, $s\ne\La$ then we define   
\index{splitting!zzT->s@$\raw T s$}%
\index{zzT->s@$\raw T s$}%
\index{zzu.s@$u[s]$}%
\index{zzu.sT@$u[s,T]$}%
$$
\raw Ts=\raw{\raw{\raw{\raw T{s(0)}}{s(1)}}{s(2)}\dots}{s(n-1)}
\in\lt
\,.
$$

\bpri
\lam{clop}
If $s\in\bse$ then the tree   
$T[s]=\ens{t\in\bse}{s\sq t\lor t\su s}$ 
belongs to $\pel$, 
\index{zzT:s:@$T[s]$}%
\index{tree!zzT:s:@$T[s]$}%
$\roo{T[s]}=s$, and  
$\qq ik(T[s])=\ang i$.
In particular $T[\La]=\bse$ and   
%  and $\req{T[\La]}s=\req{(\bse)}s=T[s]$ for all $s\in\bse$.
\lam{cloL}
%Under \ref{clop}, 
$T[s]=\raw{(\bse)}s=\req{(\bse)}s$ 
for all $s$.
\epri

\ble
\lam{s2u}
Let\/ $T\in\lt$. 
If\/ $s\in\bse$ then\/  $\raw Ts=\req T{u[s]}$, where\/ 
$u[s]=u[s,T]=\roo{\raw Ts}=
\roo T\we\qq{s(0)}0 \we\qq{s(1)}1 
\we \dots\we \qq{s(n-1)}{n-1}\in T$.  
Conversely if\/ $u\in T$ then
there is a string\/ $s=s[u]\in\bse$ such that\/
$\req Tu=\raw T{s}$. 
\ele
\bpf
To prove the converse, we put $s(k)=u(\oin Tk)$ 
for all $k$ such that $\oin Tk<\lh u$.
\epf

\ble
\lam{tran}
Let\/ $R\in\pel$, $n<\om$, $h=\oin Tn$. 
Then$:$
\ben
\renu
\itlb{tran1}% 
if\/ $u,v\in R\cap2^h$ then\/  
$\req Tu=(u\ap v)\ap (\req Tv)\;;$  

\itlb{tran2}% 
if\/ $s,t\in2^n$ then\/
$\raw Rs=\sg\ap(\raw Rt)$, where\/ 
$\sg=u[s,R]\app u[t,R]\;;$

\itlb{tran3}% 
if\/ $u,v\in R\cap 2^j$, $j<\om$, then\/  
$\req Tu=\sg\ap (\req Tv)$ for some\/ $\sg\in\bse$.   
\een
%\qed
\ele
\bpf
%\ref{tran1} is obvious.
To prove \ref{tran2} use  Lemma \ref{s2u}.
To prove \ref{tran3} 
take the least number $h\in\oi T$ with 
$j\le h$. 
There is a unique pair of strings $u',v'\in2^h$ 
satisfying $u\sq u'$, $v\sq v'$. 
Then $\req Tu=\req T{u'}$, $\req Tv=\req T{v'}$, and 
$\req T{u'}=({u'}\ap{v'})\ap (\req T{v'})$.
\epf

{\ubf Refinement.}
If $R,T\in\pel$ and $n\in\om$ then define $R\nq n T$
(\rit{refinement}), 
\index{zzRsubnT@$R\nq n T$}%
\index{zzsubn@$\nq n$}%
if 
%$R\sq T$ and $\oin Rk=\oin Tk$ for all $k<n$;
$\raw Rs\sq\raw Ts$ for all $s\in2^n$;
$R\nq 0 T$ is equivalent to $R\sq T$.
Clearly $R\nq{n+1} T$ implies $R\nq n T$
(and $R\sq T$).
Moreover, if $n\ge1$  
then\/ $R\nq n T$ is equivalent to  
$\roo R=\roo T$,  
$\qq ik[R]=\qq ik[T]$ for all\/ $i=0,1$ and\/ $k<n-1$, 
and\/ $\qq i{n-1}[T]\sq\qq i{n-1}[R]$ for all\/ $i=0,1$.

\ble
\lam{nadd}
If\/ $T\in\pel$, $s_0\in 2^n$, and\/ 
$U\in\pel$, $U\sq \raw T{s_0}$, then  
there is a unique\/ $T'\in\pel$ satisfying\/
$T'\nq{n} T$ and\/ $\raw{T'}{s_0}=U$. 
We have then 
\ben
\renu
\itlb{nadd1}% 
$\raw{T'}s= u[s_0,T]\ap u[s,T]\ap\raw{T'}{s_0}$ 
for all\/ $s\in2^n\;;$

\itlb{nadd2}% 
if\/ $[U]$ is clopen in\/ $[\raw T{s_0}]$
then\/ $[T']$ is clopen in\/ $[T]\;.$
\een
\ele
\bpf
If $s\in2^n$ then 
$\raw Ts= u[s_0,T]\ap u[s,T]\ap\raw T{s_0}$ 
by Lemma \ref{tran}. 
Put $U_s=u[s_0,T]\ap u[s,T]\ap U$ for all $s\in2^n$, 
in particular, $U_{s_0}=U$. 
The tree $T'=\bigcup_{u\in2^n}U_s$ is as required.
\epf

The next lemma is a more complex version 
of \dd{\nq n}refinement.
The proof see Lemma 4.1(iv) in \cite{kl25}.

\ble
\lam{nadd+}
If\/ $T\in\lt$, $s_0,s_1\in 2^n$, and\/ 
$U,V\in\pel$, $U\sq \raw T{s_0\we0}$, $V\sq \raw T{s_1\we1}$, 
%\kmar{сноска eeo}%
and\/ 
$U\eeo V$ 
{\rm(see footnote \ref{eo} on $\eeo$)}, 
then there exists a tree\/ $T'\in\pel$ 
satisfying\/ $T'\nq{n+1} T$ and\/ 
$\raw{T'}{s_0\we0}\sq U$, 
$\raw{T'}{s_1\we1}\sq V$.\qed 
\ele

\ble
\lam{fus}
Let\/
$\dots \nq 4 T_3\nq 3 T_2\nq 2 T_1\nq 1 T_0$ be an 
infinite sequence of trees\/ $\pel$.
Then\/  
$T=\bigcap_nT_n\in\pel$ and\/ $T\nq{n+1}T_n$, 
$\kaz n$.
\ele
\bpf
Note that  
$\oi T=\ens{\oin{T_n}{n}}{n<\om}$; 
this implies both claims.
\epf

%\pagebreak[3]

\parf{Large tree forcings}
\las{ltf}

\bdf
\lam{sitf}
Let a \rit{\sfo} be  
\index{forcing!\sfo}%
\index{forcing!\sfo!regular}%
\kmar{sfo}%
any set $\rP\sq\pel$ satisfying  
\ben
%\snenu
\Aenu
\itlb{sitf2}%
if $u\in T\in\rP$ then $\req Tu\in \rP$,   
or equivalently, if $T\in\rP$ and $s\in\bse$
then $\raw Ts\in \rP$ ;

\itlb{sitf3}%
$\rP$ is \rit{\dd\Eo invariant}, \ie, 
if $T\in\rP$ and $\sg\in\bse$ then $\sg\app T\in \rP$.  
\index{invariant!$\Eo$-invariant}%
\index{e0invariant@$\Eo$-invariant}%
\index{zze0invariant@$\Eo$-invariant}%. 
\een
If in addition $\bse\in\rP$ then $\rP$  
is a \rit{regular} \sfo.
\edf

Any \sfo\ $\rP$ can be considered as a forcing 
notion (a set of forcing conditions), ordered so that   
if $T\sq T'$ then $T$ is a stronger \usl. 
Such a forcing $\rP$ adjoins a real $x\in\dn$. 
That is, if a set $G\sq\rP$ is  
\dd\rP generic over a ground model $M$, then the 
intersection $\bigcap_{T\in G}[T]$ contains a 
unique real $x=x[G]\in \dn$, and this real  
\index{zzxGII@$x[G]$}%
\index{generic real!$x[G]$}%
satisfies $M[G]=M[x[G]]$ and $G=\ens{T\in\rP}{x\in[T]}$. 
Reals $x[G]$ of this form are called   
\dd\rP\rit{generic}. 

\bpri
\lam{xcoh}
The set $\pel$ of all large trees     
is clearly a \sfo.
Another example of a \sfo\ is the countable set  
$\dpo=\ens{T[s]}{s\in\bse}$ 
of all trees $T[s]$ of Example \ref{clop}, \ie\ 
\rit{Cohen's forcing}.
\index{forcing!Cohen, $\dpo$}%
\index{zzPcoh@$\dpo$}%
\kmar{dpo}%
Finally if $\pu\ne\rQ\sq\pel$ then 
$$
\rP
=\ens{\sg\app(\req Tu)}{u\in T\in\rQ\land\sg\in\bse}
=\ens{\sg\app(\raw Ts)}{T\in\rQ\land s,\sg\in\bse}
$$
is a \sfo\ by Lemma 5.4 in \cite{kl30}.
\epri

A tree $T\in\lt$ is an \rit{\dd ncollage} 
\index{collage}%
over a \sfo\ $\rP$, if we have $\raw Tu\in\rP$ for all 
$u\in2^n$. 
Thus a 0-collage is just a tree in $\rP$, 
and every \dd ncollage is an \dd{n+1}collage as well.

\ble
\lam{c2n}
If\/ $T\in\lt$, $\rP$ is a\/ \sfo, 
$u\in2^n$, and\/ $\raw Tu\in\rP$, 
then\/ $T$ is a \dd ncollage over\/ $\rP$.
In particular, under the conditions of Lemma~\ref{nadd}, 
if\/ $U\in\rP$ then the tree\/ $T'$ obtained is a\/ 
\dd ncollage over\/ $\rP$. 
\ele
\bpf
If $v\in2^n$ then $\raw{T}v=\ta\ap\raw{T}u$ for  
a string $\ta\in\bse$ by Lemma~\ref{tran}. 
Thus $\raw{T}v\in\rP$ since $\raw{T}u\in\rP$.
\epf

If $T\in\pel$ and $D\sq\pel$ then 
$X\sqf\bigcup D$ means 
\index{zzsqf@$\sqf$}%
that there is a finite set $D'\sq D$ satisfying 
$T\sq\bigcup D'$, or equivalently, 
$[T]\sq\bigcup_{S\in D'}[S]$.

\bdf
[extensions]
\lam{fm}
Let $P,Q\sq\lt$ be \sfo s. 
The forcing $Q$ is an \rit{extension} of $P$, 
%\index{perfect tree forcing, $\ptf$!refinement, $\dP\ssq\dQ$}%
\index{forcing!extension, $P\rf Q$}%
\index{extension, $P\rf Q$}%
\index{zzPsqQ@$P\rf Q$}%
\kmar{P rf Q}%
in symbol $P\rf Q$, if
\ben
\nenu
\itlb{fm1}% 
$Q$ is dense in $P\cup Q$: 
if $T\in P$ then $\sus S\in Q\,(S\sq T)$;

\itlb{fm2}% 
If $S\in Q$ then $S\sqf\bigcup P$.
\een
If $\cM$ is any set, 
and, in addition to $P\rf Q$, $S\sqf\bigcup D$ holds 
for all $S\in Q$ and all sets $D\in\cM$, 
$D\sq P$, pre-dense in $P$, then define 
$P\rfa\cM Q$, \dd\cM\rit{extension}.  
\kmar{P rfa MQ}%
\index{forcing!Mextension@\dd\cM extension $P\rfa\cM Q$}%
\index{extension!Mextension@\dd\cM extension, $P\rfa\cM Q$}%
\index{zzPsqMQ@$P\rfa\cM Q$}%
\edf

\ble
\label{pqr'}
\ben
\renu
\itlb{pqr'1}%
If\/ $Q\sq Q'$ and\/ 
$S\sqf\bigcup Q$ for all\/ $S\in Q'$  
then\/ $Q\rf Q'\,;$

\itlb{pqr'3}%
if\/ $P\rfa\cM Q\rf R$ 
{\rm(the second relation is $\rf$, not $\rfa\cM$!)} 
then\/ $P\rfa\cM R\;;$
%\\if\/ $\dP\rsa D\dQ\ssq\dR$  then\/ $\dP\rsa D\dR\;;$ 

\vyk{
\itlb{pqr'4}%
если\/ $\cM$ --- транзитивная модель\/ $\zfc'$, 
содержащая\/ $P$,  
$P\rfa\cM Q$, и деревья\/ 
$T,T'\in P$ несовместны в\/ $ P$    
то\/ $T,T'$ несовместны в\/ $ P\cup Q\,;$ 
}

\vyk{
\itlb{pqr'5}%
если\/ $P\rfa\cM Q$ и множество\/ $D\in\cM$, 
$D\sq P$ предплотно\snos
{Предплотность означает, что каждое дерево $T\in P$ 
совместно в $P$ с каким-то $S\in D$, \te\ найдется 
такое дерево $R\in P$, для которого $R\sq T$ и $R\sq S$.} 
в\/ $P$, то\/ $D$ 
предплотно и в\/ $P\cup Q\,;$
}

\itlb{pqr'6}%
if\/ $\sis{P_\al}{\al<\la}$ is an\/  
%\kmar{pqr2}
\dd\rf increasing sequence of\/ \sfo s 
and\/ $0\le\mu<\la$ then 
\vyk{выполнено следующее$:$
\rm
\ben
\itlb{pqr'6a}%
\it
if\/ $\mu>0$ then\/ $Q=\bigcup_{\al<\mu}P_\al\:\rf\:
R=\bigcup_{\mu\le\al<\la}P_\al\,;$ 

\rm
\itlb{pqr'6b}%
\it
если\/ $\mu>0$ и\/ $Q\rfa\cM P_\mu$ то и\/
$Q\rfa\cM R\,;$

\rm
\itlb{pqr'6c}%
\it
}%
the set\/ 
$P_\mu$ is pre-dense\snos
{Pre-density means that every tree $T\in P$ 
is compatible in $P$ with some $S\in D$, \ie\ there is 
a tree $R\in P$ satisfying $R\sq T$ and $R\sq S$.}  
in\/ 
$P=\bigcup_{\al<\la}P_\al\,.$
%\een
\een
\ele

\bpf
\ref{pqr'3}
$P\rf R$ is clear. 
Assume that a set $D\in\cM$, 
$D\sq P$ is pre-dense in $P$, and $S\in R$. 
Then $S\sqf\bigcup Q$ (since $Q\rf R$), thus  
$S\sq T_1\cup\dots\cup T_n$, 
where $T_1,\dots,T_n\in Q$.
Now $T_i\sqf\bigcup D$, $i=1,\dots,n$, 
since $P\rfa\cM Q$.
It follows that $S\sqf\bigcup D$ holds as well.

\vyk{
\ref{pqr'4}
Пусть напротив, $U\in Q$ и $U\sq T\cap T'$. 
Заметим, что $\cM$ содержит все множества вида 
$
\dpt P{F}
=\ens{S\in P}{\kaz T\in F\,(S\sq T\lor[S]\cap[T]=\pu)} 
$, 
где $F\sq P$ конечно, и все эти множества плотны в $P$.
Отсюда $U\sqf\bigcup\dpt P{\ans{T,T'}}$, 
\te\ $U\sq S_1\cup\dots\cup S_n$, 
где $S_1,\dots,S_n\in\dpt P{\ans{T,T'}}$.
Каждое $S_i$ удовлетворяеи $S_i\sq T$ 
или $[S_i]\cap[T]=\pu$. 
Значит, так как $U\sq T$, мы можем удалить все 
деревья $S_i$ с $[S_i]\cap[T]=\pu$. 
Аналогично, удаляем все такие $S_i$, что  
$[S_i]\cap[T']=\pu$. 
Теперь $S_i\sq T\cap T'$ для оставшихся $i$. 
Но $S_i\in P$.
}

\vyk{
\ref{pqr'5}
Если $T\in Q$ то 
%аналогично 
$T\sq T_1\cup\dots\cup T_n$, 
где $T_1,\dots,T_n\in D$. 
Тогда $\req Ts\sq T_i$ для каких-то $s\in T$ и $i$. 
Но $T'=\req Ts\in Q$. 
Итак, $T$ совместно с $T_i\in D$, значит, 
$D$ предплотно в $P\cup Q$. 
}

\ref{pqr'6}
Let $S\in P_\al$. 
If $\al\le\mu$ then $T\in P_\mu$, $T\sq S$ 
holds by \ref{fm}\ref{fm1}. 
If $\mu<\al$ then  
$S\sq T_1\cup\dots\cup T_n$,  
where $T_1,\dots,T_n\in P_\mu$.
Then $\req St\sq T_i$ for some $t\in S$ and $i$. 
But $S'=\req St\in P_\al$.
\epf

\parf{\Mut s}
\las{mud}

Let a \rit{\mut} be any function  
\index{multitree}%
\index{multitree!$\abs\zp$}%
\index{base!$\abs\zp$}%
\index{zzTII@$\abs\zp$}%
$\zp:\abs\zp\to\pel$, where  
$\abs\zp=\dom\zp\sq\omi$ is at most countable 
and every value $\zc\zp\xi$,
%=\zp(\xi)$, 
$\xi\in\abs\zp$, 
is a tree in $\lt$. 
Let $\mlt$ denote the set of all \mut s. 
\index{multitree!$\mlt$}%
\index{zzMT@$\mlt$}%
If $\zp\in\mlt$ then we define a brick 
\index{brick} 
in $2^{\abs\zp}$,
\index{set!$[\zp]$}%
\index{zztII@$[\zp]$}% 
$$
\bay{rclc}
[\zp]
&=& \ens{x\in2^{\abs\zp}}
{\kaz \xi\in\abs\zp\,(x(\xi)\in [\zc\zp\xi])}&=\\[1ex]
&=&
\ens{x\in2^{\abs\zp}}
{\kaz \xi\,\kaz m\,
(x(\xi)\res m\in \zc\zp\xi)}&, 
\eay
$$
naturally identified with the cartesian  
product $\prod_{\xi\in\abs\zp}[\zp(\xi)]$.

If $B\sq\omi$ is at most countable then let   
$\mlt_B=\ens{\zp\in\mlt}{\abs\zp=B}$.

The set $\mlt$ is ordered componentwise: 
\index{multitrees!ordering}%
$\zp\leq\zq$ 
($\zp$ is a \rit{stronger} \mut) 
\index{multitree!stronger}%
whenever $\abs\zq\sq\abs\zp$ and 
$\zc\zp\xi\sq\zc\zq\xi$ for all $\xi\in\abs\zq$.
Thus the ordering of \mut s corresponds to the   
componentwise inclusion.
The weakest (the largest in the sense of $\leq$) 
\usl\ in $\mlt$ 
is the empty \mut\ 
\index{multitree!$\jLa$}%
$\jLa$, satisfying $\abs\jLa=\pu$. 

It takes some effort to get right versions 
of definitions and results of Section 
\ref{tre2} in the context of \mut s.

\bdf
\lam{proj}
If $\zp\in\mlt_B$ and $C\sq B$, then $\zp\res C\in\mlt_C$ 
is the ordinary restriction. 
But if $B\sq C$ then a \mut\ $\zp\ares C\in\mlt_C$ 
\index{projection, $\res$, $\ares$}%
\index{zzresI@$\res$}%
\index{zzresIup@$\ares$}%
is defined by $(\zp\ares C)(\xi)=\zp(\xi)$ for 
$\xi\in B$, and $(\zp\ares C)(\xi)=\bse$ for 
$\xi\in C\bez B$.
\edf

\bdf
\lam{sqdv}
If $\ju$ is a \mut\ and $\zD$ is a set of \mut s, 
then $\ju\sqfv \zD$ 
\kmar{sqfv}%
means that there is a finite set   
$\zD'\sq \zD$ such that   
1) $\abc \jv\sq C=\abc \ju$ 
for all $\jv\in \zD'$, and 
2) $[\ju]\sq\bigcup_{\jv\in \zD'}[\jv\ares C]$ 
(see Definition \ref{proj} on $\ares$). 
If in addition  
\kmar{sqdv}%
\index{zzsqfin@$\sqfv$}%
\index{zzsqfd@$\sqdv$}%
3) $[\jv\ares C]\cap[\jv'\ares C]=\pu$ for all 
$\jv\ne\jv'$ in $\zD'$, 
then we define $\ju\sqdv \zD$. 
\edf

\bdf
\lam{phi}
Let $B\sq\omi$ be finite or countable.
Fix a function $\phi:\om\onto B$ that takes each 
value infinitely many times, so that 
if $\xi\in B$ then the set
$$
\pro\xi=\ens{k}{\phi(k)=\xi}=
\ans{\bk0\xi<\bk1\xi<\bk2\xi<\ldots<\bk l\xi<\ldots}
$$
is infinite. 
Such a function will be called  
\rit{\dd Bcomplete}.
If $m<\om$ then let $\qi m\xi$ be equal  
\kmar{qi m al}%
\index{zznmalpha@$\qi m\xi$}%
to the number of indices $k<m$, $k\in\pro\xi$. 
Then $\sum_{\xi\in B}\qi m\xi=m$, 
and $\qi m\xi>0$ holds  
for all $\xi\in\img m=\ens{\phi(k)}{k<m}$. 

Let $m<\om$ and $\sg\in2^m$. 
If $\xi\in\img m$ then the set $\pro\xi$ 
\kmar{pro al}%
%\index{zzFalpha@$\pro\xi$}%
cuts a substring $\qa\sg\xi\in2^{\qi m\xi}$ of 
\kmar{qa sg al}%
\index{zzSg.al@$\qa\sg\xi$}%
length $\lh{\qa\sg\xi}=\qi m\xi$ off $\sg$, 
defined by $\qa\sg\xi(j)=\sg(\bk j\xi)$ for all 
$j<\qi m\xi$. 
Thus the string $\sg\in2^m$ splits in a system  
of strings $\qa\sg\xi\in2^{\qi m\xi}$ 
($\xi\in\img m$) 
of total length $\sum_{\xi\in\img m}\qi m\xi=m$.

If $\zp\in\mlt_B$ then define 
$\rau\zp\sg\in\mlt_B$ so that 
${\rau\zp\sg}(\xi)=\raw{\zp(\xi)}{\qa\sg\xi}$
for all $\xi\in B$.
In particular, if $\xi\in B\bez \img m$
then
${\rau\zp\sg}(\xi)=\zp(\xi)$, where $m=\lh\sg$, 
because $\lh{\qa\sg\xi}=\qi m\xi=0$ holds 
provided $\xi\nin\img m$.

If $m<\om$ and $\sg,\ta\in2^m$ then define    
$\zaf\sg\ta=
B\bez\ens{\phi(i)}{i<m\land \sg(i)\ne\ta(i)}$.% 
\index{zzDsta@$\zaf\sg\ta$}% 
%Then by definition 
%$\rau\zp\sg\res{\zaf\sg\ta}=\rau\zp\ta\res{\zaf\sg\ta}$.

Let $\zp,\zq\in\mlt_B$. 
Define $\zp\leq_m\zq$, if  
$\zp(\xi)\nq{\qi m\xi}\zq(\xi)$ for all $\xi\in B$. 
This is equivalent to 
$\rau\zp\sg\sq\rau\zq\sg$ for all $\sg\in2^n$.
\edf

\ble
\lam{tu}
Under the conditions of Definition \ref{phi}, 
let\/ $\zp\in\mlt_B$. 
Then$:$
\ben
\renu 
\itlb{tu1}%
if\/ $\sg\in\bse$ then\/ $\rau{\zp}{\sg}\in\mlt_B$ and 
the set\/ 
$[\rau{\zp}{\sg}]$ is clopen in\/ 
$[\zp];$ 

\itlb{tu1+}%
if\/ $m<\om$ and\/ $\sg,\ta\in2^m$ then\/ 
$\rau\zp\sg\res{\zaf\sg\ta}=\rau\zp\ta\res{\zaf\sg\ta}\,;$ 

\itlb{tu2}%
if\/ $x\in[\zp]$, and\/ $U$ is an open  
nbhd of\/ $x$,  
then there exists a string\/ $\sg\in 2^m$ satisfying\/ 
$x\in[\rau{\zp}{\sg}]\sq U\,;$ 

\itlb{tu3}%
if\/ $m<\om$, $\sg\in 2^m$, and\/ 
$\ju\in\mlt_B$, $\ju\leq\rau\zp{\sg}$, then  
there exists a unique \mut\/ $\zq\in\mlt_B$ such that\/
$\zq\leq_{m}\zp$ and\/ $\rau{\zq}{\sg}=\ju$, 
and then if\/ $[\ju]$ is clopen in\/ 
$\rau\zp{\sg}$ then\/ $\zq$ is clopen in\/ $\zp\,;$ 

\itlb{tu4}%
if\/ $\zD$ is a set of \mut s and\/ $\zp\sqfv \zD$,
then there is a string\/ $\sg\in \bse$ and a \mut\/ 
$\zq\in\zD$ such that\/ $\rau{\zp}{\sg}\leq\zq\,.$
\een
\ele
\bpf
\ref{tu1} is clear. 
\ref{tu2} 
We have $\ans{x}=\bigcap_m[\rau\zp{a\res m}]$ for  
a suitable sequence $a\in\dn$.
By the compactness, 
there is $m$ such that $\rau\zp{a\res m}\sq U$. 

\ref{tu3}
If $\xi\in B$ then
$\ju(\xi)\sq \rau\zp\sg(\xi)=\raw{\zp(\xi)}{s}$, 
where $s=\qa\sg\xi$. 
By Lemma \ref{nadd} there is a tree  
$S_\xi\in\lt$ satisfying $S_\xi\nq n\zp(\xi)$, 
where $n=\qi m\xi=\lh s$, and $\raw{S_\xi}s=\ju(\xi)$. 
Let $\zq(\xi)=S_\xi$, $\kaz \xi$.

\ref{tu4}
There is a \mut\ $\zq\in\zD$ such that 
$\abs\zq\sq B=\abs\zp$ and the intersection 
$U=[\zp]\cap [\zq\ares B]$ has a non-empty interior in  
$[\zp]$. 
It remains to refer to \ref{tu2}.
\epf

\ble
\lam{fusm}
Under the conditions of Definition \ref{phi}, let
$\dots\leq_5\zp_4\leq_4\zp_3\leq_3\zp_2\leq_2\zp_1\leq_1\zp_0$ 
be a sequence of \mut s in\/ $\mlt_B$.
Then the \mut\/ $\zp=\bigwedge_n\zp_n$, 
\index{zzwedgenPn@$\bigwedge_n\zp_n$}%
defined by  
$\zp(\xi)=\bigcap_n\zp_n(\xi)$ for all\/ $\xi\in B$, 
belongs to\/ $\mlt_B$ and\/ $\zp\leq_{n+1}\zp_n$ 
for all\/ $n$.
\ele
\bpf
Apply Lemma \ref{fus} componentwise.
\epf

\parf{Continuous maps and reducibility}
\las{nf}

We consider here some details related to  
continuous maps defined on bricks emerged from 
\mut s, similar to some results obtained 
in \cite{nge,nwf} in the context of perfect sets 
and trees.

Let $B\sq\omi$ be countable, 
$\zp\in\mlt_B$, maps $f,g:[\zp]\to\bn$ 
are continuous.
\bit
\item 
$f$   
\rit{is reduced to\/ $C\sq B$ on\/ $[\zp]$},
if $f(x)=f(y)$ holds whenever 
$x,y\in[\zp]$ and $x\res C=y\res C$.

\item 
$f$   
\rit{is reduced to\/ $g$ on\/ $[\zp]$},
if $f(x)=f(y)$ holds whenever  
$x,y\in[\zp]$ and $g(x)=g(y)$. 
%--- 
%тогда аналогично, существует непрерывная 
%$h:\bn\to\bn$, для которой 
%$f(x)=h(g(x))$ для всех $x\in[\zp]$.

\item
$f$ 
\rit{captures\/ $\al\in B$ on\/ $[\zp]$}, 
if the \rit{co-ordinate map} $c_\al(x)=x(\al)$ 
is reduced to $f$, so that $x(\al)=y(\al)$ holds whenever  
$x,y\in[\zp]$ and $f(x)=f(y)$.
%$x(\al)=h(f(x))$ для всех $x\in[\zp]$.
\eit

\ble
\lam{cap}
If\/ $\zp\in\mlt$, $C_0,C_1,\dots\sq B=\abs\zp$, 
$f:[\zp]\to\bn$ is continuous and reduced  
to every\/ $C_k$ on\/ $[\zp]$, 
then\/ $f$ is reduced to\/ 
$\bigcap_kC_k$ on\/ $[\zp]$.
\ele
\bpf
For just two sets, if $C=C_0\cap C_1$ and $x,y\in[\zp]$, 
$x\res C=y\res C$, then, using the product structure, 
find a point $z\in[\zp]$ with 
$z\res C_0=x\res C_0$ and $z\res C_1=y\res C_1$.
Then $f(x)=f(z)=f(y)$. 
The case of finitely many sets follows by  
simple induction. 
As for the general case, we can assume that 
$C_0\qs C_1\qs C_2\qs \ldots$, by the above.
Let $C=\bigcap_kC_k$, 
$x,y\in[\zp]$, $x\res C=y\res C$. 
There is a sequence of points $x_k\in[\zp]$  
satisfying $x_k\res C_k=x\res C_k$ and 
$x_k\res{(B\bez C_k)}=y\res {(B\bez C_k)}$. 
Then immediately $f(x_k)=f(x)$, $\kaz k$.
On the other hand, clearly $x_k\to y$, hence, 
$f(x_k)\to f(y)$ as $f$ is continuous. 
Thus $f(x)=f(y)$.
\epf

\bte
\lam{nwft}
Let\/ $\zp\in\mlt_B$, $B\sq\omi$ is at most countable, 
and\/ $f,g:[\zp]\to\bn$ continuous. 
Then there is a \mut\/ $\zq\in\mlt_B$, 
$\zq\leq\zp$, satisfying 
either\/ {\rm(i)} 
$f$ is reduced to\/ $g$ on\/ $[\zq]$, or\/ 
{\rm(ii)} 
there is an ordinal\/ $\et\in B$ such that\/
$f$ captures\/ $\et$ on\/ $[\zq]$  
while\/ $g$ is reduced to\/ $B\bez\ans\et$ 
on\/ $[\zq]$.  
\ete

The co-ordinate map $c_\et(x)=x(\et)$ is obviously  
not reducible to $B\bez\ans\et$. 
Thus the theorem essentially says that the  
{\ubf non}reducibility of 
$f$ to $g$ is detected via co-ordinate maps.

\bpf
We argue in terms of Definition \ref{phi}. 
The plan is to define a sequence of \mut s  
as in Lemma \ref{fusm}, with some extra properties. 
Let $m<\om$. 
A \mut\ $\zr\in\mlt_B$ is \dd m\rit{good}, 
if $\zr\leq\zp$ and in addition 
\ben
\item[{\rm(1)$f$\hspace*{0.39ex}}]
if $\sg\in2^m$ and $\al=\phi(m)$ then either 
$f$ is reduced to $B\bez\ans\al$ on $[\rau\zr\sg]$, 
or there is no \mut\ 
$\zr'\in\mlt_B$, $\zr'\leq\rau\zr\sg$, such that 
$f$ is reduced to $B\bez\ans\al$ on $[\zr']$;  

\item[{\rm(1)$g$}]
the same for $g$;  

\item[{\rm(2)$f$\hspace*{0.39ex}}]
if $\sg,\ta\in2^m$, then either 
(i) 
$f$ is reduced to 
$\zaf\sg\ta=B\bez\ens{\phi(i)}{i<m\land \sg(i)\ne\ta(i)}$ 
on  
$[\rau\zr\sg]\cup[\rau\zr\ta]$, 
or 
(ii) 
$\imb f{[\rau\zr\sg]}\cap\imb f{[\rau\zr\ta]}=\pu$;  

\item[{\rm(2)$g$}]
the same for $g$.
\een

\ble
\lam{ngel}
Under the conditions of the theorem, if\/ $m<\om$  
and a \mut\/ $\zr\in\mlt_B$, $\zr\leq\zp$, 
is\/ \dd mgood, then there is an\/ 
\dd{m+1}good 
\mut\/ $\jq\in\mlt_B$, $\jq\leq_{m+1}\zr$. 
\ele
\bpf[Lemma]
Consider a string $\sgi\in2^{m+1}$, 
and first define a    
\mut\ $\zq\in\mlt_B$, $\zq\leq_{m+1}\zr$, 
satisfying (1)$f$ relatively to this string only.
Let $\al=\phi(m+1)$. 
If there exists a \mut\ 
$\zr'\in\mlt_B$, $\zr'\leq\rau\zr\sgi$, such that  
$f$ is reduced to $B\bez\ans\al$ on $[\zr']$ then  
let $\ju$ be one. 
If there is no such $\zr'$ then sipmly put 
$\ju=\rau\zr\sgi$. 
By Lemma \ref{tu}\ref{tu3}, there is a \mut\ 
$\zq\in\mlt_B$ such that $\zq\leq_{m+1}\zr$ 
and $\rau\zq\sgi=\ju$. 
Thus the \mut\ $\zq$ satisfies (1)$f$   
with respect to $\sgi$.
Now take $\zq$ as the \lap{new} \mut\ $\zr$, 
consider another string $\sgi\in2^{m+1}$, and 
do the same as above. 
Consider all strings $2^{m+1}$ 
consecutively, with the same procedure.
In the end, this yields a \mut\ $\zq\in\mlt_B$, 
$\zq\leq_{m+1}\zr$, satisfying (1)$f$  
for all strings in $2^{m+1}$.

Now take care of (2)$f$ . 
Let $\sgi,\tai\in2^{m+1}.$ 
Note that if $\sg'(m)=\ta'(m)$ then 
$\zaf{\sg'}{\ta'}=\zaf{\sg'\res m}{\ta'\res m}$, 
so that (2)$f$  relatively to $\sg',\ta'$ 
follows from (2)$f$  relatively to 
${\sg'\res m},{\ta'\res m}$. 
Thus it suffices to consider only pairs in $2^{m+1}$ 
of the form $\sg\we 0,\ta\we 1$, where $\sg,\ta\in2^m$.  
Consider one such a pair  
$\sgi=\sg\we 0$, $\tai=\ta\we 1$, 
and define a    
\mut\ $\jq\in\mlt_B$, $\jq\leq_{m+1}\zq$, 
satisfying (2)$f$  relatively to this pair.
The sets $\cp=\zaf{\sgi}{\tai}$ 
and $C=\zaf{\sg}{\ta}$ satisfy 
$\cp=C\bez\ans{\eto}$, where $\eto=\phi(m)$, 
while \mut s $\rau\zq\sgi$, $\rau\zq\tai$ 
satisfy 
$\rau\zq\sgi\res{\cp}=\rau\zq\tai\res{\cp}$. 

However by (2)$f$  for the pair $\sg,\ta$, 
either $f$ is reduced to $C$ on 
$[\rau\zq\sg]\cup[\rau\zq\ta]$, 
or  
$\imb f{[\rau\zq\sg]}\cap\imb f{[\rau\zq\ta]}=\pu$. 
In the second case, immediately  
$\imb f{[\rau\zq\sgi]}\cap\imb f{[\rau\zq\tai]}=\pu$. 
Thus we can wlog assume that  
{\ubf $f$ is reduced to $C$ on 
$[\rau\zq\sg]\cup[\rau\zq\ta]$}.  

If now $f$ is reduced to $B'=B\bez\ans{\eto}$ 
on $[\rau\zq\sgi]\cup[\rau\zq\tai]$, 
then $f$ is reduced to $C'=C\cap B'$ as well on 
$[\rau\zq\sgi]\cup[\rau\zq\tai]$ by Lemma \ref{cap}, 
as required. 

Thus suppose that   
{\ubf $f$ is {\bf not} reduced to 
$B'$ on $[\rau\zq\sgi]\cup[\rau\zq\tai]$}.
Then there are points 
$x_0\in[\rau\zq\sgi]$, $y_0\in[\rau\zq\tai]$ with
$x_0\res B'=y_0\res B'$ and $f(x_0)\ne f(y_0)$, \ie, 
$f(x_0)(k)=p\ne q=f(y_0)(k)$ for some $k$; 
$\ans{p,q}=\ans{0,1}$. 
As $f$ is continuous, there exist relatively open 
subsets $X\sq [\rau\zq\sgi]$, $Y\sq [\rau\zq\tai]$, 
such that $x_0\in X$, $y_0\in Y$, 
$f(x)(k)=p$ and $f(y)(k)=q$ for all $x\in X$, 
$y\in Y$. 
We assume wlog that there exists a finite set 
$H\sq B$ containing $\eto$, 
and for each $\et\in H$ strings 
$u_\et\in U_\et=\rau\zq\sgi(\et)=\raw{\zq(\et)}{\qa\sgi\et}$, 
$v_\et\in V_\et=\rau\zq\tai(\et)=\raw{\zq(\et)}{\qa\tai\et}$  
%$v_\et\in\jv(\et)$ 
of equal length 
$\lh{u_\et}=\lh{v_\et}=\ell_\et\ge\qi{m+1,}\et$, 
such that 
$\qa\sgi\et\sq u_\et$, 
$\qa\tai\et\sq v_\et$, 
$$
X=\ens{x\in[\zq]}{\kaz\et\in H\,(x(\et)\in [U'_\et])},  
\text{ and }\,
%,\;
Y=\ens{y\in[\zq]}{\kaz\et\in H\,(y(\et)\in [V'_\et])},
$$ 
where $U'_\et=\req{U_\et}{u_\et}$ and 
$V'_\et=\req{V_\et}{v_\et}$ ($\et\in H$) 
are trees in $\lt$. 
Note that $\qa\sgi\et=\qa\tai\et$, 
$u_\et=v_\et$, $U_\et=V_\et$, and $U'_\et=V'_\et$ 
for all $\et\in H$, $\et\ne\eto$, 
as 
$x_0\res{B'}=y_0\res{B'}$. 
This allows us to define a required \mut\ $\jq$ as follows. 

If $\et\in B\bez H$ then simply put $Q_\et=\zq(\et)$.

Let $\et\in H$, $\et\ne\eto$. 
Then $\qa\sgi\et=\qa\tai\et$, and this string 
$s=\qa\sgi\et=\qa\tai\et$ of length 
$\qi m\et=\qi{m+1,}\et$ satisfies 
$\rau\zq\sgi(\et)=\rau\zq\tai(\et)=\raw{\zq(\et)}s$. 
The string $u=u_\et=v_\et$ belongs to 
$\raw{\zq(\et)}s$, so that the subtree  
$W_\et=\req{\raw{\zq(\et)}s}u=U'_\et=V'_\et$ 
belongs to $\lt$ and $W_\et\sq \raw{\zq(\et)}s$. 
By Lemma \ref{nadd} there exists a tree  
$Q_\et\in\lt$ satisfying
$Q_\et\nq{\qi{m+1,}\et}\zq(\et)$ and $\raw{Q_\et}{s}=W_\et$.

Finally let $\et=\eto\in H$.
The strings $\qa\sgi\eto\ne\qa\tai\eto$ 
of equal length $\qi{m+1,}\eto=\qi{m}\eto+1$
are different in this case since 
$\qa\sgi\eto(\qi m\eto)=\sgi(m)=0\ne1=\tai(m)=
\qa\tai\eto(\qi m\eto)$. 
The equalities $\qa\sgi\eto=\qa\tai\eto$, $U_\eto=V_\eto$, 
$U'_\eto=V'_\eto$ generally speaking, fail as well. 
But still 
$u_\eto\in U_\eto=\raw{\zq(\eto)}{\qa\sgi\eto}$, 
$v_\eto\in V_\eto=\raw{\zq(\eto)}{\qa\tai\eto}$, 
$\lh{u_\eto}=\lh{u_\eto}=\ell_\eto$, and 
$U'_\eto=\req{U_\eto}{u_\eto}\sq U_\eto$, 
$V'_\eto=\req{V_\eto}{v_\eto}\sq V_\eto$. 
Thus we have $U'_\eto\eeo V'_\eto$ by Lemma \ref{tran}. 
Therefore by Lemma \ref{nadd+}, 
there exists a tree $Q_\eto\in\lt$ such that  
$Q_\eto\nq{\qi{m+1,}\eto}\zq(\eto)$, 
$\raw{Q_\eto}{\qa\sgi\eto}\sq U'_\eto$, and 
$\raw{Q_\eto}{\qa\tai\eto}\sq V'_\eto$. 

Thus a tree $Q_\et\in\lt$, satisfying
$Q_\et\nq{\qi{m+1,}\et}\zq(\et)$, 
has been defined for all $\et\in B$  
so that if $\et\in H$ 
then $\raw{Q_\et}{\qa\sgi\et}\sq U'_\et$ and  
$\raw{Q_\et}{\qa\tai\et}\sq V'_\et$. 
This allows us to define a required \mut\ $\jq\in\mlt_B$ 
by $\jq(\et)=Q_\et$ for all $\et\in B$. 
Then $\jq\leq_{m+1}\zq$, and by construction  
$\rau\jq\sgi\sq X$ and $\rau\jq\tai\sq Y$, so that 
$\imb f{[\rau\jq\sgi]}\cap\imb f{[\rau\jq\tai]}=\pu$.

To conclude, a \mut\ $\jq\in\mlt_B$, 
$\jq\leq_{m+1}\zq$, satisfying (2)$f$  relative to   
the pair $\sgi,\tai\in2^{m+1}$ considered, 
is defined.
Going over all pairs in $2^{m+1}$ consecutively, 
we obtain a \mut\ $\jq\in\mlt_B$, 
$\jq\leq_{m+1}\zq$, satisfying (2)$f$  with respect to 
all pairs $\sgi,\tai\in2^{m+1}$.

Then proceed with the same procedure for $g$.
\epF{Lemma}

Come back to the proof of the theorem. 
Lemma \ref{ngel} yields  
an infinite sequence 
$\dots\leq_3\zq_2\leq_2\zq_1\leq_1\zq_0=\zp$
of \mut s $\zq_m\in\mlt_B$, such that each $\zq_m$ 
is a \dd{m}good.
The limit \mut\ 
$\zq=\bigwedge_m\zq_m\in\mlt_B$ satisfies  
$\zq\leq_{m+1}\zq_m$ for all $m$ by Lemma \ref{fusm}. 
Therefore $\zq$ is \dd{m}good for every $m$, hence 
we can freely use (1)$f,g$ and (2)$f,g$ in the 
arguments below.\vom

{\it Case 1\/}: 
if $m<\om$; $\sg,\ta\in2^m$; and 
$\imb f{[\rau\zq\sg]}\cap\imb f{[\rau\zq\ta]}=\pu$; 
then 
$\imb g{[\rau\zq\sg]}\cap\imb g{[\rau\zq\ta]}=\pu$.
Prove that $f$ is reduced to $g$ on $[\zq]$ 
in this case, as required by (i) of the theorem.
Let $x,y\in [\zq]$ and $f(x)\ne f(y)$; 
show that $g(x)\ne g(y)$.
Pick $a,b\in\dn$ satisfying 
$\ans x=\bigcap_m[\rau\zq{a\res m}]$ and 
$\ans y=\bigcap_m[\rau\zq{b\res m}]$. 
As $x\ne y$, we have 
$\imb f{[\rau\zq{a\res m}]}\cap
\imb f{[\rau\zq{b\res m}]}=\pu$
for some $m$ 
by the continuity and compactness. 
Then by the Case 1 assumption, 
$\imb g{[\rau\zq{a\res m}]}\cap
\imb g{[\rau\zq{b\res m}]}=\pu$ holds, 
hence $g(x)\ne g(y)$.\vom

{\it Case 2\/}: 
not Case 1, that is, by (2)$g$, 
there is a pair of strings
$\sgi=\sg\we i, \:\tai=\ta\we j\in2^{m+1}$, $m<\om$, 
such that 
$\imb f{[\rau\zq\sgi]}\cap\imb f{[\rau\zq\tai]}=\pu$, 
but $g$ is reduced to $C'=\zaf\sgi\tai$ on 
$[\rau\zq\sgi]\cup[\rau\zq\tai]$.
%(We refer to .)
We assume that $m$ is the least possible for this case.
We are going to prove that the \mut\ $\rau\zq\sg$ 
satisfies (ii) of the theorem with the ordinal 
$\eto=\phi(m)$, that is, (*) $g$ is reduced to 
$B\bez\ans\eto$ on $[\rau\zq\sg]$, and 
(**) $f$ captures $\eto=\phi(m)$ on $[\rau\zq\sg]$. 

To prove (*) note that 
(\ensuremath\dag)
$f$ is reduced to $C=\zaf\sg\ta$ on 
$[\rau\zq\sg]\cup[\rau\zq\ta]$. 
Indeed otherwise  
$\imb f{[\rau\zq\sg]}\cap\imb f{[\rau\zq\ta]}=\pu$ 
by (2)$f$, 
hence $g$ is \rit{not} reduced to $C$ on 
$[\rau\zq\sg]\cup[\rau\zq\ta]$ by the minimality of $m$, 
thus 
$\imb g{[\rau\zq\sg]}\cap\imb g{[\rau\zq\ta]}=\pu$ 
by (2)$g$, then obviously 
$\imb g{[\rau\zq\sgi]}\cap\imb g{[\rau\zq\tai]}=\pu$, 
which contradicts to the assumption that 
$g$ is reduced to $C'$ on $[\rau\zq\sgi]\cup[\rau\zq\tai]$ 
because   
$[\rau\zq\sgi]\res C'=[\rau\zq\tai]\res C'$ 
by Lemma \ref{tu}\ref{tu1+}.

%$\imb f{[\rau\zq\sgi]}\cap\imb f{[\rau\zq\tai]}=\pu$, 

Note that $i\ne j$, as otherwise $C=C'$, and 
a contradiction easily follows. 
Thus let, \eg, 
$\sgi=\sg\we0,\;\tai=\ta\we1$. 
Then $C'=C\bez\ans\eto$, and the \mut\ $\rau\zq\sgi$
witnesses that $g$ is reduced to 
$B\bez\ans\eto$ on $[\rau\zq\sg]$ by (1)$g$.
Thus we have (*).

We further claim that (\ensuremath\ddag) 
\rit{$f$ is not reduced to $B\bez\ans\eto$ on any  
\mut\ $\ju\in\mlt_B$, $\ju\leq\rau\zq\sg$}.
Indeed otherwise $f$ is reduced to $B\bez\ans\eto$ 
on $[\rau\zq\sg]$ by (1)$f$. 
Then $f$ is reduced to $C'=C\cap (B\bez\ans\eto)$    
on $[\rau\zq\sg]$ by Lemma \ref{cap}.
It follows that $f$ is reduced to $C'$    
on the union $W=[\rau\zq\sg]\cup[\rau\zq\ta]$,\snos 
{Let $x,y\in W$ and $x\res C'=y\res C'$. 
As $\rau\zq\sg\res C=\rau\zq\ta\res C$ 
by Lemma \ref{tu}\ref{tu1+}, 
there are $x',y'\in[\rau\zq\sg]$ with  
$x\res C=x'\res C$ and $y\res C=y'\res C$.
But $f$ is reduced to $C$ on $W$ by (\ensuremath\dag). 
Thus $f(x)=f(x')$, $f(y)=f(y')$. 
Finally $f(x')=f(y')$ since $f$ is reduced to $C'$    
on $[\rau\zq\sg]$.}
hence on $W'=[\rau\zq\sgi]\cup[\rau\zq\tai]\sq W$ 
as well.
But this contradicts to 
$\imb f{[\rau\zq\sgi]}\cap\imb f{[\rau\zq\tai]}=\pu$, 
since 
%\mut s $\rau\zq\sgi$ and  $\rau\zq\tai$ satisfy  
$\rau\zq\sgi\res C'=\rau\zq\tai\res C'$ by 
Lemma \ref{tu}\ref{tu1+}. 
As required. 

%In continuation of Case 2, 
We now prove (**), that is, 
{\ubf $f$ captures $\eto$ on $[\rau\zq\sg]$}. 
Let $x,y\in [\rau\zq\sg]$ and $f(x)=f(y)$; 
prove that $x(\eto)=y(\eto)$. 
Indeed    
$\ans x=\bigcap_n[\rau\zq{a\res n}]$ and  
$\ans y=\bigcap_n[\rau\zq{b\res n}]$, 
where $a,b\in\dn$, $\sg\su a$, $\sg\su b$. 
Let $\zaf ab=\bigcap_n\zaf{a\res n}{b\res n}$. 
Then $x\res\zaf ab=y\res\zaf ab$, since 
$\rau\zq{a\res n}\res{\zaf{a\res n}{b\res n}}=
\rau\zq{b\res n}\res{\zaf{a\res n}{b\res n}}$
for all $n$. 
Thus it remains to check that 
$\eto\in\zaf{a\res n}{b\res n}$ for all $n$.

Suppose towards the contrary that 
$\eto=\phi(m)\nin\zaf{a\res n}{b\res n}$ 
for some $n$. 
Note that $n>m$ because $a\res m=b\res m=\sg$. 
However $f$ is reduced to $\zaf{a\res n}{b\res n}$ 
on $[\rau\zq{a\res n}]$
by (2)$f$, since $f(x)=f(y)$. 
Yet we have $\eto\nin\zaf{a\res n}{b\res n}$, 
therefore, $\zaf{a\res n}{b\res n}\sq B\bez\ans\eto$. 
It follows that $f$ 
is reduced to $B\bez\ans\eto$ on $[\rau\zq{a\res n}]$. 
But this contradicts to (\ensuremath\ddag)  above with 
$\ju=\rau\zq{a\res n}$.  

To conclude Case 2, we have checked (*) and (**), 
as required.  
\epf

\vyk{
\bcor
%[ср.\ теорема 20 в \cite{nge}]
\lam{nget}
Пусть\/ $\zp\in\mlt$, $E\sq B=\abs\zp$, 
$f:[\zp]\to\bn$ непрерывна. 
Найдется такое \mutо\/ $\zq\in\mlt_B$, 
$\zq\leq\zp$, 
что либо\/ $f$ сводится к\/ $E$ на\/ $\zq$ либо\/ 
$f$ захватывает некоторый ординал\/ $\al\in B\bez E$ 
на\/ $\zq$. 
\ecor
\bpf
Возьмем любую непрерывную функцию 
$g:[\zp]\to\bn$, удовлетворяющую 
$x\res E=y\res E\eqv g(x)=g(y)$. 
Также следует заметить, что если $\al\in B\bez E$ 
то координатная функция $h(x)=x(\al)$ не сводится 
ни к какому $\zq\in\mlt_B$.
\epf
}

\parf{\Muf s and \smuf s}
\las{muf}
             
Let a \rit{\muf} be any function 
\index{multiforcing}%
\index{multiforcing!small}%
\index{multiforcing!regular}%
\index{multiforcing!$\abs\jpi$}%
\index{base!$\abs\jpi$}%
\index{zzPII@$\abs\jpi$}%
$\jpi$, such that   
$\abs\jpi=\dom\jpi\sq\omi$ and every value 
$\jpi(\xi)$, $\xi\in\abs\jpi$, is a \sfo. 
Thus a \muf\ is a partial \dd\omi sequence of \sfo s.
A \muf\ $\jpi$ is \rit{small}, 
if the base 
%(domain) 
$\abs\jpi$  
and each forcing $\pc\jpi\xi$, $\xi\in\abs\jpi$, 
are at most countable sets, and 
\rit{regular}, if $\bse\in\jpi(\xi)$ for all 
$\xi\in\abs\jpi$.

If $\jpi$ is a \muf\ then let $\mt\jpi$ denote
the set of all \mut s $\zp$ such that
$\abs\zp\sq\abs\jpi$ and $\zc\zp\xi\in\pc\jpi\xi$   
\index{multitrees!$\mt\jpi$}%
\index{multiforcing!$\mt\jpi$}
\index{zzMTP@$\mt\jpi$}%
for all $\xi\in\abs\jpi$.
The set $\mt\jpi$ can be identified with 
the \rit{countable base product\/ 
$\prod_{\xi\in\abs\jpi}\pc\jpi\xi$}.

The next definition introduces a type of sets containing  
\mut s and satisfying some minimal closure conditions.
 
\bdf
\lam{mfd}
Let $\jpi$ be a regular \muf. 
A set $\mP\sq\mt\jpi$ is a \rit{\smuf}, 
\index{submultiforcing}%
if it satisfies the following: 
\ben
\Renu
\itlb{mfd1}%
if $\zp\in\mP$, $\xi\in\abs\zp$, and $T\in\jpi(\xi)$, 
then the \mut\ $\zq$, defined by $\abs\zq=\abs\zp$, 
$\zq(\xi)=T$, and 
$\zq(\et)=\zp(\et)$ for $\et\ne\xi$, 
also belongs to $\mP$; 

\itlb{mfd2}%
if $\zp\in\mP$, $\xi\in\abs\jpi\bez\abs\zp$, and 
$T\in\jpi(\xi)$, 
then the \mut\ $\zq$, defined by 
$\abs\zq=\abs\zp\cup\ans\xi$, 
$\zq(\xi)=T$, and $\zq\res{\abs\zp}=\zp$, 
also belongs to $\mP$; 

\itlb{mfd3}%
if $\zp,\zq\in\mP$ then 
the \mut\ $\zp'=\zp\ares{(\abs\zp\cup\abs\zq)}$, 
defined by 
$\abs{\zp'}=\abs\zp\cup\abs\zq$, 
$\zp'(\xi)=\zp(\xi)$ for $\xi\in\abs\zp$, and  
$\zp'(\xi)=\bse$ for $\xi\in\abs\zq\bez\abs\zp$, 
also belongs to $\mP$.\qed 
\een
\eDf 

\bpri
\lam{mcoh}
Let $\jpi$ be a regular \muf, $B=\abs\jpi$.
Then $\mt\jpi$ is the largest \smuf\ 
in $\mt\jpi$, while the smallest \smuf\ in $\mt\jpi$ 
is the countable set $\mpo B$ 
\kmar{mpo B}%
of all \mut s $\zp\in\mt\jpi$ such that 
$\abs\zp\sq B$ is finite and $\zp(\xi)\in \dpo$ 
(Example \ref{xcoh}) for all $\xi\in\abs\zp$,  
\rit{Cohen's forcing} in $(\dn){}^B$.
\index{forcing!Cohen, $\mpo B$}%
\index{zzScohB@$\mpo B$}%
\epri         

\Mut s $\zp,\zq$ in a \smuf\ $\mP\sq\mtt$ are 
\rit{compatible in\/ $\mP$}, if there is a \mut\ 
\index{multitrees!compatible}%
$\ju\in \mP$ satisfying $\ju\leq\zp$ and 
$\ju\leq\zq$.
A set $\zD\sq\mP$ is:
\bde
\item[{\it dense in $\mP$}\rm,]
\index{multitrees!set!dense}%
\index{set!dense}%
when $\kaz\zp\in\mP\:\sus\zq\in \zD\:(\zq\leq\zp)$;

\item[{\it open dense in $\mP$}\rm,]
if in addition  
\index{multitrees!set!open dense}%
\index{set!open dense}%
$\,$\vphantom{.}\ $\kaz\zp,\zq\in\mP\:
(\zp\leq\zq\in \zD\imp\zp\in \zD)$;

\item[{\it pre-dense in $\mP$}\rm,]
if the set \
\index{multitrees!set!pre-dense}%
\index{set!pre-dense}%
%$\,$\vphantom{.}\hfill 
$\zD^+=\ens{\zp\in\mP}{\sus\zq\in \zD\:(\zp\leq\zq)}$  
\ is dense in $\mP$.
\ede

In the context of Definition \ref{phi}, 
a \mut\ $\zp$ (not necessarily $\zp\in\mP$!) 
is called \rit{an \dd mcollage over\/ $\mP$},  
\index{multitree!collage}%
\index{collage}%
if $\rau\zp u\in\mP$ for all strings $u\in2^m$. 
Thus a 0-collage is any  
\mut\ in $\mP$, while every \dd mcollage is    
an \dd{m+1}collage as well by the closure properties  
in Definition~\ref{mfd}. 

\ble
\lam{mokr}
Let\/ 
$\jpi$ be a \muf, $\mP\sq\mt\jpi$ be a \smuf, 
$\zp\in\mlt_B$. 
% and\/ $\phi:\om\to B$ be a\/ \dd Bcomplete map.
Then, in terms of Definition \ref{phi}, 
the following holds$:$
\ben
\renu
\itlb{mokr1}%
if\/ $\sg\in\bse$ and\/ $\zp\in\mP$ 
then\/ $\rau\zp\sg\in\mP$$;$

\itlb{mokr3}%
if\/  
$\sg\in2^n$ and\/ $\rau \zp\sg\in\mP$, 
then\/ $\zp$ is an \dd ncollage over\/ $\mP$$;$

\itlb{mokr4}%
if\/ $\zp$ is an\/ \dd mcollage over\/ $\mP$, 
and\/ $\zD\sq\mP$ is clopen in\/ $\mP$, 
then there is a \mut\/ $\zq\in\mlt_B$,  
which is an\/ \dd mcollage over\/ $\mP$ and satisfies\/
$\zq\leq_{m}\zp$ and\/ $\rau{\zq}{\sg}\in\zD$ 
for all\/ $\sg\in2^m$$;$ 

\itlb{mokr2}%
if\/ $U\sq [\zp]$ is a nbhd of\/ $x_0\in[\zp]$ 
in\/ $[\zp]$ then there is a \mut\/ $\zq\in\mP$ 
such that\/ $\abs\zq=B$, $x_0\in[\zq]\sq U$, 
and\/ $\zq\leq\zp$. 
\een
\ele
\bpf
\ref{mokr1}
Use property \ref{sitf}\ref{sitf2} of \sfo s with  
the closure properties of Definition \ref{mfd}.
Further, splitting the operation 
$\rau{}\sg$ to components as in  
Definition \ref{phi},  immediately reduces \ref{mokr3}  
to Lemma \ref{c2n}.

\ref{mokr4}
If $\sg\in2^m$ then by Lemma \ref{tu}\ref{tu3} 
there exists a \mut\ $\zq\in\mlt_B$, $\zq\leq_{m}\zp$, 
satisfying $\rau{\zq}{\sg}\in\zD$ 
for this $\sg$.
And $\zq$ is still an \dd mcollage over $\mP$ 
by \ref{mokr3}.
Iterate this procedure, going over all strings in  
$\sg\in2^m$.

\ref{mokr2}
We refer to \ref{mokr1} and Lemma \ref{tu}\ref{tu2}.
\epf

\parf{On subsets with the Baire property}
\las{sb}

This and the next section present two applications  
of Lemma \ref{fusm} to the construction of \mut s 
with certain properties. 
Comparing to Theorem \ref{nwft}, where Lemma \ref{fusm}
was also used in the course of the proof, here 
by necessity 
we'll have to consider intermediate \mut s  
related to some \muf.

\ble
\lam{ssq}
Let\/ $\zp\in\mlt$, $B=\abs\zp$. 
If the set\/ $X\sq[\zp]$ has the Baire property  
inside\/ $[\zp]$ then there is a \mut\/ $\zq\in\mtt_B$ 
such that\/ 
$[\zq]\sq X$ or\/ $[\zq]\sq [\zp]\bez X$. 
\ele
\bpf
Fix a \dd Bcomplete function $\phi:\om\onto B$. 
In our assumptions, $X$ or $[\zp]\bez X$ 
is co-meager on a non-empty clopen $U\sq[\zp]$. 
The cases are symmetric, hence we can assume that 
$X$ is co-meager on $U$. 
Note that $[\rau\zp\sg]\sq U$ for some $\sg\in\bse$ 
by Lemma \ref{tu}\ref{tu2}.
Yet the set $[\rau\zp\sg]$ itself is clopen in  
$[\zp]$, and $X'=X\cap[\rau\zp\sg]$ is co-meager in  
$[\rau\zp\sg]$. 
Thus the task is reduced to the case when the set $X$ 
is co-meager in $[\zp]$, and this will be assumed below. 
In this assumption, we can further suppose that  
$X=\bigcap_nU_n$, where every set 
$U_n\sq[\zp]$ is topologically open and dense 
in $[\zp]$.\vom

{\sl Case 1\/}:
there exists a \mut\ $\zq\in\mtt_B$ such that 
$\zq\leq\zp$ and $[\zq]\cap U_n=\pu$ for some $n$. 
Then $[\zq]\sq [\zp]\bez X$, as required.\vom

{\sl Case 2\/}:
if $\zq\in\mtt_B$ and $\zq\leq\zp$ then 
$[\zq]\cap U_n\ne\pu$ for all $n$. 
Define a regular \muf\ $\jpi$ such that $\abs\jpi=B$  
and if $\xi\in B$ then
$$
\jpi(\xi)=
\ens{s\ap(\raw{\zp(\xi)}t)}{s\in\bse\land t\in\zp(\xi)}
\cup\dpo\;
\text{ (see Definition \ref{xcoh})}.
$$
Consider the \smuf\ 
$\mP=\ens{\zq\in\mt\jpi}{\abs\zp=B}$; 
$\zp\in\mP$.
We claim that for every $m$ the set
$$
\zD_m=\ens{\zq\in\mP}{[\zq]\cap[\zp]=\pu \,\text{ or }\,
\zq\leq\zp\land[\zq]\sq U_m }
$$
is open dense in $\mP$ (in the sense of Section \ref{muf}).
The openness is obvious. 
To prove the density let $\zp'\in\mP$. 
If $[\zp']\not\sq[\zp]$ then $U=[\zp']\bez[\zp]$
is topologically open in $[\zp']$ and non-empty. 
By Lemma \ref{tu}\ref{tu2}, there exists a  
\mut\ $\zq\in\mP$ such that $[\zq]\sq U$, 
\ie, $\zq\leq\zp'$ and $\zq\in \zD_m$.
Thus assume that $\zp'\leq\zp$. 
Then $[\zp']\cap U_m\ne\pu$ by the Case 2 assumption. 
Applying Lemma \ref{tu}\ref{tu2}, we 
find a \mut\ $\zq\in\mP$ satisfying $[\zq]\sq U_m$, 
that is,  $\zq\in \zD_m$. 
The density is proved. 
 
Now Lemma \ref{mokr}\ref{mokr4} implies a  
sequence 
$\dots\leq_4\zp_3\leq_3\zp_2\leq_2\zp_1\leq_1\zp_0\leq\zp$ 
of  \mut s $\zp_m\in\mlt_B$ with    
$\rau{\zp_m}\sg\in\zD_m$ for all $m$ and $\sg\in2^m.$ 
The \mut\ $\zq=\bigwedge_m\zp_m$ (Lemma \ref{fusm}) 
then satisfies $[\zq]\sq U_m$, $\kaz m$, 
hence $[\zq]\sq X$.
\epf

\parf{Separating image from preimage}
\las{opo}

If $x_0\in X\sq\dn$, $f:X\to\dn$ is continuous, 
and $f(x_0)\ne x_0$, then there exists a nbhd  
$U$ of $x_0$ in $X$ whose \dd fimage $\imb fU$ 
does not intersect $U$.
The next theorem is a version of this claim. 
% в контексте \mutьев.
 
\bdf
\lam{prof}
Let $\zp\in\mtt_B$ and $\xi\in B$.  
A continuous map $f:[\zp]\to\dn$ is called  
\index{map!simple}%
\rit{simple on\/ $[\zp]$ for\/ $\xi$,}
if there exists a string 
$\sg\in\bse$ such that $f(x)=\sg\ap x(\xi)$ 
holds for all $x\in[\zp]$. 
\edf    

\bte
\lam{nf'}
Under the conditions of Definition\/ \ref{phi}, let\/ 
$\xi\in B=\abs\jpi$, 
$\mP\sq\mt\jpi$ is a \smuf, $m,n<\om$,
$\zp\in\mtt_B$ is an\/ \dd mcollage over\/ $\mP$, and  
$f:[\zp]\to\dn$ is continuous.
Then$:$
\ben
\renu
\itlb{nf'1}%
if\/ $U\in\lt$ is an\/ \dd ncollage 
over a \sfo\/ $\rP$,
then there exists a \mut\/ $\zp'\in\mtt_B$  
and a tree\/ $U'\in\lt$ such that\/ 
$\zp'\leq_m\zp$, $U'\nq n U$, 
$\zp'$ is an \dd mcollage over\/ $\mP$, 
$U'$ is a \dd ncollage over\/ $\rP$, 
and\/ $[U']\cap\imb f{[\zp']}=\pu\,;$

\itlb{nf'2}%
if\/ $\xi\in B=\abs\jpi$, and if\/ $r\in\bse$ 
then\/ $f$ is not simple for\/ $\xi$ on\/ 
$\rau\zp r$,   
then there is a \mut\/ $\zp'\in\mtt_B$ such that\/ 
$\zp'\leq_m\zp$,  
$\zp'$ is an\/ \dd mcollage over\/ $\mP$, 
and\/ $[\zp'(\xi)]\cap\imb f{[\zp']}=\pu\,.$ 
\een
\ete
\bpf
\ref{nf'1}
To begin with consider a pair of strings $u\in2^m$, 
$s\in 2^n$. 
Let $x_0\in[\rau\zp u]$. 
Pick $y_0\in[\raw Us]$, $y_0\ne f(x_0)$. 
As $f$ is continuous, there exists an open nbhd  
$G\sq[\zp]$ of $x_0$ in $\rau\zp u$ and a string  
$t\in\raw Us$, satisfying  
$t\su y_0$, and $t\not\su x(\xi)$ for all $x\in G$. 
Put $V=\req Ut$. 
Then $V\in\rP$ and $V\sq\raw Us$. 
By Lemma \ref{nadd}, there exists a tree $U'\in\lt$, 
such that $U'\nq n U$ and $\raw{U'}s=V$. 
Note that $U'$ is an \dd ncollage over $\rP$ by  
Lemma~\ref{c2n}.
 
On the other hand, by Lemma \ref{mokr}\ref{mokr2}, 
there is a \mut\ $\zq\in\mP$ such that $\abs{\zq}= B$ 
and $[\zq]\sq G$.
By Lemma \ref{tu}\ref{tu3}, 
there is a \mut\ $\zp'\in\mtt_B$ satisfying 
$\zp'\leq_m\zp$ and $\rau{\zp'}u=\zq$. 
Note that $\zp'$ is an \dd mcollage over $\mP$ by 
Lemma~\ref{mokr}\ref{mokr3}. 
Thus $\zp'$ and $U'$  
witness \ref{nf'1} partially: 
$[\raw{U'}s]\cap\imb f{[\rau{\zp'}u]}=\pu$ holds,  
not yet $[U']\cap\imb f{[\zp']}=\pu$.
However this procedure can be iterated, going over  
all pairs of strings $u\in2^m$, $s\in 2^n$. 
This leads to the result required.

\ref{nf'2}
As in the first part, it suffices,  
given a pair of strings $r,s\in2^m$ 
(possibly $r=s$), 
to find an \dd mcollage $\zp'\in\mtt_B$ over $\mP$, 
satisfying $\zp'\leq_m\zp$   
and $[\rau{\zp'}s(\xi)]\cap\imb f{[\rau{\zp'}r]}=\pu$.
The tree $T=\zp(\xi)$ belongs to $\jpi(\xi)\sq\lt$, 
and  $\rau\zp s(\xi)=\raw T{s'}$, 
$\rau\zp r(\xi)=\raw T{r'}$, where  
$s'=\qa s\xi$, $t'=\qa t\xi$ are strings of length 
$n=\qi m\xi$, see Definition \ref{phi}. 
Now 
%найдется такой кортеж $\sg\in\bse,$ что 
$\raw T{s'}=\ta\ap\raw T{r'}$ 
by Lemma \ref{tran}, where $\ta=u[{s'},T]\ap u[{r'},T]$. 
But $f$ is not simple  
on $\rau\zp r$, hence there exists a point 
$x_0\in\rau\zp r$ such that $f(x_0)\ne \ta\ap x_0(\xi)$. 
We have two strings $v\ne w$ in $\bse$ 
of equal length $\lh v=\lh w>\lh\ta$, satisfying   
$v\su f(x_0)$ and $w\su\ta\ap x_0(\xi)$. 
We put $w'=\ta\ap w$; then $w'\su x_0(\xi)$. 

But $f$ is continuous, hence using Lemma \ref{mokr} 
as above, we find a \mut\ $\zq\in\mP$ 
such that $\abs\zq=B$, 
$\zq\leq\rau\zp r$, and if $x\in[\zq]$ then 
$v\su f(x)$, $w\su\ta\ap x(\xi)$, $w'\su x(\xi)$.
And further we find a \mut\ $\zp'\in\mtt_B$  
satisfying $\zp'\leq_m\zp$ and $\rau{\zp'}r=\zq$, 
and being an \dd mcollage over $\mP$. 

We claim that 
$[\rau{\zp'}s(\xi)]\cap\imb f{[\rau{\zp'}r]}=\pu$.
Indeed by construction if 
$x\in[\zq]=[\rau{\zp'}r]$ thrn $v\su f(x)$. 
Thus it remains to check that $w\su b$ for all
$b\in [\rau{\zp'}s(\xi)]$. 
Note that $\rau{\zp'}s(\xi)=\raw{T'}{s'}$ and  
$\rau{\zp'}r(\xi)=\raw{T'}{r'}$, where  
$T'=\zp'(\xi)\in\jpi(\xi)$.
On the other hand $T'$ is a tree in $\lt$ and 
$T'\nq n T$, hence $\raw{T'}{s'}=\ta\ap\raw{T'}{r'}$ 
by Lemma \ref{nadd}. 
Thus if $b\in [\rau{\zp'}s(\xi)]$ then 
$a=\ta\ap b\in [\rau{\zp'}r(\xi)]=[\raw{T'}{r'}]$. 
It follows that $w'\su a$ by the choice of 
$\zq=\rau{\zp'}r$. 
Then $w\su b=\ta \ap a$ (since $w=\ta\ap w'$), 
as required. 
\epf

\parf{Extension of \muf s} 
\las{gex}

The forcing notion for the proof 
of Theorem \ref{Tun} will be defined as  
an \dd\omi union of an increasing \dd\omi 
sequence of \muf s.
Definition \ref{dfex} below contains 
conditions which every step of the construction 
will have to obey.
We begin with the following definition.

\bdf
[coding continuous maps]
\lam{nfk}
Let $B\sq\omi$ be at most countable.
A \rit{code of a continuous map} $\ddn B\to \dn$
\index{code of a continuous map}%
\index{map!code}%
is an indexed family 
$\bc=\sis{U^\bc_i(k)}{k<\om,i=0,1}$ of finite sets 
$U^\bc_i(k)\sq\mpo B$ (see Example \ref{mcoh}), 
such that for all $k$: \vom 

(1) 
if $\zp\in U^\bc_0(k)$ and $\zq\in U^\bc_0(k)$ then 
$[\zp\ares B]\cap[\zq\ares B]=\pu$, and \vom 

(2)
$\bigcup_{k<\om,i=0,1}\bigcup_{\zp\in U^\bc_i(k)}
[\zp\ares B]=\ddn B$.\vom

\noi
Let $\knf B$ denote the set of all such codes. 
\index{zzccfb@$\knf B$}%
\index{code!$\knf B$}%
\index{zzccf@$\kng$}%
\index{code!$\kng$}%
\index{zzcII@$\abs\bc$}%
\index{code!$\abs\bc$}%

We set  
$\kng=\bigcup_{B\sq\omi,\,\card B\le\alo}\knf B$, 
and if 
$\bc\in\knf B$ then $\abs\bc= B$.

The coded map $f=f^\bc:\ddn B\to\dn$ itself is  
defined as follows in this case: 
\index{zzfc@$f^\bc$}%
\index{map!$f^\bc$}%
$f^\bc(x)(k)=i$, if there is a \mut\ $\zp\in U^\bc_i(k)$  
such that $x\in[U\ares B]$.
Make use of (1) to show that the definition is sound.
\edf

We skip a routine proof of the following lemma, based 
on the compactness of the spaces considered.

\ble
\lam{knf}
If\/ $B\sq\omi$ is countable, $X\sq\ddn B$ closed, 
and a map\/ $f:X\to\dn$ is continuous, then there is 
a code\/ $\bc\in\knf B$ such that\/ $f=f^\bc\res X$.\qed
\ele

\bdf
[in $\rL$]
\lam{dfex}
Let $\cM$ be a countable transitive model of 
theory $\ZFC'$,  
\index{zzzfcp@$\zfc'$}%
\index{theory!zfcp@$\zfc'$}%
which includes all $\ZFC$ axioms except for the 
power set axiom, but with the axiom 
which claims the existence of $\pws\om$. 
(This implies the existence of the ordinal $\omi$ 
and sets like $\dn$, $\pet$, $\pel$ of cardinality 
$\mathfrak c=2^{\alo}$.)

Recall that $\rL_\al$ is \dd\al th level of the G\"odel  
constructible hierarchy. 

Let $\jpi\in\cM$ be a regular (small) \muf. 
Then $\abs\jpi=\bb\in\cM$ and $\al=\tsup\bb=\bigcup\bb<\omi$. 
\kmar{smf jpi}% 
\index{submultiforcing!$\smf\jpi$}%
\index{zzSP@$\smf\jpi$}%
We let $\smf\jpi$ denote the closure of  
$\mt\jpi\cap\rL_\al$ in $\mt\jpi$ 
with respect to the three operations of 
Definition \ref{mfd}. 
Thus $\smf\jpi\in\cM$, $\smf\jpi\sq\mt\jpi$,  
$\smf\jpi$ is a countable \smuf. 

Note that $\smf\jpi$ does not depend on $\cM$.

A \muf\ $\jqo$ (not necessarily in $\cM$)
is an \dd\cM\rit{extension} of $\jpi$, 
in symbol $\jpi\rfa\cM\jqo$, 
\kmar{rfa M}%
if 
\index{zzPsuMQ@$\jpi\rfa\cM\jqo$}%
\index{multiforcing!extension}%
\index{Mextension@\dd\cM extension of \muf, $\rfa\cM$}%
the following holds:
\ben
\Aenu
\itlb{dfex1}%
$\abs\jqo=\abs\jpi$ and $\jqo$ is a small \muf; 

\itlb{dfex3}%
if\/ $\xi\in\abs\jpi$ then
$\jpi(\xi)\rfa\cM \jqo(\xi)$ in the sense of 
Definition \ref{fm};

\itlb{dfex5}%
if\/ $\zp\in\smf\jpi$  
then there is a \mut\ $\zq\in\mt\jqo$ satisfying  
$\zq\leq\zp$ and $\zq\sqdv \zD$ for all 
open dense sets $\zD\sq\smf\jpi$, $\zD\in\cM$;

\itlb{dfex6}%
if $\zp\in\smf\jpi$, $\xi\in \abs\zp$, 
a map $f:\ddn{\abs\zp}\to\dn$ is continuous and 
has a code in $\knf{\abs\zp}\cap\cM$, 
then there exists 
a \mut\  
$\zq\in\mt\jqo$ such that  
$\abs\zq=\abs\zp$, $\zq\leq\zp$, 
and either (i) 
there is a string $\sg\in\bse$ such that 
$f(x)=\sg\ap x(\xi)$ for all $x\in[\zq]$,  
or (ii) 
$f(x)\nin[U]$ for all $x\in[\zq]$ and 
$U\in\jqo(\xi)$.
\qed
\een
\eDf

\vyk{
Смысл условия \ref{dfex6} можно в целом 
пояснить так:  
если $\etz<\omi$ и $\rc$ есть \dd{\mt\jpi}имя точки из $\dn$ 
то продолженный форсинг $\mt\jpq$ вынуждает, что
$\rc$ не принадлежит к множествам вида $[U]$, где 
$U$ --- дерево в $\pc\jqo\etz$ --- если $\rc$ не есть 
имя $\rpi_{\etz}$ одной из точек вида $\sg\app x_{\etz}[G]$, 
где $\sg\in\bse$.
}

\bte
[in $\rL$]
\lam{tex}
Let\/ $\cM$ be a countable transitive model of\/ $\ZFC'$,  
and\/ $\jpi\in\cM$ be a regular (small) \muf. 
% and\/ $\abs\jpi=\bb=\ens{\xi}{\xi<\bb}\,;$ $\bb<\omi$.
%so that\/ $\bse\in\jpi(\xi)$ для всех\/ $\xi\in\abs\jpi$. 
% 
Then there is 
%a small\/ \muf\ $\jqo$ which is 
an\/
\dd\cM extension\/ $\jqo$ of\/ $\jpi$. 
% в смысле\/ \ref{dfex}.
\ete

The proof of the theorem follows in the two next  
sections. 
The construction of $\jqo$ is presented in  
Section~\ref{jex}, the proof of its    
properties follows in Section \ref{pres}.

\parf{The construction of extending \muf} 
\las{jex}

The following definitions formalize construction   
of generic \mut s for the proof of Theorem~\ref{tex}, 
by means of Lemma~\ref{fusm}. 
\bit
\item
{\ubf Arguing under the conditions of Theorem \ref{tex}}, 
we let $\bb=\abs\jpi$ and $\mP=\smf\jpi$, 
so that $\bb<\omi$ and $\mP\sq\mt\jpi$ is a countable 
\smuf. 

\item 
During the course of the proof of Theorem \ref{tex}, 
\ie, to the end of Section\ref{pres}, we 
{\ubf fix a \dd\bb complete function $\phi:\om\onto\bb$}. 
This allows to use  
{\ubf the notation of Definition \ref{phi}}.
\eit 
To begin with, we reduce all \mut s $\zp\in\mP$ 
to the domain $\bb$, substituting each of them  
by its copy $\uar\zp=\zp\ares\bb$ 
\index{zzTI@$\uar\zp$}%
\index{multitree!$\uar\zp$}%
\kmar{uar zp}%
(see Definition \ref{proj}). 
Thus, by the regularity of $\jpi$, we have 
$\uar\zp\in\mt\jpi$ and  
$\abs{\uar\zp}=\bb$, 
%(счетное подмножество $\omi$), 
and by definition $\uar\zp(\xi)=\zp(\xi)$ for 
$\xi\in\abs\zp$, but $\uar\zp(\xi)=\bse$ for 
$\xi\in\bb\bez\abs\zp$. 
We put $\ump=\ens{\uar\zp}{\zp\in\mP}$, 
\index{zzSI@$\ump$}%
\index{submultiforcing!$\ump$}%
this is a \smuf, too. 
%, и если 
\kmar{ump}%
%$\zp'\in\udp$ то пусть 
%$\abt{\zp'}=\ens{\xi\in\bb}{\zp'(\xi)\ne\bse}$. 
%Имеем $\abt{\uar\zp}\sq\abs\zp\sq \bb$. 

\bdf
\lam{mut}
A \rit{\mus\ (over\/ $\ump$)} is any function  
$\vpi:\dom\vpi\to\mtt_\bb$ 
%of the form 
\index{system}% 
\index{system!$\dom\vpi$}% 
where $\dom\vpi\sq\om\ti\om$ is finite, 
and if $\ang{k,m}\in\dom\vpi$ then 
\ben
%\atc\atc
\nenu
\itlb{mut1}%
if $n<m$ then $\ang{k,n}$ also belongs to $\dom\vpi$;    
\index{zzfII@$\dom\vpi$}%
%\index{база!$\abs\vpi$}%

\itlb{mut2}%
$\vpi(k,m)$ is a tree in $\mtt_\bb$ 
and a \dd{m}collage over $\ump$, 
and $\abs{\vpi(k,m)}=\bb$;

\itlb{mut3}%
if $m>0$ then $\vpi(k,m)\leq_m\vpi(k,m-1)$.
\een
In this case, let $\zn\vpi k$ denote the largest  
\kmar{zn vpi k}%
number $m$ satisfying $\ang{k,m}\in \dom\vpi$,  
but $\zn\vpi k=-1$ if there is no such $m$. 
Let $\abs\vpi=\ens{k}{\zn\vpi k\ge0}$, a finite set. 

Let $\ms\ump$ denote the set of all \mus s.
\kmar{ms udp}%
\index{zzsysP@$\ms\ump$}%
\index{system!$\ms\ump$}%

A \mus\ $\vpi$ \rit{extends} a \mus\ $\psi$, 
\index{system!extension, $\psi\cle\vpi$}% 
in symbol $\psi\sq\vpi$, 
if $\dom \psi\sq\dom\vpi$ and $\psi=\vpi\res{\dom\psi}$; 
while $\psi\su\vpi$ will denote strict extension.
\edf  

\ble
[elementary]
\lam{nwm}
Suppose that\/ $\vpi\in\ms\ump$. 
Then
\ben
\renu
\itlb{nwm1}%
if\/ $k\in\abs\vpi$ and\/ $m=\zn\vpi k$ then 
the extension\/ $\vpi'$ of \mus\/ $\vpi$ by\/
$\zn{\vpi'}k=m+1$ and\/ $\vpi'(k,m+1)=\vpi'(k,m)$ 
is a \mus\ extending\/ $\vpi\;;$ 

\itlb{nwm-1}%
if\/ $k\nin\abs\vpi$ and\/ $\zp\in\ump$, then 
the extension\/ $\vpi'$ of \mus\/ $\vpi$ by\/
$\dom{\vpi'}=\dom\vpi\cup\ans{\ang{k,0}}$  
and\/ $\vpi'(k,0)=\zp$,  
is a \mus\ extending\/ $\vpi$.\qed 
\een
\ele

\bdf
\lam{muss}
(A) 
Let $\df$ denote the set of all sets 
$X\sq\hc$, definable in $\hc$ 
\kmar{df, \ hc}%
\index{zzDEF@$\df$}%
\index{zzHC@$\hc$}%
\index{hereditarily countable sets, $\hc$}%
\index{set!hereditarily countable, $\hc$}%
(= all hereditarily countable sets) 
by \dd\in formulas with parameters in $\cM\cup\ans{\cM,\phi}$. 
As $\df$ is countable, Lemma \ref{nwm} 
allows to define an infinite system 
$\Phi:\om\ti\om\to\mtt_\bb$  
satisfying the requirements \ref{mut2} and \ref{mut3} 
of Definition \ref{mut} on the whole domain $k,m<\om$, 
and also satisfying the following 
\rit{genericity condition}: 
every set $\Da\in\df$ is \rit{blocked} by 
\index{blocks} 
one of \mus s $\vpi\in\ms\ump$, $\vpi\su\Phi$, 
in the sense that:\vom 

--- either 
(I) 
$\vpi\in\Da$,\vom 

--- or   
(II) 
there is no \mus\
$\psi\in\ms\ump\cap\Da$ extending $\vpi$.\vom 

\noi
We let $\ft km=\Phi(k,m)$ for all $k,m<\om$.\vom

(B)
The \rit{limit trees}
$\fl k=\bigwedge_m\ft km$, 
defined by $\abs{\fl k}=\bb$ 
\kmar{fl k}%
and $\fl k(\xi)=\bigcap_m\ft km(\xi)$ for all $\xi\in\bb$, 
belong to $\mlt_\bb$ and satisfy 
$\fl k\leq_{m+1}\ft km$ for all $k,m$ 
by Lemma~\ref{fusm}. 
Accordingly if $\xi\in\bb$ then $\fl k(\xi)\in\pel$ 
and $\fl k(\xi)\nq n\ft km(\xi)$ for all $m$, 
where $n=\qi m\xi$ (Definition \ref{phi}). 
This means $\raw{\fl k(\xi)}s\sq \raw{\ft km(\xi)}s$
for all $s\in2^n$.\vom 

(C) 
If $\xi\in\bb$ then the set   
$\rQ_\xi=
\ens{\sg\app\raw{\fl k(\xi)}s}{k<\om\land \sg,s\in\bse}$ 
is a countable \sfo, see Example \ref{xcoh}. 
We define a small \muf\ $\jqo$ 
by $\abs\jqo=\bb$ and 
$\jqo(\xi)=\rQ_\xi$ for all $\xi\in\bb$.
\edf
 
We'll check that the 
\muf\ $\jqo$ satisfies all conditions of 
Definition \ref{dfex}. 
Note that \ref{dfex}\ref{dfex1}
directly holds by construction. 
The following lemma is obvious since  
option (II) of Definition \ref{muss}(A)  
is impossible for dense sets $\Da$. 
It will be a key ingredient in the verification 
of other conditions below.

\ble
\lam{12g}
Let a set\/ $\Da\in\df$, $\Da\sq\ms\ump$, 
be dense in\/ $\ms\ump$, that is, every \mus\ 
in\/ $\ms\ump$ is extendable to a \mus\ in\/ $\Da$. 
Then there exists a \mus\/ $\vpi\in\Da$  
satisfying\/ $\vpi\su\Phi$.\qed
\ele

\bcor
\lam{12c}
If\/ $\zp\in\ump$ then there is an index\/ $k$ 
such that\/ 
$\fl k\leq\ft k0=\zp$. 
If\/ $\xi\in\bb$ and\/ $T\in\jpi(\xi)$ 
then there is an index\/ $k$ 
such that\/ $\fl k(\xi)\sq\ft k0(\xi)=T$.  
\ecor
\bpf
Consider the set $\Da$ of all  
\mus\ $\vpi\in\ms\ump$ such that $\vpi(k,0)=\zp$    
holds for at least one $k\in\abs\vpi$.
As $\zp\in\ump\in\cM$, the set $\Da$ 
belongs to $\df$. 
We claim that $\Da$ is dense in $\ms\ump$. 
Indeed let $\vpi\in \ms\ump$. 
Take any $k\nin\abs\vpi$. 
By Lemma \ref{nwm}\ref{nwm-1} there is a \mus\ 
$\psi\in\ms\ump$ extending $\vpi$ and satisfying 
$\ang{k,0}\in\dom\psi$ and $\psi(k,0)=\zp$. 
Thus $\psi\in\Da$, and the density is proved.

By Lemma \ref{12g}, there is a \mus\  
$\vpi\in\Da$, $\vpi\su\Phi$. 
Then $\ft k0=\vpi(k,0)=\zp$ 
for some $k$. 
But $\fl k$ satisfies $\fl k\leq\ft k0$ 
by \ref{muss}(B), as required. 

To reduce the second claim to the first one, 
note that  
if $\xi\in\bb$ and $T\in\jpi(\xi)$ then by 
definition there is a \mut\ $\zp\in\ump$ 
satisfying $\zp(\xi)=T$.
\epf

\parf{Verification of requirements} 
\las{pres}

We check conditions of 
Definition \ref{dfex} for $\jqo$ 
in the context of Section \ref{jex}.\vom

{\ubf Validation of \ref{dfex}\ref{dfex3}.}
%\lam{uu3}
Fix $\xi\in\bb$. 
To check \ref{fm1} of Definition \ref{fm} 
(the density of  
$\pc\jqo\xi$ in $\pc\jqo\xi\cup\pc\jpi\xi$),
let $T\in\pc\jpi\xi$. 
Then $\fl k(\xi)\sq T$ for some $k$
by Corollary \ref{12c}. 
But the tree $S= \fl k(\xi)$ belongs to  
$\rQ_\xi=\jqo(\xi)$ 
by \ref{muss}(C), as required.
%\vom

Now assume that $\xi\in\bb$, a set $D\in\cM$, 
$D\sq\pc\jpi\xi$ is pre-dense in $\pc\jpi\xi$, 
and $U\in\pc\jqo\xi$. 
Prove $U\sqf\bigcup D$. 
By definition, $U=\sg\app\fls k\xi s$, 
where $k<\om$, $\xi\in\bb$, and $s,\sg\in\bse.$ 
We can assume that $\sg=\La$, 
\ie, in fact just $U=\fls k\xi s$. 
(The general case is reduced to $U=\fls k\xi s$ 
by the substitution of $\sg\app D$ for $D$.) 
Furthermore, we can assume that $s=\La$, 
\ie, $U=\fl k(\xi)$, because  
$\fls k\xi s\sq\fl k(\xi)$.
Thus let $U=\fl k(\xi)$. 
The index $k$ will be fixed.

It follows, from the pre-density of $D$ and property   
\ref{mfd}\ref{mfd1} of the \smuf\ $\ump$, 
that the set $\zD\in\cM$, of all  
\mut s $\zp\in\ump$ satisfying $\zp(\xi)\sq V$ 
for some $V\in D$, is itself open dense in $\ump$.

We claim that the set $\Da\in\cM$ 
of all \mus\ $\vpi\in \ms\ump$, such that  
$k\in\abs\vpi$, and for every string $t\in 2^{n}$, 
where $n=\zn\vpi k$, the \mut\ $\rau{\vpi(k,n)}t$ 
belongs to $\zD$, is dense in $\ms\ump$. 
Indeed let $\vpi\in \ms\ump$. 
By Lemma \ref{nwm}\ref{nwm-1}, we assume that  
$k\in\abs\vpi$, \ie, $n'=\zn\vpi k\ge0$. 
By definition the \mut\ $\zp=\vpi(k,n')$ is an  
\dd{n'}collage over $\ump$, and then, 
by Lemma \ref{mokr}\ref{mokr1}, \dd{n}collage, 
too, where $n=n'+1$. 
Then by Lemma \ref{mokr}\ref{mokr4} there is a  
\mut\ $\zp'\in\mlt_\bb$, which is an  
\dd{n}collage over $\ump$ and satisfies  
$\zp'\leq_n\zp$ and $\rau{\zp'}t\in\zD$ 
for all $t\in2^n$. 
Extend $\vpi$ to a \mus\ $\psi$ by  
$\dom\psi=\dom\vpi\cup\ans{\ang{k,n}}$ 
and $\psi(k,n)=\zp'$; 
we have $\psi\in\Da$. 

Now by Lemma~\ref{12g} there is a \mus\  
$\vpi\in\Da$ satisfying $\vpi\su\Phi$. 
Then $\rau{\vpi(k,n)}t=\rau{\ft kn}t\in\zD$ for  
all $t\in 2^{n}$, where $n=\zn\vpi k$, thus  
$\ft kn\sqdv\zD$, hence  
$\fl k\sqdv\zD$. 
Therefore $U=\fl k(\xi)\sqf\bigcup D$
by the definition of $\zD$.\vom    

{\ubf Validation of  \ref{dfex}\ref{dfex5}.}
Assume that   
$\zD\in\cM,$ $\zD\sq\mP$ is open dense in $\mP$. 
Accordingly the set 
$\uzd=\ens{\uar\zp}{\zp\in\zD}\sq\ump$ 
is dense in $\ump$.\snos
{To prove the openness let $\zp\in\zD$.  
Then $\uzp\in\uzd$, 
$\zq\in\mP$, and $\uzq\leq\uzp$. 
We cannot assert directly that $\zq\leq\zp$. 
However the \mut\ $\zq'=\zq\ares{(\abs\zp\cup\abs\zq)}$ 
also belongs to $\mP$ by Definition~\ref{mfd}\ref{mfd3}. 
Note that $\uzq\leq\uzp$ easily implies $\zq'\leq\zp$. 
Therefore $\zq'\in\zD$, since $\zD$ is open. 
We conclude that 
$\uzq=\uar{\zq'}\in\uzd$.}
By Corollary \ref{12c}, it suffices to prove  
that $\fl k\sqdv\uzd$ for all $k<\om$.

By the open-density of $\uzd$, the set $\Da_k\in\cM$ 
of all \mus s $\vpi\in \ms\ump$ such that   
$k\in\abs\vpi$,and for every string $t\in 2^{n}$, 
where $n=\zn\vpi k$, the \mut\ $\rau{\vpi(k,n)}t$  
belongs to $\uzd$, is dense in $\ms\ump$. 
(See the verification of \ref{dfex}\ref{dfex3} above.) 
By Lemma \ref{12g} there exists a \mus\  
$\vpi\in\Da_k$ satisfying $\vpi\su\Phi$. 
Then $\rau{\vpi(k,n)}t=\rau{\ft kn}t\in\zD$ for  
all $t\in 2^{n}$, where $n=\zn\vpi k$, that is, 
$\ft kn\sqdv\uzd$ holds, hence $\fl k\sqdv\uzd$, 
as required.\vom
 
{\ubf Validation of \ref{dfex}\ref{dfex6}.}
%
%\lam{K}
Let $\zp\in\mP$, $\xi\in C=\abs\zp$, 
$\bc\in\knf C\cap\cM$, and $f=f^\bc$ 
(a continuous map $\ddn C\to\dn$). 
The \mut\ $\uzp=\zp\ares\bb$ belongs to $\ump$, 
and the map $\uff(x)=f(x\res C):\ddn \bb\to\dn$ 
is continuous.
In terms of Section \ref{opo}, we can assume that  
(*) there is no \mut\ $\zp'\in\ump$, $\zp'\leq\uzp$, 
such that $\uff$ is simple for $\xi$ on $\zp'$. 
Indeed otherwise using Corollary \ref{12c} 
we get a \mut\ $\zq$ of the form  
$\fl k$, satisfying $\fl k\leq\zp'$, and hence 
(i) of \ref{dfex}\ref{dfex6}. 

Now assuming (*) we accordingly prove that  
any \mut\ $\zq=\fl k$ with $\fl k\leq\ft k0=\uzp$  
satisfies (ii) of \ref{dfex}\ref{dfex6}.
%$f(x\res C)\nin[U]$ for all $x\in[\zq]$ and 
Let $U\in\jqo(\xi)=\rQ_\xi$, 
and we have to prove that $\uff(x)\nin[U]$ for all  
$x\in[\fl k]$.
By definition, $U=\ta\ap\rau{\fl\ell(\xi)}s$, 
where $\ta,s\in\bse$ and $\ell<\om$.
Now, as $\rau{\fl\ell(\xi)}s\sq\fl\ell(\xi)$, 
we can assume that $s=\La$, 
that is, $U=\ta\ap\fl\ell(\xi)$. 
Moreover we can assume that $\ta=\La$, \ie, 
$U=\fl\ell(\xi)$; otherwise consider  
the map $f'(x)=\ta\ap\uff(x)$ instead of $\uff$. 

Thus we fix an index  $\ell<\om$ and prove that
%$\uff(x)\nin[\fl\ell(\xi)]$ для всех $x\in[\fl k]$. 
%
$[\fl\ell(\xi)]\cap\imb\uff{[\fl k]}=\pu$.\vom

{\ubf Case 1:} $\ell\ne k$. 
Consider the set $\Da$ of all \mus s  
$\vpi\in \ms\ump$ such that $k,\ell\in\abs\vpi$, 
that is, 
%так что 
$m=\zn\vpi k\ge0$ and $n=\zn\vpi\ell\ge0$, and 
$[\vpi(\ell,n)(\xi)]\cap\imb\uff{[\vpi(k,m)]}=\pu$. 

\ble
\lam{ls1}
The set\/ $\Da$ is dense in\/ $\ms\ump$. 
\ele
\bpf[Lemma]
Let $\vpi\in \ms\ump$. 
By Lemma \ref{nwm}\ref{nwm-1}, we can assume that 
$k,\ell\in\abs\vpi$, that is, $n'=\zn\vpi\ell\ge0$ 
and $m'=\zn\vpi k\ge0$. 
By definition, the \mut\ $\zr'=\vpi(k,m')$ is an  
\dd{m'}collage over $\ump$, and then \dd{m}collage, 
too, by Lemma \ref{mokr}\ref{mokr1}, where $m=m'+1$. 

Further, we can assume that $\phi(n')=\xi$, 
for if not then take the least number  
$n''>n'$ satisfying $\phi(n'')=\xi$, and  
trivially extend the \mus\ $\vpi$ by 
$\vpi(\ell,j)=\vpi(\ell,n')$ for all $\ell$ with  
$n'<\ell\le n''$.
As above, the \mut\ $\bZ'=\vpi(\ell,n')$ is an  
\dd{n'}collage over $\ump$, and hence \dd{n}collage, 
where $n=n'+1$.
It follows that $\rau{\bZ'}\sg\in\ump$ 
for all $\sg\in2^n.$
In particular $\rau{\bZ'}\sg(\xi)\in\jpi(\xi)$ for   
$\sg\in2^n.$
Yet $\rau{\bZ'}\sg(\xi)=\raw{{\bZ'}(\xi)}{\qa\sg\xi}$ 
by Definition \ref{phi}, where $\qa\sg\xi\in2^\nu$ and  
$\nu=\qi m\xi$. 
Therefore the tree $Z'={\bZ'}(\xi)$ 
is an \dd\nu collage over 
$\jpi(\xi)$.

By Theorem \ref{nf'}\ref{nf'1}, there  
exist a  \mut\ $\zr\in\mtt_\bb$  
and a tree  $Z\in\lt$ such that 
$\zr\leq_m\zr'$, $Z\nq\nu Z'$, 
$\zr$ is an \dd mcollage over $\ump$, 
$Z$ is an \dd\nu collage over $\mP(\xi)$, 
and $[Z]\cap\imb \uff{[\zr]}=\pu$.
Define a \mut\ $\bZ\in\mtt_\bb$ so that  
$\bZ(\xi)=Z$ and $\bZ(\et)=\bZ'(\et)$ for 
all $\et\in\bb$, $\et\ne\xi$.

\bsl
\lam{znc}
$\bZ$ is an\/ \dd ncollage over\/  $\ump$ 
and\/ $\bZ\leq_n\bZ'$.
\esl
\bpf 
Let $\sg=\ta\we i\in2^n,$ where 
$\ta\in2^{n'}$ and $i=0,1$. 
The strings $\qa\sg\et\in2^{\qi{n}{\et}}$ and 
$\qa\ta\et\in2^{\qi{n'}{\et}}$ 
(Definition \ref{phi}) 
are connected so that: 
$\qa\sg\et=\qa\ta\et$ for $\et\ne\xi$, but  
$\qa\sg\xi=(\qa\ta\xi)\we i$, since $\phi(n')=\xi$ 
and $n=n'+1$.
It follows that 
$$
\rau{\bZ}\sg(\et)=\raw{\bZ(\et)}{\qa\sg\et}=
\raw{{\bZ'}(\et)}{\qa\sg\et}=\rau{\bZ'}\sg(\et) 
$$ 
for $\et\ne\xi$, that is, 
$\rau{\bZ}\sg\res{(\bb\bez\ans\xi)}=
\rau{\bZ'}\sg\res{(\bb\bez\ans\xi)}$.
Further, 
$\rau{\bZ}\sg(\xi)=\raw{\bZ(\xi)}{\qa\sg\xi}=
\raw{Z}{\qa\sg\xi}=
\raw{\raw{Z}{\qa\ta\xi}}i \in\jpi(\xi)
$, 
since $Z$ is a \dd\nu collage over $\mP(\xi)$.
This implies $\rau{\bZ}\sg\in\ump$ by the property
\ref{mfd}\ref{mfd1} of \smuf s.
As $\sg\in2^n$ is arbitrary, 
$\bZ$ is a \dd ncollage over $\ump$. 

To establish $\bZ\leq_n\bZ'$, we need 
(in the same notation)  
to prove $\rau{\bZ}\sg\leq\rau{\bZ'}\sg$ 
for all $\sg\in2^n$, that is, 
$\rau{\bZ}\sg(\et)\sq\rau{\bZ'}\sg(\et)$ for all $\et\in\bb$.
If $\et\ne\xi$ then simply  
$\rau{\bZ}\sg(\et)\sq\rau{\bZ'}\sg(\et)$, 
as above. 
Further, we have $\rau{\bZ}\sg(\xi)=\raw Z{s}$ and 
$\rau{\bZ'}\sg(\xi)=\raw {Z'}{s}$, 
where $s=\qa\sg\xi\in2^\nu$, $\nu=\qi m\xi$. 
But $Z\nq\nu Z'$ by construction, hence 
$\raw Z{s}\sq\raw {Z'}{s}$, or equivalently, 
$\rau{\bZ}\sg(\xi)\sq\rau{\bZ'}\sg(\xi)$. 
Thus $\rau{\bZ}\sg(\et)\sq\rau{\bZ'}\sg(\et)$ 
for all $\et\in\bb$, 
that is, $\rau{\bZ}\sg\leq\rau{\bZ'}\sg$, 
as required.
\epF{Sublemma}
 
Coming back to the lemma, 
we extend $\vpi$ to a \mus\ $\psi$ with  
$\dom\psi=\dom\vpi$, 
$\zn\vpi k=m$, $\zn\vpi \ell=n$, 
$\psi(k,m)=\zr$, and $\psi(\ell,n)=\bZ$ 
(just two new values). 
Thus $\psi$ is a \mus\ in $\ms\ump$. 
Indeed $\psi(k,m)=\zr$, one of the two new  
terms relatively to $\vpi$, is an \dd mcollage 
over $\ump$, and $\zr\leq_{m}\zr'=\vpi(k,m')$, 
where $m=m'+1$, as required by \ref{mut}\ref{mut3}. 
Similarly for $\psi(\ell,n)=\bZ$ the other new term. 
Thus $\psi\in\ms\ump$ and clearly 
$\vpi\cle\psi$.
Finally $[Z]\cap\imb \uff{[\zr]}=\pu$  by construction, 
hence $\psi\in\Da$. 
This ends the proof of the density of $\Da$.
\epF{Lemma}

Now Corollary \ref{12g} yields a \mus\ 
$\vpi\in\ms\ump$, $\vpi\su\Phi$. 
Then $k,\ell\in\abs\vpi$, hence   
$m=\zn\vpi k\ge0$ and $n=\zn\vpi\ell\ge0$, 
and \mut s $\ft km=\vpi(k,m)$, $\ft\ell n=\vpi(\ell,n)$ 
satisfy 
$[\ft\ell n(\xi)]\cap\imb\uff{[\ft km]}=\pu$ 
by the definition of $\Da$, therefore 
$[\fl k(\xi)]\cap\imb\uff{[\fl k]}=\pu$, 
because $\fl k\sq\ft km$, 
as required.\vom

{\ubf Case 2:} $\ell=k$. 
Consider the set $\Da$ of all \mus s 
$\vpi\in \ms\ump$ such that $k\in\abs\vpi$    
(and then $m=\zn\vpi k\ge0$) and
$[\vpi(k,m)(\xi)]\cap\imb\uff{[\vpi(k,m)]}=\pu$. 
We don't claim that $\Da$ is dense. 
However by Definition~\ref{muss} there is a \mus s 
$\vpi\in\ms\ump$, $\vpi\su\Phi$, \rit{blocking}
$\Da$ in the sense of \ref{muss}(A), (I)$\lor$(II).

We now assert that \ref{muss}(A)(II) is impossible 
for $\vpi$. 
Indeed let $m'=\zn\vpi k$ and 
$\zr'=\vpi(k,m')=\Phi(k,m')=\ft k{m'}$. 
Then $\zr\sq\ump=\ft k0$, and hence by (*) 
(in the beginning of validation of \ref{dfex}\ref{dfex6}), 
if $\zp'\in\ump$, $\zp'\leq\zr'$, then the map 
$\uff$ is not siple for $\xi$ on $\zp'$. 
Therefore by Theorem \ref{nf'}\ref{nf'2} there is a  
\mut\ $\zr\in\mt\jpi$, which is a \dd mcollage over  
$\ump$, where $m=m'+1$, and satisfies $\zr\leq_m\zr'$ 
and $[\zr(\xi)]\cap\imb{\uff}{[\zr]}=\pu$.
As in Case 1, we can extend $\vpi$ to a \mus\  
$\psi\in\ms\ump$ 
with the only one new term $\psi(k,m)=\zr$, 
and then $\psi\in\Da$ by the choice of $\zr$. 
This proves that \ref{muss}(A)(II) 
cannot happen for $\vpi$. 

Thus \ref{muss}(A)(I) takes place, that is, 
$\vpi\in\Da$. 
It follows that 
$[\vpi(k,m)(\xi)]\cap\imb\uff{[\vpi(k,m)]}=\pu$, 
hence $[\ft km(\xi)]\cap\imb\uff{[\ft km]}=\pu$. 
This implies $[\fl k(\xi)]\cap\imb\uff{[\fl k]}=\pu$, 
since $\fl k\leq\ft km$, as required.\vom
\qeD{Theorem~\ref{tex}}

\parf{The forcing} 
\las{jfor}

{\ubf We argue in the constructible universe $\rL$}
in this section. 

We begin with some definitions related to 
sequences of \muf s.

First of all, we somewhat generalize the 
definition of $\rfa\cM$ in \ref{dfex}. 
Given small \muf s $\jpi,\jqo$ and a model $\cM$, 
we define $\jpi\rfp\cM\jqo$, when $\abs\jpi\sq\abs\jqo$ 
\kmar{rfp cM}%
and $\jpi\rfa\cM(\jqo\res\abs\jpi)$ 
in the sense of \ref{dfex}.   
\index{zzPsu+MQ@$\jpi\rfp\cM\jqo$}%
\index{multiforcing!extension}%
%\index{mextension@\dd\cM extension of a \muf}%
\index{Mextension@\dd\cM extension of \muf, $\rfp\cM$}%
%а если $\xi\in\abs\jqo\bez\abs\jpi$ то $\jqo(\xi)$ --- 
%регулярный форсинг (\te\ содержит дерево $\bse$).
If 
$\vjpi=\sis{\jpi_\al}{\al<\la}$ ($\la<\omi$) 
is a sequence of small \muf s $\jpi_\al$ then: 
\ben
\aenu
\itlb{jfor1}% 
$\cM(\vjpi)$ will denote  
\index{model!$\cM(\vjpi)$}%
\index{zzMP@$\cM(\vjpi)$}%
the least transitive model of $\zfc'$ 
(see Definition \ref{dfex}) 
of the form $\rL_\ga$,
\kmar{cM(vjpi)}%
containing $\vjpi$ 
(and then all \muf s $\jpi_\nu$), 
in which $\la$ and all sets $\abs{\jpi_\nu}$ and 
forcings $\jpi_\nu(\xi)$ 
($\xi\in\abs{\jpi_\nu}$) 
are at most countable, 

\itlb{jfor2}% 
a \muf\ $\jpi=\bkw\vjpi=\bkw_{\nu<\la}\jpi_\nu$ 
(componentwise union)
\kmar{ bkw vjpi}%
\index{zzcupkw@$\bkw$}% 
is defined by $\abs\jpi=\bigcup_{\nu<\la}\abs{\jpi_\nu}$ and  
$\jpi(\xi)=\bigcup_{\xi<\nu<\la,\,\xi\in\abs{\jpi_\nu}}
\jpi_\nu(\xi)$ 
for all $\xi\in\abs\jpi$.
\een

\bdf
[in $\rL$]
\lam{vmf}
Let $\la\le\omi$. 
$\vmf_\la$ is the set of all  
\dd\la sequences $\vjpi=\sis{\jpi_\nu}{\nu<\la}$ of 
small \muf s $\jpi_\nu$, such that 
for each $\nu<\la$:\vom   

1) 
$\abs{\jpi_\nu}=\nu+1$,\vom 

2) 
$\jpi_{\nu}(\nu)$ contains the tree $\bse$ 
(regularity), \ 
and\vom  

3) 
$\bkw_{\mu<\nu}\jpi_\mu\rfp{\cM(\vjpi\res\nu)}\jpi_\nu$.\vom

\noi
We put $\vmf=\bigcup_{\la<\omi}\vmf_\la$. 
\edf

The set $\vmf\cup\vmo$ is ordered by the extension  
relations $\su$, $\sq$.

\ble
[in $\rL$]
\lam{112}
Assume that\/ $\ka<\la<\omi$, and\/ 
$\vjpi=\sis{\jpi_\nu}{\nu<\ka}$ 
is a sequence in\/ $\vmf_\ka$. 
Then$:$ 
\ben
\renu
\itlb{112a}% 
$\jpi=\bkw\vjpi$
is a small regular \muf\ and\/ 
$\abs\jpi=\ka$;  

\itlb{112b}% 
there is a sequence\/ $\vjqo\in\vmf$  
such that\/ $\dom\vjqo =\la$ and\/ $\vjpi\su\vjqo$.
\een
\ele
\bpf
\ref{112a}
By definition, 
$\jpi(\xi)=\bigcup_{\xi\le\nu<\ka}\jpi_\nu(\xi)$. 
The first term $\jpi_\xi(\xi)$ in the union contains 
$\bse$, so that the regularity follows.

\ref{112b}
We define \muf s $\jpi_\al$, $\ka\le\al<\la$,  
by induction on $\al$.
Assume that all terms $\jpi_\nu$, $\ka\le\nu<\al$, 
are defined, and the sequence obtained  
$\vjqo=\sis{\jpi_\mu}{\mu<\al}$ belongs to $\vmf_\al$.
Then $\jpi'=\bkw\vjqo=\bkw_{\mu<\al}\jpi_\mu$ 
is a small regular  \muf\ with $\abs{\jpi'}=\al$
by \ref{112a}, and $\jpi'\in\cM=\cM(\vjqo)$. 
Theorem \ref{tex} gives a small \muf\ $\jqo$ 
satisfying $\abs\jqo=\al$ and $\jpi'\rfa\cM\jqo$. 
Define a small \muf\ $\jpi_\al$ so that 
$\abs{\jpi_\al}=\al+1$, $\jpi_\al(\xi)=\jqo(\xi)$ 
for all$\xi<\al$, and, 
to fix the regularity, $\jpi_\al(\al)=\dpo$
(see Example \ref{xcoh}), 
hence $\bse\in\jpi_\al(\al)$.
\epf

\bdf
[key definition]
\lam{dzap}
A sequence\/ $\vjpi\in\vmf$
\rit{blocks} a set $W\sq\vmf$, if either  
\index{blocks}%
\index{blocks!positive}%
\index{blocks!negative}%
$\vjpi\in W$ (positive block)
or there is no sequence  
$\vjqo\in W$ with $\vjpi\sq\vjqo$ 
(negative block).
\edf 

Approaching the blocking sequence theorem, 
we recall that $\hc$ is the set of all   
\index{set!hereditarily countable, $\hc$}% 
\index{hereditarily countable, $\hc$}% 
\index{zzhc@$\hc$}% 
\rit{hereditarily countable\/} sets, so that  
$\hc=\rL_{\omi}$ in $\rL$.
See \cite[Part 2, Chapter 5.4]{skmlL} 
on definability classes $\is X n,\,\ip X n,\,\id X n$ 
\index{definability classes!$\is\hc n,\,\ip\hc n,\,\id\hc n$}% 
\index{zzSHC@$\is\hc n,\,\ip\hc n,\,\id\hc n$}% 
for any set $X$,
in particular, on $\is\hc n,\,\ip\hc n,\,\id\hc n$ for  
$X=\hc$ in \cite[sections 8, 9]{skmlD} 
or elsewhere.

\bte
[in $\rL$]
\lam{ep}
If\/ $\nn\ge3$ then there is a sequence\/ 
$\vdp=\sis{\gp_\al}{\al<\omi}\in\vmo$,  
satisfying the following two conditions$:$
\ben
\renu
\itlb{ep1}%
$\vdp$ itself, 
as a set of pairs\/ $\ang{\al,\gp_\al}$, 
belongs to\/ $\id\hc{\nn-1}\;;$

\itlb{ep2}%
{\rm(genericity of $\vdp$ with respect to 
$\is\hc{\nn-2}(\hc)$ sets)} 
\ if\/ $W\sq\vstf$ is a $\is\hc{\nn-2}(\hc)$ set\/ 
{\rm(\ie, parameters in $\hc$ admitted in the  
defining formula)}, 
\index{definability classes!$\is\hc n(\hc),\,\ip\hc n(\hc)$}% 
\index{zzSHChc@$\is\hc n(\hc),\,\ip\hc n(\hc)$}% 
then there is\/ $\ga<\omi$ such that  
the restricted sequence\/ 
$\vdp\res\ga=\sis{\gp_\al}{\al<\ga}\in\vmf$
blocks\/ $W$.
\een
\ete
\bpf
Let $\lel$ denote a canonical well-ordering  
of $\rL$; 
its restriction to $\hc=\rL_{\omi}$ is a  
$\id\hc1$ relation.
There exists a universal $\is\hc{\nn-2}$ set  
$\gU\sq\omi\ti\hc$. 
Thus $\gU$ belongs to $\is\hc{\nn-2}$ 
(parameter-free $\is{}{\nn-2}$ definability in $\hc$), 
and for any $\is\hc{\nn-2}(\hc)$ set $X\sq\hc$    
(definable in $\hc$ by a $\is{}{\nn-2}$ formula with  
parameters in $\hc$) 
there is an ordinal $\al<\omi$ satisfying $X=\gU_\al$, 
where $\gU_\al=\ens{x}{\ang{\al,x}\in\gU}$. 
The choice of $\omi$ as the domain of parameters 
is validated by the hypothesis $\rV=\rL$, which 
is accepted in this section and implies 
the existence of a 
$\id\hc1$ surjection $\omi\onto\hc$.

Coming back to Definition~\ref{dzap}, note that   
if $\vjpi\in\vstf$ and $W\sq\vstf$ is any set 
then there is a sequence $\vjqo\in\vstf$, 
satisfying $\vjpi\su\vjqo$ 
and blocking $W$.
We define $\vjqo_\al\in\vstf$ 
by induction on $\al<\omi$ so that $\vjqo_0=\pu$, 
$\vjqo_\la=\bigcup_{\al<\la}\vjqo_\al$ for limit $\la$, 
and each   
$\vjqo_{\al+1}$ is the \dd\lel least 
sequence $\vjqo\in\vstf$ satisfying 
$\vjpi\su\vjqo$ and blocking $\gU_\al$. 
Then $\vdp=\bigcup_{\al<\omi}\vjqo_\al\in\vsto$. 

Now \ref{ep2} holds by construction, while  
\ref{ep1} admits a routine verification based  
on the fact that $\vstf\in\id\hc1$.
\epf

\bdf
[in $\rL$]
\lam{vdp}
Fix a number $\nn\ge3$, for which   
\index{zzn@$\nn\ge 3$}%
\index{n@number $\nn\ge 3$}%
\index{number $\nn\ge 3$}%
Theorem~\ref{Tun} is proved.
\kmar{vdp, gp al}%
Fix a sequence  
\index{zzP-@$\vdp$}%
\index{zzP-al@$\gp_\al$}%
$\vdp=\sis{\gp_\al}{\al<\omi}\in\vmo$ 
which Theorem~\ref{ep} provides for this $\nn$.

We put 
$\dP=\bkw_{\al<\omi}\gp_\al$. 
\index{zzPd@$\dP$}%  
\index{forcing!zzPd@$\dP$}%  
\kmar{yp xi}%
Thus $\gp$ is a \muf, $\abs\gp=\omi$, and 
\index{multiforcing!$\yp\xi$}%  
\index{multiforcing!$\yyp\al\xi$}%  
\index{zzPxi@$\yp\xi$}%  
\index{zzPalxi@$\yyp\al\xi$}%  
$\yp\xi=\bigcup_{\xi\le\al<\omi}\yyp\al\xi$ 
for all $\xi<\omi$.
\kmar{yyp al xi}%
\index{multiforcing!$\gp_\al$}%  
By construction, each set $\gp_\al$ is a small 
\muf\ satisfying $\abs{\gp_\al}=\al+1$, while  
each component $\yyp\al\xi$ ($\xi\le\al<\omi$) 
is a countable \sfo. 
It follows that if $\al<\omi$ then the \muf\ 
$\mdp\al=\bkw_{\nu<\al}\gp_\nu$ satisfies \ 
$\abs{\mdp\al}=\al$.
%и \ $\mdp{\al}=\bkw_{\nu<\al}\gp_\nu$, 
\index{multiforcing!$\mdp\al$}%  
\index{zzP.al@$\mdp\al$}% 
In addition, since $\vdp\in\vmo$, we have
\ben
\fenu
\itlb{vdp*}
$\mdp\al\rfp{\mm\al}\gp_\al$, \ that is, \
$\mdp\al\rfa{\mm\al}\gp_\al\res\al$ \ --- \ 
for all $\al$,
\een
where $\mm\al=\cM(\vdp\res\al)$. 
\kmar{mm al, mdp al}%
\index{zzMalpha@$\mm\al$}%  
\index{model!zzMalpha@$\mm\al$}%  
The \smuf\ $\gs\al=\smf{\mdp\al}$   
\index{submultiforcing!$\gs\al$}%  
\index{zzSal@$\gs\al$}%  
\kmar{gs al}%
in $\mt{\mdp\al}$ (see Definition \ref{dfex}) 
will also be considered. 
\edf 

The set $\fP=\mt\gp$
\index{forcing!$\fP$}%  
\index{zzPf@$\fP$}%  
\kmar{fP}%  
will be used in the proof of Theorem~\ref{Tun}  
as a forcing notion. 
It is naturally identified with the countable-support  
product $\prod_{\xi<\omi}{\pc\gp\xi}$ (in $\rL$).
The sets $\gp$ and $\fP$ belong to $\rL$ by construction. 

The next theorem shows that \dd\fP generic 
extensions of $\rL$ are models for Theorem \ref{tun}. 
Therefore Theorem \ref{tuk} implies 
Theorem \ref{tun} (and Theorem \ref{Tun} as well). 

\bte
\lam{tuk}
Under the conditions of Definition~\ref{vdp}, 
let\/ $\zG\sq\fP$ be a generic filter over\/ $\rL$. 
Then the following holds in\/ $\rL[\zG]$$:$
\ben
\renu
\itlb{tuk1}% 
condition\/ \ref{tun1} of Theorem\/ \ref{tun}$;$
%существует\/ \dd{\ip1\nn}множество\/ $W\sq\dn\ti\dn$, 
%все вертикальные сечения\/ $W_x=\ens{y}{\ang{x,y}\in P}$ 
%которого являются\/ \dd\Eo классами, 
%и которое нельзя униформизовать никаким ROD множеством$;$

\itlb{tuk2}% 
condition\/ \ref{tun2} of Theorem\/ \ref{tun}$.$
%если\/ $p\in\dn$, то 
%каждое\/ \dd{\is1\nn(p)}множество\/ $W\sq\dn\ti\dn$, 
%все вертикальные сечения которого не более чем счетны, 
%униформизуется множеством класса\/ $\id1{\nn+1}(p)$, 
%в частности, ROD множеством. 
\een
\ete

To prove Theorem \ref{tuk}, 
we explore properties of the forcing notion $\fP$ 
and related generic extensions in Sections 
\ref{pres+}--\ref{bex2}, 
then establish \ref{tuk1} of Theorem \ref{tuk}  
in Section \ref{nounc}, and finally \ref{tuk2} 
in Section \ref{zav} with the help of a special  
approximating forcing relation $\fo$.

\parf{Key forcing properties} 
\las{pres+}

Here we study $\fP$ as the forcing notion.
{\ubf We argue under the conditions and 
notation of Definition~\ref{vdp}}.

\bdf
[in $\rL$]
\lam{ppro}
\sloppy
If $C\sq\omi$ then we define the subproduct  
\index{forcing!$\fP\res C$}%  
\index{forcing!restricted, $\fP\res C$}%  
\index{zzPfС@$\fP\res C$}%  
$\fP\res C=\mt{\gp\res C}=\ens{\zp\in\fP}{\abs\zp\sq C}=
\prod_{\xi\in C}\yp\xi$ 
with countable support. 
Then $\fP$ can be identified with  
$(\fP\res{C})\ti\big(\fP\res{(\omil\bez C)}\big)$.

If $C\sq\omi$ is at most countable (in $\rL$), then  
by the regularity of $\dP$ the set  
$\fP\res C$ can be identified with  
\index{forcing!$\fP_C$}%  
\index{forcing!restricted, $\fP_C$}%  
\index{zzPf-С@$\fP_C$}%  
$\fP_C=\ens{\zp\in\fP}{\abs\zp=C}$.  

If $C=\ans\xi$, $\xi<\omil$, then $\fP\res\ans\xi$ 
is naturally identified with $\yp\xi$, and then 
$\fP$ is identified with  
$\yp\xi\ti\fP\res{\cne\xi}$, where  
$\cne\xi=\omil\bez\ans{\xi}$.
\kmar{cne xi}%
\index{zzC=xi@$\cne\xi$}%
\edf

\ble
\lam{pqa}
If\/ $\xi\le\al<\ga<\omi$ then\/ 
$\yyp\al\xi\rf\yyp\ga\xi$ in the sense of\/ \ref{fm}. 
Therefore each\/ $\yyp\al\xi$ is pre-dense in\/ 
$\yp\xi=\bigcup_{\al\ge\xi}\yyp\al\xi$ 
by Lemma \ref{pqr'}\ref{pqr'6}.
\ele
\bpf
\sloppy
Arguing by induction, suppose that   
$\yyp\mu\xi \rf \yyp\nu\xi$ 
is established for all $\xi\le\mu<\nu<\ga$. 
Lemma \ref{pqr'}\ref{pqr'6} implies that the set  
$\yyp\al\xi$  
is pre-dense in $\bigcup_{\xi\le\nu<\ga}\yyp\nu\xi$. 
The \muf\ $\jqo=\gp_\ga\res\ga$ satisfies 
$\mdp\ga\rfa{\cM_\ga}\jqo$ by \ref{vdp}\ref{vdp*}. 
By Definition \ref{dfex}, this includes the condition   
$\mdp\ga(\xi)\rfa{\cM_\ga}\jqo(\xi)$.  
Then immediately  $\jqo(\xi)$ 
is dense in $\mdp\ga(\xi)\cup\jqo(\xi)$. 
However $\jqo(\xi)=\gp_\ga(\xi)$ while 
$\mdp\ga(\xi)= \bigcup_{\xi\le\nu<\ga}\gp_\nu(\xi)$. 
Therefore, first, $\yyp\ga\xi$ is dense in  
$\yyp\al\xi\cup\yyp\ga\xi$, thus we have \ref{fm1} 
of Definition \ref{fm}.
And second, as the set $\yyp\al\xi$ is dense in 
$\mdp\ga(\xi)$ by the above,
and clearly $\yyp\al\xi\in\cM_\ga$, we obtain  
$S\sqfv\yyp\al\xi$ for each tree 
$S\in\jqo(\xi)=\yyp\ga\xi$, thus we have \ref{fm2} 
of Definition \ref{fm}. 
%Итак, $\gp_\al(\xi)\rf \gp_\ga(\xi)$ установлено. 
\epf

\ble
[in $\rL$]
\lam{wcc}
Assume that, for each\/ $n$, 
$\zD_n\sq\fP$ is open dense in\/ $\fP$, 
and let\/ $\zp\in\fP$. 
There is a \mut\/ $\zq\in\fP$, satisfying\/ 
$\zq\leq\zp$ and\/ $\zq\sqdv \zD_n$ for all\/ $n$.
Therefore\/ \dd\fP generic extensions of\/ $\rL$
preserve $\omil$.
\ele
\bpf
There is a countable elementary submodel $M$ 
of $\stk{\rL_{\om_2}}{\in}$,  
containing $\zp$ and all sets $\zD_n$. 
Then $M$ also contains $\omi$, 
as it is a definable set, and contains the  
sequence $\vdp$ along with the derived  
sets $\gp=\bkw\vdp$, $\fP=\mt\gp$, 
by the same reason.
The set $M\cap\rL_{\omi}$ is transitive. 
Indeed if $X\in M\cap\rL_{\omi}$ then $X$ is at most 
countable, hence there exist functions $f:\om\onto X$. 
Let $f_X$ be the least of them in the sense of  
the G\"odel wellordering $\lel$ of $\rL$. 
Then $f_X\in M$ since $X\in M$ and the ordering             
${\lel}\res {\rL_{\om_2}}$ is definable in ${\rL_{\om_2}}$. 
It follows that each $x\in X$ belongs to 
$M$ because $x=f_X(k)$ for some $k$.

Let $\phi:M\onto\rL_\la$ be the Mostowski collapse function, 
and $\al=\phi(\omi)$. 
Then $\al<\la<\omi$ and, by the transitivity, 
it holds (*) $\phi(x)=x$ for all $x\in M\cap\rL_{\omi}$.
Thus $\phi(\xi)=\xi$, $\phi(T)=T$, $\phi(\zq)=\zq$ 
for each ordinal $\xi\in M\cap\omi$, 
tree $T\in M\cap\lt$, \mut\ $\zq\in M\cap\mtt$. 
We conclude that
$\phi(\vdp)=\vdp\cap\rL_\al=\vdp\res\al$, 
$\phi(\gp)=\mdp\al=\bkw_{\ga<\al}\gp_\al$ 
(a \muf\ with $\abs{\mdp\al}=\al$), and 
%и наконец, благодаря простому рассуждению, 
$\phi(\fP)=\fP\cap\rL_\al
=\mt{\mdp\al}\cap\rL_\al$. 

We assert that moreover
$\phi(\fP)=\mP_\al$, where, we recall, 
$\mP_\al=\smf{\mdp\al}$.
Indeed by Definition \ref{dfex}, 
$\smf{\mdp\al}$ is equal to the closure of 
$\mt{\mdp\al}\cap\rL_\al$ relatively to the three 
operations of Definition~\ref{mfd}.
But $\vpi(\fP=\mt{\mdp\al}\cap\rL_\al$, 
thus 
$\mt{\mdp\al}\cap\rL_\al$ is already closed  
under the operations, 
since so is $\fP=\mt{\gp}$. 
We conclude that $\smf{\mdp\al}=\mt{\mdp\al}\cap\rL_\al$.

Furthermore, a similar argument allows to prove 
that if $n<\om$ then the set
$\phi(\zD_n)=\zD_n\cap\rL_\al=\zD_n\cap\mP_\al\in\rL_\la$
is open dense in $\smf{\mdp\al}$. 
In addition, $\phi(\zp)=\zp\in\mP_\al$.
On the other hand, by the elementarity, 
the ordinal $\al$ is uncountable in $\rL_\la$.  
It follows that $\rL_\la\sq\mm\al$.
However we have $\mdp\al\rfa{\mm\al}\gp_\al\res\al$ 
by \ref{vdp}\ref{vdp*}, and also 
$\zp\in\mP_\al=\smf{\mdp\al}$.
Therefore, by Definition~\ref{dfex}\ref{dfex5}, 
there exists a \mut\ $\zq\in\mt{\gp_\al}$ 
satisfying  $\zq\leq\zp$ and $\zq\sqdv \phi(\zD_n)$ 
for all $n$ 
Finally, $\mt{\gp_\al}\sq\fP$ and $\phi(\zD_n)\sq \zD_n$. 
This ends the proof of the first claim.

To prove the second claim of the lemma, suppose 
towards the contrary that  
$\dof$ is a name of a function from $\om$ to $\omil$, 
and some $\zp\in\fP$ forces $\ran\dof=\omil$. 
Let $\zD_{n\al}$ be the set of all \mut s 
$\zr\in\fP$, either 
(1) incompatible with $\zp$ in $\fP$, 
or 
(2) satisfying $\zr\leq\zp$ and \dd\fP forcing 
$\dof(n)=\al$. 
A simple argument shows that every set  
$\zD_n=\bigcup_\al\zD_{n\al}$ is dense in $\fP$. 
By the first claim of the lemma, there exists 
a \mut\ $\zq\in\fP$ satisfying
$\zq\leq\zp$ and $\zq\sqdv\zD_n$, $\kaz n$. 
Let the relations $\zq\sqdv\zD_n$ be witnessed 
by finite sets $\zD'_n\sq\zD_n$. 
Accordingly, the sets 
$A_n=\ens{\al}{\zD'_n\cap\zD_{n\al}\ne\pu}$
are finite, hence the union  
$A=\bigcup_nA_n$ is countable in $\rL$, 
\ie, $\omil\not\sq A$. 
On the other hand, we assert that  
$\zq$ forces $\dof(n)\in A_n$, for each $n$. 
This implies a contradiction and accomplishes the proof.

To finally prove that $\zq$ forces $\dof(n)\in A_n$, 
suppose to the contrary that $\zr\in\fP$, 
$\zr\leq\zq$, and $\zr$ forces $\dof(n)=\al$, 
where $\al<\omil$, $\al\nin A_n$. 
Then $\zr\sqdv\zD_n$ by means of the same   
finite set $\zD'_n\sq\zD_n$.
%Это означает, что каждое $\jU\in \zD'_n$ удовлетворяет 
Lemma \ref{tu}\ref{tu4} provides a string  
$\sg\in\bse$ and a \mut\ $\ju\in \zD'_n$ such that 
$\zr'=\rau\zr\sg\leq\ju$.  
Note that $\zr'\in\fP$ by Lemma \ref{mokr}\ref{mokr1}. 
Thus the \mut s $\zr$ and $\ju$ are compatible in $\fP$. 
Finally $\ju\in \zD'_n\sq\zD_n$, therefore  
$\ju\in\zD_{n\ga}$ for some $\ga$. 
Then by definition $\ju$ forces $\dof(n)=\ga$, where 
$\ga\in A_n$, that is, $\ga\ne\al$. 
However $\zr$ forces $\dof(n)=\al$, where $\al\nin A_n$,   
which is a contradiction.
\epf

\ble
[in $\rL$]
\lam{blo}
If a set of \mut s\/ $Q\sq\mtt$ belongs to\/ 
$\is\hc{\nn-2}(\hc)$ and\/ 
$Q^-=\ens{\zp\in\mtt}{\neg\:\sus\zq\in Q\,(S\leq T)}$,  
then the set\/ $\fP\cap{(Q\cup Q^-)}$ 
is dense in\/ $\fP$.
In particular if\/ $Q$ is dense in\/ $\mtt$ then\/ 
$Q\cap\fP$ is dense in\/ $\fP$.
\ele
\bpf
Consider a \mut\ $\zp_0\in\fP=\mt\gp$, thus 
$\zp_0\in\mt{\mdp{\al_0}}$, $\al_0<\omi$. 
The set $\Da$ of all sequences  
$\vjpi\in\vmf$, such that $\vdp\res\al_0\sq\vjpi$ and 
$\sus\zp\in Q\cap(\mt{\bkw\vjpi})\,(\zp\leq\zp_0)$, 
belongs to $\is\hc{\nn-2}(\hc)$ 
as so does $Q$. 
Therefore there is an ordinal $\al<\omi$   
such that the sequence $\vdp\res\al$ 
blocks $\Da$.\vom 

{\ubf Case 1:\/} 
$\vdp\res\al\in\Da$; let this be witnessed by 
$\zp\in Q\cap(\mt{\bkw(\vdp\res\al)})$.
Then $\al_0\le\al$ and  
the \mut\ $\zp$ belongs to   
$Q\cap\fP$ and satisfies $\zp\leq\zp_0$.\vom

{\ubf Case 2:\/} 
no sequence in $\Da$ extends $\vdp\res\al$.
Let $\ga=\tmax\ans{\al,\al_0}$. 
% and $\cM=\cM(\vdp\res\ga)$.
Then $\mdp\ga\prol{\cM_\ga}\gp_\ga\res\ga$ by 
\ref{vdp}\ref{vdp*}. 
As $\al_0\le\ga$, there exists a \mut\ $\zp\in\mt{\gp_\ga}$, 
$\zp\leq\zp_0$. 
We can wlog assume that $\abs\zp=\abs{\gp_\ga}$, 
that is $\,= \ga+1$.
Then $\zp(\xi)\in\gp_\ga(\xi)$ for all $\xi\le\ga$.
It remains to prove that $\zp\in Q^-$. 

Suppose to the contrary that $\zp\nin Q^-$. 
By definition there is a \mut\ $\zq\in Q$, $\zq\leq\zp$.
Then $\ga+1=\abs\zp\sq\abs\zq$. 
We can assume that  
$\abs\zq=\la<\omi$, $\la\ge\ga+1$.
We are going to define a sequence  
$\vjpi=\sis{\jpi_\al}{\al<\la}\in\vmf$, which extends  
$\vdp\res\ga$, that is, 
$\jpi_\al=\gp_\al$ for all $\al<\ga$,
and satisfies $\zq\in\mt{\bkw\vjpi}$. 
This implies $\vjpi\in\Da$ by the choice of $\zq$, 
which contradicts to the Case 2 hypothesis and  
completes the proof of $\zp\in Q^-$ and the proof 
of the lemma. 

Thus we have to appropriately define \muf s 
$\jpi_\al$, $\ga\le\al<\la$. 
We begin with $\jpi_\ga$. 
This is based on the \muf\ $\gp_\ga$. 
Note that $\zq(\xi)\sq\zp(\xi)\in\gp_\ga(\xi)$ 
for all $\xi\le\ga$.
We put 
$\jpi_\ga(\xi)=\gp_\ga(\xi)\cup
\ens{\sg\app (\raw{\zq(\xi)}t)}{t,\sg\in\bse}$ 
for all $\xi\le\ga$.
Every \lap{new} tree 
$S=\sg\ap{(\raw{\zq(\xi)}t)}$ satisfies 
$S\sq \sg\ap \zp(\xi)$, 
where $\sg\ap \zp(\xi)\in\gp_\ga(\xi)$. 
However  
$\mdp\ga\rfp{\mm\ga}\gp_\ga$ 
by Definition \ref{vdp}\ref{vdp*}. 
Therefore $\mdp\ga\rfp{\mm\ga}\jpi_\ga$ as well. 
Thus the term $\jpi_\ga$ extends the \mus\   
$\vdp\res\ga=\sis{\gp_\al}{\al<\ga}
=\sis{\jpi_\al}{\al<\ga}\in\vmf_\ga$ 
to a \mus\ in   
$\vmf_{\ga+1}$, and we have $\zq(\xi)\in\jpi_\ga(\xi)$ 
for all $\xi\le\ga$.
The extended \mus\ can be further extended to a \mus\  
in $\vmf_{\la}$ by terms $\jpi_\al$, $\ga<\al<\la$, 
by induction as in the proof of \ref{112}\ref{112b}, 
with the amendment that  
$\jpi_\al(\al)=
\dpo\cup\ens{\sg\app(\raw{\zq(\al)}t)}{t,\sg\in\bse}$, 
rather than just $\jpi_\al(\al)=\dpo$, 
for all $\al$.
\epf

\parf{Generic  extension} 
\las{bex1}

Here we study \dd\fP generic extensions $\rL[\zG]$ 
of $\rL$ obtained by adjoining \dd\fP generic sets 
$\zG\sq\fP$ to $\rL$. 
We will use the forcing notion $\fP=\mt\gp\in\rL$ 
and other notation of Definition~\ref{vdp}, 
with the difference that the reasoning will 
not be relativized to $\rL$ by default, 
and accordingly the first uncountable cardinal in $\rL$ 
will be denoted by $\omil$ instead of $\omi$.

\bdf
[generic reals]
\lam{gg}
Let a set $\zG\sq\fP$ be \dd\fP generic over $\rL$.
Note that $\omi^{\rL[\zG]}=\omil$ 
by Lemma~\ref{wcc}.

If $\xi<\omil$ then  
$\zG(\xi)=\ens{\zc\zp\xi}{\xi\in\abs\zp\land\zp\in \zG}$ is 
\index{zzGx@$\zG(\xi)$}%
a set   
\dd{\yp\xi}generic over $\rL$, the intersection 
\kmar{yp xi}% 
$X_{\xi}=\bigcap_{T\in \zG(\xi)}[T]$ contains a single 
real $\bx_{\xi}=\bx_{\xi}[\zG]\in\dn$, 
and this real  
\index{zzxxiG@$x_{\xi}[\zG]$}%
is \dd{\yp\xi}generic over $\rL$.
These reals are assembled into a \lap{multireal} 
$\xx \zG=\sis{\bx_{\xi}[\zG]}{\xi<\omil}\in
\ddn\omil$.
\edf

\bcor
[of \ref{ppro} and the product forcing theorem]
\lam{ggaa}
If \/ $B\in\rL$, $B\sq\omil$ is at most countable in\/ $\rL$, 
and\/  
$\zG\sq\fP$ is \dd\fP generic over\/ $\rL$, 
then the set\/ $\zG_B=\ens{\zp\in\zG}{\abs\zp=B}$ 
is\/ \dd{\fP_B}generic over\/ $\rL$.\qed
\ecor

Recall that $\cne\xi=\omil\bez\ans\xi$.

\bpro
[in terms of Definition~\ref{gg}] 
\lam{sym}
If\/ $\xi<\omil$ then the real\/ $\bx_\xi[\zG]$ 
is not\/ 
\dd{(\ans{\zG\res\cne\xi}\cup\Ord)}definable 
in\/ $\rL[\zG]$, in particular, 
$\bx_\xi[\zG]\nin\rL[\zG\res\cne\xi]$. 
\epro
\bpf
See the proof of Lemma 14.5 in \cite{kl30},  
based on the product forcing theorem and the 
\dd\Eo invariance of each component 
$\yp\xi$ in the sense of \ref{sitf}\ref{sitf3}.
\epf

The next theorem belongs to the type of  
\lap{continuous reading of names} theorems
in theory of forcing extensions. 
It involves the coding of continuous maps by   
Definition \ref{nfk}, and asserts that 
reals $x\in \dn$ in  
\dd\fP generic extensions are obtained by applications 
of continuous maps coded in $\rL$ to suitable sequences 
of generic reals.
To render the notation less cumbersome, 
if $\bc\in\rL$ and 
$\bc\in\kng$ in $\rL$, and $\zG\sq\fP$ is   
generic over $\rL$, 
then we put $\fg\bc:=f^\bc(\xx G \res B)$, 
where $B=\abs\bc$.

\ble
\lam{repd}
If\/ $C\in\rL$, $C\sq\omil$, 
$\zG\sq\fP$ is generic\/ $\rL$, 
and\/ $x\in\dn\cap\rL[\zG\res C]$, 
then there is a code\/ $\bc\in\kng\cap\rL$, 
%где\/ $B=\abs\zp$, 
such that\/ $\abs\bc\sq C$ and\/ $x=\fg\bc$.
\ele

\bpf
Let $\dox$ be a name for $x$ in the forcing language  
related to the forcing notion $\fP$. 
Thus the indexed family of sets
$$
A_{ki}=\ens{\zp\in\fP}
{\zp\,\text{ forces that }\,\dox(k)=i}
\,,\quad k<\om,\;\,i=0,1,
$$
belongs to $\rL$ and we have \ 
(A) $x(k)=i\eqv \zG\cap A_{ki}\ne\pu$, \ 
(B) $A_{k0}\cap A_{k1}=\pu$, \ and 
(C) each set $A_k=A_{k0}\cup A_{k1}$
is open dense in $\fP$.
We can assume that $\dox$ contains  
an explicit effective construction of $x$ 
from $\zG\res C$, and then: \ 
(*) if $\zq\in A_{ki}$ then 
$\zq\res{(C\cap\abs\zq)}\in A_{ki}$ as well. 

The set  
$D=\ens{\zp\in\fP}{\kaz k\,(\zp\sqdv A_k)}$ 
also is dense in $\fP$ by Lemma \ref{wcc}. 
Therefore, by the genericity, there is a 
\mut\ $\zp'\in \zG$, such that 
$\zp'\sqdv A_k$ for all $k$. 
In addition, (*) implies that the \mut\ 
$\zp=\zp'\res{(C\cap\abs\zp)}\in \zG$ also satisfies  
$\zp\sqdv A_k$, $\kaz k$, but now $\abs\zp\sq C$. 

This means (Definition \ref{sqdv}) that, 
in $\rL$, there exists a sequence  
of finite sets $F_k\sq A_k$ which witness $\zp\sqdv A_k$ 
in the sense that: \ 
1) $\abc \ju\sq B=\abc\zp$ 
for all $\ju\in F_k$, \  
2) $[\zp]\sq\bigcup_{\ju\in F_k}[\ju\ares B]$, \ and
3) $[\ju\ares B]\cap[\jv\ares B]=\pu$ for all 
$\jv\ne\ju$ in $F_k$. 
We put $F_{ki}=F_k\cap A_{ki}$, $i=0,1$.

Now arguing in $\rL$ we define a continuous 
$f:[\zp]\to\dn$ as follows: 
$f(y)(k)=i$, if there is a \mut\ $\zq\in F_{ki}$ 
with $y\res{\abs \zq}\in[\zq]$. 
Then $f=f^\bc\res[\zq]$ by Lemma \ref{knf}, 
where $\bc$ is a suitable code $\knf B\cap\rL$. 
One easily verifies that $x=\fg\bc$
\epf

By the next theorem, the relation $y\nin\rL[x]$ 
between reals $x,y\in\dn$ in $\rL[\zG]$ is fully 
determined by a generic real $\bx_\xi[\zG]$, 
so that $\bx_\xi[\zG]$ belongs to $\rL[y]$ but 
$x$ belongs to $\rL[\xx \zG\res \cne\xi]$, where 
definitely $\bx_\xi[\zG]\nin\rL[\xx \zG\res \cne\xi]$ 
by Proposition~\ref{sym}.

\bte
\lam{fg}
If a set\/ $\zG\sq\fP$ is generic over\/ $\rL$, 
%$C\in\rL$, $C\sq\omil$, 
$x,y\in\dn\cap\rL[\zG]$, 
%\res C]$,  
and\/ $y\nin\rL[x]$,
then there is an ordinal\/ $\xi<\omil$ such that\/ 
$x\in\rL[\xx \zG\res \cne\xi]$ but\/ $\bx_\xi[\zG]\in\rL[y]$, 
and in addition, $\bx_\xi[\zG]=g(y)$, where\/ $g:\dn\to\dn$ 
is a continuous map coded in\/ $\rL$.
\ete  

\bpf
By Lemma \ref{repd}, there exist codes
$\bc,\bd\in\kng\cap\rL$, 
such that $x=\fg\bc$ and $y=\fg\bd$.
Let $B=\abs\bc\cup\abs\bd$; 
we can assume that $\abs\bc=\abs\bd=B$.

\rit{We argue in\/ $\rL$}. 
The set $\zD$ of all \mut s $\zq\in\mlt_B$ such that 
either\/ 
{\rm(i)} 
$f^\bd$ is reduced to\/ $f^\bc$ on\/ $[\zq]$, 
or\/ 
{\rm(ii)} 
$f^\bd$ captures some ordinal\/ $\xi\in B$ on\/ $[\zq]$  
and\/ $f^\bc$ is reduced to the set\/ $B\bez\ans\xi$ 
on\/ $[\zq]$, is dense in $\mlt_B$ by Theorem \ref{nwft}. 
% ниже $\zp$ 
%(\te\ плотно в множестве $\ens{\zp'\in\mlt_B}{\zp'\leq\zp}$).  
It follows by Lemma \ref{blo} that the set  
$\zD'=\zD\cap{\fP_B}$ is dense in 
$\fP_B=\ens{\zr\in\fP}{\abs\zr=B}$. 
% ниже $\zp$. 

\rit{We argue in\/  $\rL[\zG]$}. 
We have $\zG\cap\zD\ne\pu$ by Corollary \ref{ggaa}. 
%(так как $\zp\in \zG$). 
Let $\zq\in \zG\cap\zD$; then $\xx \zG\res B\in[\zq]$. 
Note that (i) fails for this $\zq$, 
since (i) implies $f^\bd(z)=g(f^\bc(z))$ 
for all $z\in[\zq]$, where $g:\dn\to\dn$ is a continuous  
map coded in $\rL$, thus (with $z=\xx \zG\res B$) we get 
$y=g(x)$, and further $y\in\rL[x]$ 
(as $g$ is coded in $\rL$), 
a contradiction to the condition of the theorem. 
Thus (ii) holds, \ie, \rit{still in $\rL$}, 
$f^\bd$ captures an ordinal $\xi\in B$ on $[\zq]$, while  
$f^\bc$ is reduced to $B\bez\ans\xi$ on $[\zq]$. 

By the compactness of the spaces considered, this  
implies the existence of continuous maps  
$f:\ddn{B\bez\ans\xi}\to\dn$ and $g:\dn\to\dn$, 
both coded in $\rL$ and satisfying  
$f^\bc(z)=f(z\res{(B\bez\ans\xi)})$ and 
$z(\xi)=g(f^\bd(z))$
for all $z\in\zq$. 
In particular, for $z=\xx \zG\res B$, we have  
$x=f^\bc(\xx \zG\res{(B\bez\ans\xi)})$, hence  
$x\in\rL[\xx \zG\res \cne\xi]$, 
and $\bx_\xi[\zG]=g(y)$, hence 
$\bx_\xi[\zG]\in\rL[y]$.
\epf

\vyk{
Перед формулировкой следствия, если $C\in\rL$, $C\sq\omi$, 
то $\is1n(\xx \zG\res C)$ будет обозначать класс всех множеств, 
определимых \dd{\is1n}формулами с параметрами из 
$\dn\cap\rL[\xx \zG\res C]$. 

\bcor
\lam{fgc}
В условиях теоремы \ref{fg}, 
если\/ $y$ принадлежит счетному множеству\/ $Y\sq\dn$ 
класса\/ $\is1n(\xx \zG\res \cne\xi)$ то и\/ 
$\bx_\xi[\zG]$ принадлежит счетному множеству\/ $X\sq\dn$ 
класса\/ $\is1n(\xx \zG\res \cne\xi)$.
\ecor
\bpf
$X$ есть \dd fобраз множества $Y$, где, как в теореме, 
$f$ --- непрерывная функция с кодом из $\rL$.
\epf
}

\parf{Definability of generic reals} 
\las{bex2}

We continue to argue in terms of definitions   
\ref{vdp} and \ref{gg}. 
Now the main goal will be to study \dd{\yp\xi}generic  
reals $x\in\dn$ in \dd\fP extensions of $\rL$.
 
\bte
\lam{X=Y}
In any\/ \dd\fP generic extension\/ $\rL[\zG]$ 
of\/ $\rL$, 
it is true that$:$ 
if\/ $\xi<\omil$ then the set\/ 
$X_\xi=\eko{\bx_\xi[\zG]}=\ens{\sg\ap\bx_\xi[\zG]}{\sg\in\bse}$ 
is equal to the set\/ 
$Y_\xi=
\bigcap_{\xi\le\al<\omil}\bigcup_{U\in\yyp{\al}\xi}[U]$. 
\ete
\bpf
The real $x=\bx_\xi[\zG]\in\dn$ is  
\dd{\yp\xi}generic, while every set of the form 
$\yyp{\al}\xi$ is pre-dense in $\yp\xi$ by Lemma \ref{pqa}.
Therefore $x\in Y_\xi$. 
Moreover all sets $\yyp{\al}\xi$ are \sfo s by construction, 
hence they are \dd\Eo invariant in the sense of 
\ref{sitf}\ref{sitf3}. 
It follows that $X_\xi\sq Y_\xi$. 

To establish the inverse, 
assume that $y_0\in Y_\xi$ in $\rL[\zG]$. 
By Lemma \ref{repd}, there is a code 
$\bc\in\kng\cap\rL$ such that 
$y_0=\fg\bc=f^\bc(\xx \zG\res B)$, where $B=\abs\bc$.
Consider the set $\zD$ of all  \mut s 
$\zq\in\fP_B$ such that either\/ 
{\rm(i)} 
there is a string $\sg\in\bse$ such that 
$f^\bc(x)=\sg\ap x(\xi)$ for all $x\in[\zq]$,  
or  
{\rm(ii)} 
there exists an ordinal $\al$, $\xi\le\al<\omi$, 
such that $f^\bc(x)\nin\bigcup_{U\in\yyp{\al}\xi}[U]$
for all $x\in[\zq]$. 

\ble
\lam{XY*}%
The set\/ $\zD$ is dense in\/ $\fP_B$. 
\ele
\bpf
Let $\zp\in\fP_B$,  
then $\abs{\zp}=B$. 
%and $\zp'\leq\zp$. 
There is an ordinal $\al<\omil$ such that 
1) $B\sq\al$ --- hence $\xi<\al$,
2) $\zp\in\mP_\al=\smf{\mdp\al}$, and 
3) $\bc\in\mm\al$.  
Note that 
$\mdp\al\rfp{\mm\al}\gp_\al$ 
holds by \ref{vdp}\ref{vdp*}.
Therefore by Definition \ref{dfex}\ref{dfex6} there
is a \mut\ $\zq\in\mt{\gp_\al}$ such that  
$\abs\zq=\abs{\zp}=B$, $\zq\leq\zp$, 
and either (i) 
there is a string $\sg\in\bse$ 
satisfying $f^\bc(x)=\sg\ap x(\xi)$ for all $x\in[\zq]$,  
or (ii) 
$f^\bc(x)\nin\bigcup_{U\in\yyp{\al}\xi}[U]$ 
for all $x\in[\zq]$.
Thus $\zq\in\zD$, getting the density.
\epF{Lemma} 

We return to the theorem. 
Corollary \ref{ggaa} implies $\zG\cap\zD\ne\pu$ 
by the lemma. 
%(так как $\zp\in \zG$). 
Let $\zq\in \zG\cap\zD$. 
In particular $x_0=\xx \zG\res B\in[\zq]$. 
It follows that $\zq$ does not satisfy  
(ii) of the definition of $\zD$, 
since $y_0=f^\bc(x_0)\in Y_\xi$. 
Therefore $\zq$ satisfies (i) of the definition of $\zD$ 
with some $\sg\in\bse.$
Then  
$y_0=f^\bc(x_0)= \sg\ap x_0(\xi)= \sg\ap \xx \zG(\xi)
=\sg\ap\bx_\xi[\zG]$, 
that is, $y_0\in X$, as required.
\epf       

One easily proves that, under the conditions of 
the theorem, the set 
$X_\xi=Y_\xi$ is equal to the set of all  
\dd{\yp\xi}generic reals $y\in\dn$, see 
\cite{kl22}.

\parf{Non-uniformizable set} 
\las{nounc}

Here we prove claim \ref{tuk1} of Theorem \ref{tuk}.
To begin with, we define a non-uniformizable set in 
the \lap{rectangle} $\omil\ti\dn$.

\ble
\lam{lnu}
%Пусть множество\/ $\zG\sq\fP$ является\/ 
%\dd\fP генерическим над\/ $\rL$. 
%Тогда 
Under the conditions of Theorem \ref{tuk}, the set\/ 
$
\gfu=
%\gfu[\zG]=
\ens{\ang{\xi,x}}{\xi<\omil\land x\in\eko{\bx_\xi[\zG]}}
$ 
\index{set!$\gfu=\gfu[\zG]$}%
\index{zzKKG@$\gfu=\gfu[\zG]$}%
belongs to\/ $\rL[\zG]$ and has 
the following properties in\/ $\rL[\zG]:$
\ben
\renu
\itlb{lnu1}\msur%
$\gfu$ belongs to the definability class\/ $\ip\hc{\nn-1}\;;$ 

\itlb{lnu3}%
if\/ $\xi<\omi$ then the cross-section\/ 
$\gfu_\xi=\ens{x}{\ang{\xi,x}\in\gfu}$  
is a\/ \dd\Eo class$;$

\itlb{lnu4}%
the set\/ $\gfu$ is not  
ROD-uniformizable.
\een
\ele
\bpf
%\ref{lnu2} выполнено по лемме~\ref{only} а 
\ref{lnu3} is quite obvious: $\gfu_\xi=\eko{x_\xi[\zG]}$.
To prove \ref{lnu1} we note that 
Lemma\ref{wcc} implies $\omi=\omil$ in $\rL[\zG]$. 
Therefore by Theorem~\ref{X=Y}, 
the sentence $\ang{\xi,x}\in\gfu$ 
is equivalent to
%тому, что 
$$
\xi<\omi\;\land\;
\kaz\al\:\big(\xi\le\al<\omi\imp\sus T\in \yyp{\al}\xi
\,(x\in[T])\big)\,.
$$
%по лемме~\ref{mod1}.
Yet the formula in the outer brackets here  
expresses a $\ip\hc{\nn-1}$ relation 
by condition~\ref{ep1} of Theorem \ref{ep}. 
(The quantifier $\sus T\in \yyp{\al}\xi$ is bounded, 
hence it does not affect the definability estimation.)

To prove \ref{lnu4} suppose towards the contrary 
that it is true in $\rL[\zG]$ that  
$R\sq\gfu$ is an uniformizing ROD set. 
Let $r\in\dn\cap\rL[\zG]$ be a real which witnesses 
that $R$ is $\ans r\cup\Ord$-definable in $\rL[\zG]$.
Lemma~\ref{wcc} (preservation of $\omil$) implies  
the existence of an ordinal $\xi<\omil$ such that 
$r\in \rL[\zG\res\ens{\et}{\et<\xi}]$, 
hence $r\in \rL[\zG\res\cne\xi]$, 
where $\cne\xi=\omil\bez\ans\xi$.
Therefore the unique real $x\in\dn$, satisfying  
$\ang{\xi,x}\in R$, 
is \dd{(\ans{\zG\res\cne\xi}\cup\Ord)}definable 
in $\rL[\zG]$.
However $R\sq\gfu$, thus 
%$x$ является 
%\dd{\opp\etz}генерической точкой согласно \ref{lnu2}, 
%и мы имеем 
$x\Eo \bx_\xi[\zG]$. 
% по лемме \ref{only}. 
It follows that the generic real $\bx_\xi[\zG]$ itself is
\dd{(\ans{\zG\res\cne\xi}\cup\Ord)}definable 
in $\rL[\zG]$.
But this contradicts to Proposition \ref{sym}.
\epf

To transform the set $\gfu=\gfu[\zG]$ into a similar 
non-uniformizable set in the plane $\dn\ti\dn$,  
we make use of the following rather elementary 
transformation. 
% в общем, элементарное преобразование,   
%не связанное с форсингом и моделями.

Let   
$\raz=\ens{q_n}{n<\om}$ be a recursive enumeration 
of the rationals. 
If $z\in\dn$ then let 
$Q_z=\ens{q_n}{z(n)=1}\sq\raz$, 
$Q'_z\sq Q_z$ be the largest (perhaps, empty) 
well-ordered initial segment of $Q_z$, and let 
$\abs z<\omi$ be the ordinal number of $Q'_z$; 
thus obviously $\ens{\abs z}{z\in\dn}=\omi$.

\ble
\lam{lnv}
Under the conditions of Theorem \ref{tuk}, the set\/ 
$$
\ww=\ens{\ang{z,x}\in\break 
\dn\ti\dn}{\ang{\abs z,x}\in \gfu}
$$ 
\index{set!$\ww$}%
\index{zzW@$\ww$}%
belongs to\/ $\rL[\zG]$ and has 
%Тогда в\/ $\rL[\zG]$ истинно, что множество\/ 
%$\gfu=\gfu[\zG]$  
%из\/ \ref{mod2} 
the following properties in\/ $\rL[\zG]:$
\ben
\renu
\itlb{lnv1}\msur%
$\ww$ belongs to the definability class\/ $\ip1{\nn}\;;$  

\itlb{lnv3}%
if\/ $z\in\dn$ then the cross-section\/ 
$\ww_z=\ens{x}{\ang{z,x}\in\ww}$  
is a\/ \dd\Eo class$;$

\itlb{lnv4}%
the set\/ $\ww$ is not  
ROD-uniformizable.
\een
\ele
\bpf
The set $\ww$ belongs to $\ip\hc{\nn-1}$ since so  
does $\gfu$; 
indeed the map $z\mto\abs z$ is  
$\id\hc1$. 
Therefore by the transfer theorem  
(see \eg\ 9.1 in \cite{skmlD}), $\ww$ is a   
$\ip1\nn$ set. 

Further, each cross-section $\ww_z$ coincides with 
the corresponding cross-section $\gfu_\xi$ 
of $\gfu$, where $\xi=\abs z$, thus $\ww_z$  
is a \dd\Eo class.

To prove \ref{lnv4},  
suppose to the contrary that $\ww$ is 
uniformized by a  ROD set $S\sq \ww$.
As $\omil=\omi$ holds, 
for every ordinal $\xi<\omi$ there is a real 
$z\in\dn\cap\rL$ satisfying $\abs z=\xi$. 
Let $z(\xi)$ be the \dd\lel least of such reals. 
Then
$$
R=\ens{\ang{\xi,x}\in \gfu}{\ang{z(\xi),x}\in S}  
$$
is a ROD subset of $\gfu$, 
which uniformizes the set $\gfu$, 
contrary to Lemma~\ref{lnu}.
Thus $\ww$ satisfies  
\ref{lnv1}, 
\ref{lnv3}, 
\ref{lnv4}.
\epf

\bpf[Theorem \ref{tuk}\ref{tuk1}]
Obvious by Lemma \ref{lnv}.
\epf

\parf{Auxiliary forcing relation} 
\las{fs}

Here we define a key instrumentarium for the  
proof of \ref{tuk2} of Theorem \ref{tuk}.
This is a forcing-type relation $\fo$. 
It is not directly connected with the forcing 
notion $\dP$, but rather related to the 
countable-support product $\pel^{\omi}$.  
But it happens to be compatible with the 
$\fP$-forcing relation for formulas of certain 
quantifier complexity 
(Lemma~\ref{ft.}). 
An important property of $\fo$ will be its 
permutation-invariance   
(Lemma~\ref{inv}), 
a property which the $\dP$-forcing relation  
definitely lacks.
This will be the key argument in the proof of   
Theorem~\ref{lun}.

{\ubf We argue in $\rL$.}

Let $\cL$ be a \rit{language\/} containing   
\index{formula!language $\cL$}%
variables $i,j,k,\dots$ of type 0 with the domain $\om$,  
and variables $x,y,z,\dots$ of type 1 with the domain $\dn$. 
Let \rit{terms} be variables of type  0 and    
expressions of the form $x(k)$. 
Atomic formulas are those of the form $R(t_1,\dots,t_n)$, 
where $R\sq\om^n$ is any \dd nary relation on   
$\om$ in $\rL$. 
Formulas containing no quantifiers over type 1 
variables, are \rit{arithmetic}.    
\index{formula!arithmetic}%
Formulas of the form
$$
\sus x_1\:\kaz x_2\:\sus x_3\:\dots\:\sus(\kaz)\,x_n\,\Psi
\quad\text{and}
\quad
\kaz x_1\:\sus x_2\:\kaz x_3\:\dots\:\kaz(\sus)\,x_n\,\Psi\,,
$$
where $\Psi$ is arithmetic, belong to types resp.\ 
\index{formula!$\lis1n$, $\lip1n$}%
\index{zzLS1n@$\lis1n$}%
\index{zzLP1n@$\lip1n$}%
$\lis1n$, $\lip1n$. 

Additionally, we allow codes   
$\rc\in\ccf$ to substitute free variables of type 1. 
We let $\abs\vpi=\bigcup_{\bc\in\vpi}\abs\bc$ 
for any \dd\cL formula, where $\bc\in\vpi$ means that  
a code $\bc$ occurs in $\vpi$. 
The semantics is as follows. 
Let $\vpi:=\vpi(\rc_1,\dots,\rc_k)$ be a  
\dd\cL formula, 
and all codes in $\ccf$ occurring in $\vpi$   
are explicitly indicated, 
and $\abs{\vpi}\sq B\sq\omi$. 
% --- не более чем 
%счетное множество, удовлетворяющее $\abs{\vpi}\sq B$, 
If $x\in\ddn B$  
\index{formula!$\vpi[x]$}%
\index{zzfix@$\vpi[x]$}%
then let $\vpi[x]$ denote the formula   
$\vpi(f^{\rc_1}(x\res\abs{f^{\rc_1}}),
\dots,f^{\rc_k}(x\res\abs{f^{\rc_k}}))$; 
all elements $f^{\rc_i}(x\res\abs{f^{\rc_i}})$ 
are reals in $\dn$. 
 
Arithmetic formulas and those in $\lis1n\cup\lip1n$, 
$n\ge1$, will be called \rit{normal}. 
\index{formula!normal}%
If $\vpi$ is a formula in $\lis1n$ or $\lip1n$ then       
$\vpi^-$ is the result of canonical transformation of 
\index{formula!$\vpi^-$}%
\index{zzfi@$\vpi^-$}%
$\neg\:\vpi$ to resp.\ $\lip1n$, $\lis1n$ form. 
We let $\vpi^-:=\neg\:\vpi$ for arithmetic formulas.

\bdf
[in $\rL$]
\lam{d:fo}
We define a relation $\zp\fo\vpi$ between 
\mut s $\zp\in\mtt$ 
and closed normal \dd\cL formulas:
\ben
\Renu
\itlb{fo1}%
if $\vpi$ is a closed \dd\cL formula, 
arithmetic or in $\lis11\cup\lip11$, 
and $\abs\vpi\sq B=\abs\zp$, 
then $\zp\fo\vpi$ whenever $\vpi[x]$ holds for all 
$x\in[\zp]$;

\itlb{fo2}%
if $\vpi:=\sus x\,\psi(x)$ is a closed    
$\lis1{n+1}$-formula, $n\ge 1$ 
($\psi$ belongs to $\lip1n$),
then $\zp\fo\vpi$ whenever there is a code $\bc\in\ccf$   
such that $\zp\fo\psi(\bc)$;

\itlb{fo3}%
if $\vpi$ is a closed $\lip1{n}$-formula, $n\ge 2$,
then $\zp\fo\vpi$ whenever there is no \mut\ $\zq\in\mtt$      
satisfying $\zq\leq\zp$ and $\zq\fo\vpi^-$.
\een
We define $\Fo\vpi=\ens{T\in\pes}{T\fo\vpi}$ and  
$\des\vpi=\Fo\vpi\cup\Fo{\vpi^-}$.
\edf

\ble
[in $\rL$]
\lam{deff}
If\/ $m\ge2$ and\/ $\vpi$ is a closed formula in\/ 
$\lis1m$, resp., $\lip1m$, then\/ 
$\Fo\vpi$ belongs to\/ $\is\hc{m-1}(\hc)$, 
resp., $\ip\hc{m-1}(\hc)$.
\ele
\bpf
If $\vpi$ is a $\lip11$ formula then $\Fo\vpi\in\fp11$ 
by Definition~\ref{d:fo}\ref{fo1}, and hence  
$\Fo\vpi$ belongs to $\id\hc{1}(\hc)$. 
Then argue by induction using  
\ref{d:fo}\ref{fo2},\ref{fo3}.
\epf

\parf{Auxiliary forcing relation: two lemmas} 
\las{2mm}

We here prove two key properties of the relation $\fo$. 
They will be used in the proof of 
Theorem \ref{tuk}\ref{tuk2} below.
One of them (Lemma \ref{ft.}) says that $\fo$ 
is connected with the truth in \dd\fP generic 
extensions similarly to the ordinary  
\dd\fP forcing --- for formulas of cartain complexity. 
The other one (Lemma \ref{inv}) claims the invariance  
$\fo$ relatively to the permutations of $\omi$.

Recall that a number $\nn\ge3$ is fixed by  
Definition~\ref{vdp}.

\ble
[in $\rL$]
\lam{dens}
Let\/ $\vpi$ be a closed normal\/ \dd\cL formula. 
Then the set\/ $\des\vpi$ is dense in\/ $\mtt$. 
If\/ $\vpi$ is of type\/ $\lis1m$, 
$m<\nn$, then\/ 
$\des\vpi\cap\fP$ is dense in\/ $\fP$.
\ele
\bpf
It suffices to prove the density of $\des\vpi$ 
for formulas $\vpi$ as in \ref{d:fo}\ref{fo1}. 
If $\vpi$ is such and $\zp\in\mtt$, 
$\abs\vpi\sq B=\abs\zp$, then the set 
$X(\vpi)=\ens{x\in [\zp]}{\vpi[x]}$ in  
$\ddn B$ belongs to $\fs11\cup\fp11$ and hence 
has the Baire property inside the closed    
set $[\zp]\sq\ddn B$. 
It remains to refer to Lemma~\ref{ssq}. 
The second claim follows by lemmas   
\ref{deff} and \ref{blo}.
\epf

\ble
\lam{ft.}
Assume that\/ $1\le n<\nn$, $\vpi\in\rL$ is a     
closed formula in\/ $\lip1{n}\cup\lis1{n+1}$, 
and a set\/ 
$\zG\sq\dP$ is generic over\/ $\rL$. 
Then\/ $\vpi[\xx \zG]$ holds in\/ $\rL[\zG]$,  
if and only iff\/ $\sus \zp\in \zG\,(\zp\fo\vpi)$.
\ele
\bpf
{\ubf Base of induction:} 
$\vpi$ is arithmetic or $\lis11\cup\lip11$,
as in \ref{d:fo}\ref{fo1}. 
If $\zp\in \zG$ and $\zp\fo\vpi$ then $\vpi[\xx \zG]$ 
holds in $\rL[\zG]$ by the Shoenfield absoluteness 
theorem, since  
$\xx \zG\res\abs\zp\in[\zp]$. 
In the opposite direction apply Lemma~\ref{dens}. 

{\ubf Step $\lip1n\imp\lis1{n+1}$:} 
$\vpi$ is $\sus x\,\psi(x)$, where $\psi$ 
belongs to $\lip1n$. 
Let $\zp\in \zG$ and $\zp\fo\vpi$. 
Then by Definition~\ref{d:fo}\ref{fo2}  
there exists a code $\bc\in\ccf\cap\rL$ such that  
$\zp\fo\psi(\bc)$.
By the inductive hypothesis, 
the formula $\psi(\bc)[\xx \zG]$, \ie,
$\psi[\xx\zG](f^\bc(\xx \zG\res B))$, 
where $B=\abs\bc$, holds in $\rL[\zG]$. 
Then $\vpi[\xx \zG]$ is true as well. 

Conversely assume that $\vpi[\xx \zG]$ holds in $\rL[\zG]$. 
There is a real $y\in\rL[\zG]\cap\dn$ such that 
$\psi[\xx\zG](y)$ holds.
By Lemma~\ref{repd},  
$y=\fg\bc=f^\bc(\xx \zG\res B)$, where $\bc\in\ccf\cap\rL$ 
and $B=\abs\bc$. 
Then $\psi(\bc)[\xx \zG]$ holds in $\rL[\zG]$. 
By the inductive hypothesis, 
there is $\zp\in \zG$ such that $\zp\fo\psi(\bc)$, 
hence, $\zp\fo\vpi$. 

{\ubf Step $\lis1n\imp\lip1{n}$:} 
$\vpi$ is a $\lip1n$ formula, $n\ge2$. 
Lemma~\ref{dens} yields a \mut\ $\zp\in \zG$ such   
that either $\zp\fo\vpi$ or $\zp\fo\vpi^-$.
If $\zp\fo\vpi^-$ then $\vpi^-[\xx \zG]$ holds  
in $\rL[\zG]$ by the inductive hypothesis, 
hence, $\vpi[\xx \zG]$ fails. 
Now assume that $\zp\fo\vpi$. 
We have to prove $\vpi[\xx \zG]$ in $\rL[\zG]$. 
Suppose to the contrary that $\vpi^-[\xx \zG]$ holds. 
By the inductive hypothesis, 
there exists a \mut\   
$\zq\in \zG$  such that $\zq\fo\vpi^-$.
But the \mut s $\zq,\zp$ belong to the generic set 
$\zG$, hence, they are compatible, which contradicts 
to the assumption $\zp\fo\vpi$.
\epf

{\ubf Invariance.} 
The relation $\fo$ turns out to be  
\rit{invariant} with respect to the natural action  
of the group $H$ of all self-inverse (\ie, $h=h\obr$) 
\index{invariance!permutation}%
\rit{permutations} of the set   
$\omil$ in $\rL$. 
Thus $h\in H$, iff $h\in\rL$,  
$h:\omil\onto\omil$ is a bijection, and $h=h\obr$.\vom

{\ubf We argue in $\rL$.}
Let $h\in H$. 
If $B\sq\omi$ and $F$ is a function defined on 
$B$ then a function $hF=h\ap F$ is defined on  
$\imb hB=\ens{h(\xi)}{\xi\in B}$ so that   
\index{action!$hF$}%
$(hF)(h(\xi))=F(\xi)$ for all $\xi\in B$. 
Thus $hF$ is equal to the superposition 
$F\circ h\obr$, and even $hF=F\circ h$ 
by the self-invertibility.

In particular if $x\in\ddn B$ then $hx\in\ddn{\imb hB}$, 
and if $\zp\in\mlt_B$ then $h\zp=h\ap\zp$ is a  
\mut\ in $\mlt_{\imb hB}$. 
Further, if $\bc\in\knf B$ then a code 
$h\bc=h\ap\bc\in\knf{\imb hB}$ can be  
canonically defined so that 
$f^{h\bc}(hx)=f^\bc(x)$ for all $\xi\in B$. 
Finally if $\vpi:=\vpi(\rc_1,\dots,\rc_k)$ is a   
\dd\cL formula then $h\vpi$ or $h\ap\vpi$ 
\index{action!$h\vpi$}%
denotes the formula $\vpi(h\rc_1,\dots,h\rc_k)$. 
Then  
$(h\vpi)[hx]$ coincides with $\vpi[x]$. 

\ble
[in $\rL$]
\lam{inv}
Let\/ $h\in H$, $\zp\in\mlt$, and\/ 
$\vpi$ is a closed normal\/ \dd\cL formula.
Then\/ $\zp\fo\vpi$ if and only if\/ $h\zp\fo h\vpi$.
\ele
\bpf
If $\vpi$ is a formula of type \ref{d:fo}\ref{fo1}  
then $[h\zp]=\ens{hx}{x\in[\zp]}$, 
and on the other hand, if $x\in[\zp]$ then $\vpi[x]$ 
coincides with $(h\vpi)[hx]$. 
We skip further routine indictive steps on the base 
of Definition~\ref{d:fo}\ref{fo2},\ref{fo3}.
\epf

\parf{Uniformization of sets with countable sections} 
\las{zav}

To prove claim \ref{tuk1} of Theorem~\ref{tuk} 
in the end of this section, 
we establish Theorem~\ref{lun} saying that in 
$\fP$-generic extensions any element of a countable 
$\fs1\nn$ set $X$ is constructible relative to  
the parameter of a $\fs1\nn$ definition of $X$. 
The relation $\fo$ and Lemma \ref{ft.} 
will play the key role.

\bte
%[proved in Section~\ref{zav}]
\lam{lun}
If a set\/ 
$\zG\sq\dP$ is\/ \dd\dP generic over\/ $\rL$ 
and\/ $p\in\rL[\zG]\cap\dn,$  
then it is true in\/ $\rL[\zG]$ any countable\/ 
$\is1\nn(p)$ set\/ $Y\sq\dn$ satisfies\/ 
$Y\sq\rL[p]$. 
\ete 

In fact a stronger claim $Y\in\rL[p]$ holds. 
However it requires more complex transformations 
beyond $H$, so we are going to skip this issue 
whatsoever. 

\bpf 
We argue in terms of Definition~\ref{vdp}. 
%Сначала мы докажем, что $Y\sq\rL[p]$. 
Suppose to the contrary that $Y\not\sq\rL[p]$. 
Then $Y=\ens{y\in\dn}{\vpi(p,y)}$, where   
$\vpi(p,y):=\sus z\,\psi(p,y,z)$ 
is a $\is1\nn$ formula with $p$ as the only parameter, 
and there is a real $y_0\in Y$, $y_0\nin\rL[p]$.
By Theorem \ref{fg}, there is an ordinal $\et<\omil$  
such that $p\in\rL[\xx \zG\res \cne\et]$ 
and $\bx_\et[\zG]\in\rL[y_0]$, 
and moreover, $\bx_\et[\zG]=g(y_0)$, 
where $g:\dn\to\dn$ is a continuous map coded in $\rL$.
By Lemma \ref{repd}, 
there exist codes $\bc,\bd\in\kng$ satisfying  
$p=\fg\bd=f^\bd(\xx \zG\res B)$ and  
$y_0=\fg\bc=f^\bc(\xx \zG\res B')$, 
where $B=\abs\bd\sq\cne\et$ and $B'=\abs\bc$. 
We can assume that $B\sq B'$ and $\et\in B'$. 
Note that definitely $\et\nin B$. 
%%%%%%%%%%%%%%%%%%%%%%%%%%
\vyk{
Существует \mutо $\zp\in \zG$, вынуждающее этот выбор, 
\te\ 
\ben
\fenu
\itlb{lun*}%
$\zp$ \dd\fP вынуждает, что 

1) 
множество $\ens{y\in\dn}{\vpi(f^\bd(\xx\uG\res B),y)}$ 
счетно, 

2)
выполнено $\vpi(f^\bd(\xx\uG\res B),f^\bc(\xx\uG\res B'))$, 
и 

3)
 $\bx_\et[\uG]=g_\jba(f^\bc(\xx\uG\res B'))$, 
\een
где, как обычно, $\uG$ --- имя для самого генерического 
фильтра $G$. 
По очевидным соображениям (увеличим $B'$, если нужно) 
можно предполагать, что $\abs\zp=B'$.
}%
%%%%%%%%%%%%%%%%%%%%%%%%%%

The goal is {\ubf to get a contradiction}. 

Consider the $\lis1\nn$ formula $\vpi(\bd,\bc)$. 
By the choice of the codes, $\vpi(\bd,\bc)[\xx \zG]$ 
coincides with  
%$\vpi(f^\bd(\xx \zG\res B),f^\bc(\xx \zG\res B'))$, 
$\vpi(\fg\bd,\fg\bc)$, 
therefore $\vpi(\bd,\bc)[\xx \zG]$ 
holds in $\rL[\zG]$. 
By Lemma \ref{ft.}, there is a \mut\ $\zq\in \zG$ 
satisfying $\zq\fo\vpi(\bd,\bc)$. 

{\sloppy 
Further, the equality $\bx_\et[\zG]=g(y_0)$ 
(see above) can be rewritten as\break   
$f^\je(\xx\uG\res B')=g(f^\bc(\xx \zG\res B'))$, 
where $\je\in\knf{B'}\cap\rL$ is a canonical code of 
the map $f^\je(x)=x(\et)$. 
We render this formula as 
}
$$
\sus z\,
\big(
z= f^\bc(\xx \zG\res B')\land 
f^\je(\xx\uG\res B')=g(z)
\big).
%\sus z\,\big(z=\fg\bc\land\fg\je=g(z)\big).
$$ 
%чтобы избежать явной подстановки одной кодированной 
%функции в другую.
As above by Lemma \ref{ft.} there exists a 
\mut\ $\zq'\in \zG$ satisfying
$\zq'\fo\sus z\,(z=\bc\land\je=g(z))$.  
We can assume that $\zq'=\zq$. 
(Otherwise replace both \mut s by a stronger \mut\ 
in $\zG$).
Thus we have  
\ben
\fenu
%\atc
\itlb{lun**}%
$\zq\fo\vpi(\bd,\bc)$ \ and \ 
$\zq\fo\sus z\,(z=\bc\land\je=g(z))$. 
\een
We can wlog assume that $\abs \zq= B'$, 
as otherwise we just replace  
$B'$ by $B'\cup \abs \zq$ and $\zq$ by 
$\zq\ares{(B'\cup \abs \zq)}$.

If $\vt<\omil$ then let $H_\vt$ denote the set 
of all permutations $h\in H$ such that  
$h(\xi)=\xi$ for all $\xi\in B$ and $h(\xi)>\vt$ 
for all $\xi\in B'\bez B$.

\ble
%[в $\rL$]
\lam{lun1}
If\/ $\vt<\omil$ then there is a permutation\/ 
$h\in H_\vt$ and a \mut\/ $\zq'\in \zG$ such that\/  
$\zq'\leq h\ap\zq$. \ 
{\rm (It is not assumed that $h\ap\zq\in\fP$.)}
\ele
\bpf[Lemma]
Arguing in $\rL$, consider the set $\zD_\vt$ of  
all \mut s $\zq'\in\mtt$ such that  
%либо $\zq$ несовместно с $\zp'$ в $\mtt$, либо же 
$\zq'\leq\zq$ and there exists a permutation  
$h\in H_\vt$ such that the \mut\ $h\ap\zq$ 
satisfies $\zq'\leq h\ap\zq$.
A routine estimation shows that $\zD$ is a  
$\is\hc1(\zq,\vt)$ set. 
Therefore by Lemma \ref{blo} there is a \mut\ 
$\zq'\in \zG$, 
such that  either  
(1) 
$\zq'\in\zD_\vt$, or  
(2)
there is no \mut\ $\zr\in\zD_\vt$ satisfying  
$\zr\leq\zq'$. 
And as $\zq$ also belongs to $\zG$, 
we cal wlog assume that  
$\zq'\leq\zq$.

We claim that (2) is impossible. 
Indeed let $\ga<\omil$   
satisfy $\abs{\zq'}\sq\ga$ and $\ga\ge\vt$. 
Define a permutation $h$ by $h(\xi)=\xi$ 
for $\xi\in B$, $h(\xi)=h\obr(\xi)=\ga+\xi$ for 
$\xi<\ga$, $\xi\nin B$, and still  
$h(\xi)=\xi$ for all other $\xi<\omil$. 
The \mut s $\zq'$ and $\ju=h\ap\zq'$ coincide on 
the common domain $\abs{\zq'}\cap\abs{\ju}=B$, 
hence are compatible. 
It follows that the union $\zr=\zq'\cup\ju$ 
belongs to $\mtt$ and $\zr\leq\zq',\ju$. 
And further we have $\zr\leq\ju=h\ap\zq'\leq h\ap\zq$ 
by construction, hence $\zr\in\zD$, as required. 
Thus (2) fails. 
Therefore (1) holds, that is,
$\zq'\in\zD_\vt$, as required. 
\epF{Lemma}

Coming back to Theorem \ref{lun}, recall that  
$\omil$ remains a cardinal in $\fP$-generic  
extensions by Lemma \ref{wcc}. 
Therefore Lemma \ref{lun1} allows to define 
by induction an increasing sequence  
$\sis{\vt_\nu}{\nu<\omil}$ of ordinals $\vt_\nu<\omil$  
and a sequence of \mut s  $\zq_\nu\in \zG$ 
and a sequence of permutations $h_\nu\in H_{\vt_\nu}$, 
satisfying  
$B'\sq\vt_0$ and $\zq_\nu\leq h_\nu\ap \zq$ for  
all $\nu$, and $\abs{\zq_\mu}\sq\vt_\nu$ for 
$\mu<\nu$.

Let  
$\zp_\nu=h_\nu\ap\zq$, 
$\bc_\nu=h_\nu\ap\bc$, 
$\bd_\nu=h_\nu\ap\bd$, 
$\je_\nu=h_\nu\ap\je$ 
for all $\nu$. 
Then we have   
$\zp_\nu\fo\vpi(\bd_\nu,\bc_\nu)$ 
and   
$\zp_\nu\fo\sus z\,(z=\bc_\nu\land \je_\nu=g(z))$
by \ref{lun**} and lemma \ref{inv}. 
It follows that   
\ben
\fenu
\atc
%\atc
\itlb{lun*}
$\zq_\nu\fo\vpi(\bd,\bc_\nu)$ 
\ and \ 
$\zq_\nu\fo\sus z\,(z=\bc_\nu\land \je_\nu=g(z))$,
\een
since $\zq_\nu\leq \zp_\nu$, and, with respect  
to the code $\bd$: $\bd_\nu=h_\nu\ap\bd=\bd$. 
(Indeed $h_\nu(\xi)=\xi$ whenever $\xi\in B=\abs\bd$.) 

Recall that $f^\bd(\xx \zG\res B)= p$. 
Let $B'_\nu=\imb h{B'}$, 
$z_\nu=f^{\bc_\nu}(\xx \zG\res B'_\nu)$, 
and $q_\nu=f^{\je_\nu}(\xx \zG\res B'_\nu)$. 
If $\nu<\omil$ then, by \ref{lun*} and 
Lemma \ref{ft.},  
$\vpi(p,z_\nu)$ holds in $\rL[\zG]$ --- 
hence $z_\nu\in Y$, and we have $q_\nu=g(z_\nu)$ as well. 
Further,
$$
\bay{lcccccccccccc}
q_\nu
&=&f^{\je_\nu}(\xx \zG\res B'_\nu)
= (h_\nu\ap f^{\je})(\xx \zG\res B'_\nu)
= f^{\je}(h_\nu\obr(\xx \zG\res B'_\nu))&=\\[1ex]
&=& f^{\je}((h_\nu\obr(\xx \zG)\res B')
= (h_\nu\obr(\xx \zG)) (\et) 
=(\xx \zG)(\et_\nu)
&=&\bx_{\et_\nu}[\zG]\,,
\eay
$$ 
where $\et_\nu=h_\nu(\et)$. 
Thus an uncountable sequence of the reals $z_\nu\in Y$ 
in $\rL[\zG]$ ($\nu<\omil$) is defined, 
and it satisfies  
$g(z_\nu)=\bx_{\et_\nu}[\zG]$, $\kaz \nu$.
The ordinals $\et_\nu=h_\nu(\et)$ satisfy  
$\et_\nu\ge\vt_\nu$ by the choice of $h_\nu$, 
since $\et\in B'\bez B$. 
Therefore there exist uncountably many pairwise 
different of $\et_\nu$ in $\rL[\zG]$. 
It follows that there exist uncountably many pairwise 
different generic reals $\bx_{\et_\nu}[\zG]$.
On the other hand, all reals $z_\nu$ belong  
to the countable set $Y$, 
and $\bx_{\et_\nu}[\zG]=g(z_\nu)$, 
where $g$ does not depend on $\nu$.  
This is a {\ubf contradiction} required, 
and the theorem is proved.
\epf

\vyk{ 

To prove $Y\in\rL$, a stronger statement, it suffices 
now to show that if $y_0\in\dn\cap\rL$ then 
$y_0\in Y$ iff $\sus T\in\pes\,(T\fo\vpi(\rc_0))$, where 
$\rc_0\in\ccf\cap\rL$ is the code of the constant 
function $f_{\rc_0}(x)=y_0$, $\kaz x\in\dn$.

If $y_0\in Y$ then the formula $\vpi(y_0)$, equal 
to $\vpi(\rc_0)[x]$ for any $x$, is true in $\rL[\zG]$ 
by the choice of $\vpi$. 
It follows by Lemma~\ref{ft.} that there is a tree 
$T\in G$ satisfying $T\fo\vpi(\rc_0)$, as required. 

Now suppose that $T\in\pes$ 
(not necessarily $\in\dP$!) 
and $T\fo\vpi(\rc_0)$. 
As the set 
$D=\ens{T\in\pes}{\oi T\text{ is co-infinite}}$ 
(see Lemma~\ref{jden}\ref{prop5})
is open dense in $\pes$, 
we can assume that $\oi T$ is co-infinite. 
On the other hand, it follows from 
Lemma~\ref{jden}\ref{prop5} that there is a 
tree $S\in G\cap D$, 
so that $\oi S$ is co-infinite as well. 
Now we have $S\fo\vpi(\rc_0)$ by Corollary~\ref{tt'}, 
and then $\vpi(\rc_0)[\xx G]$ is true in $\rL[G]$ 
by Lemma~\ref{ft.}, that is, $\vpi(y_0)$ holds 
in $\rL[G]$, and $y_0\in Y$, as required.  
}    

\bpf[\ubf Theorem~\ref{tuk}\ref{tuk2}]
Arguing under the requirements of Therorem \ref{tuk}, 
assume that, in $\rL[\zG]$, $p\in\dn$ and $W\sq\dn\ti\dn$  
is a $\is1\nn(p)$ set whose cross-sections  
$W_x=\ens{y}{\ang{x,y}\in W}$ 
are at most countable. 
Every set $W_x$ is  
$\is1\nn(p,x)$, whence $W_x\sq\rL[p,x]$ 
by Theorem~\ref{lun}.
If $W_x\ne \pu$ then let $q_x$ be the  
\dd{<_{px}}least real in $W_x$, 
where $<_{px}$ is the canonical G\"odel 
well-ordering of $\rL[p,x]$. 
The set  
$Q=\ens{\ang{x,q_x}}{x\in\dn\land W_x\ne\pu}$ 
then uniformizas $W$. 
Moreover
$$
\ang{x,y}\in Q
\leqv
\ang{x,y}\in W
\land
\kaz z\,(z<_{px}y\imp \ang{x,y}\nin W)\,.
$$
Therefore the set $Q$ belongs to $\id1{\nn+1}(p)$, 
or more exactly is the intersection of a $\is1\nn(p)$ 
set and a $\ip1\nn(p)$ set, because the G\"odel 
well-orderings $<_{px}$ are well-known to be   
$\is12(p,x)$-definable uniformly in $p,x$.
\epf

\qeD{Theorems \ref{tun} and \ref{Tun}}

%{\small\printindex}

\bibliographystyle{plain}
{\small
\renek{\refname} {{\large\bf References}}
\addcontentsline{toc}{subsection}{\indent\quad References}

\bibliography{38e,kle}
}

\def\indexname{Index{\normalsize\ubf\ \ 
\addcontentsline{toc}{subsection}{\indent\quad Index 
%(added for the convenience of the refereeing process)
}%
%(added for the convenience  of the refereeing process)
}}
%\footnotetext{xxx}
%\small\printindex

%
\input{38e.ind}

% added by arXiv admin:
\typeout{get arXiv to do 4 passes: Label(s) may have changed. Rerun}

\end{document}